\documentclass{amsart}
\usepackage{amssymb,euscript,amsmath, mathrsfs, amscd}
\usepackage[pdftex]{graphicx}
\usepackage[pdftex]{color}

\newcounter{ENUM}
\newcommand{\itm}{\item}
\newenvironment{ilist}[1][0]{\renewcommand{\theENUM}{\roman{ENUM}}\renewcommand{\itm}{\addtocounter{ENUM}{1}\item[(\theENUM)]}\begin{itemize}\setcounter{ENUM}{#1}}{\end{itemize}}
\newenvironment{Ilist}{\renewcommand{\theENUM}{\Roman{ENUM}}\renewcommand{\itm}{\addtocounter{ENUM}{1}\item[(\theENUM)]}\begin{itemize}\setcounter{ENUM}{0}}{\end{itemize}}

\newcommand{\margh}[1]{}

\def\risom{\overset{\sim}{\rightarrow}}

\input xy
\xyoption{all}
\CompileMatrices

\newcommand{\el}{$\ell$}
\newcommand{\e}{\varepsilon}
\def\ZZ{{\mathbb Z}}

\def\PP{{\mathbb P}}

\def\cG{{\mathcal G}}

\def\cO{{\mathcal O}}

\def\cP{{\mathcal P}}

\def\cU{{\mathcal U}}

\def\cM{{\mathcal M}}

\def\cPic{{\mathcal P}ic}
\def\sE{{\mathscr E}}
\def\sF{{\mathscr F}}
\def\sG{{\mathscr G}}

\def\sL{{\mathscr L}}
\def\sM{{\mathscr M}}
\def\sO{{\mathscr O}}
\def\sV{{\mathscr V}}
\def\sW{{\mathscr W}}

\def\fg{{\mathfrak g}}

\def\vp{\varphi}

\def\Sch{\operatorname{Sch}}

\def\Hom{\operatorname{Hom}}

\def\GL{\operatorname{GL}}

\def\Spec{\operatorname{Spec}}

\def\id{\operatorname{id}}

\def\rk{\operatorname{rk}}

\def\ord{\operatorname{ord}}

\def\EHT{\operatorname{EHT}}

\def\sw{\operatorname{sw}}
\def\nd{\operatorname{nd}}
\def\LG{\operatorname{LG}}
\def\LAG{\operatorname{LAG}}
\def\II{\operatorname{II}}
\newcommand{\eind}{\operatorname{EI}}
\newcommand{\oind}{\operatorname{OI}}

\newtheorem{thm}{Theorem}[section]
\newtheorem{prop}[thm]{Proposition}
\newtheorem{lem}[thm]{Lemma}
\newtheorem{cor}[thm]{Corollary}

\theoremstyle{definition}
\newtheorem{defn}[thm]{Definition}

\newtheorem{subex}{Example}
\numberwithin{subex}{section}
\newtheorem{sit}[thm]{Situation}
\theoremstyle{remark}
\newtheorem{notn}[thm]{Notation}
\newtheorem{rem}[thm]{Remark}

\numberwithin{equation}{section}

\begin{document}
\title[Linked symplectic forms and limit linear series]{Linked symplectic forms and limit linear series in rank $2$ with
special determinant}
\author{Brian Osserman}
\author{Montserrat Teixidor i Bigas}
\begin{abstract} 
We generalize the prior linked symplectic Grassmannian construction, applying
it to to prove smoothing results for rank-$2$ limit linear series with 
fixed special determinant on chains of curves. We apply this general 
machinery to prove new results on nonemptiness and dimension of rank-$2$ 
Brill-Noether loci in a range of degrees.
\end{abstract}

\thanks{The first author was partially supported by NSA grant H98230-11-1-0159
and Simons Foundation grant \#245939 during the preparation of this work.}

\maketitle

\tableofcontents

\section{Introduction}

At its most basic, higher-rank Brill-Noether theory addresses the
question: for a general curve of genus $g$, how many global sections
can a (semi)stable vector bundle of given rank and degree have? This
has been an active subject of study for more than 20 years, and the picture
which has emerged is one of complexity, with no comprehensive conjectures
even in the case of rank $2$. Nonetheless, those cases which have been 
understood have already been important in a number of strikingly different 
contexts, from Mukai's work \cite{mu6} on classification of Fano threefolds 
to recent work of Bhosle, Brambila-Paz and Newstead \cite{b-b-n2} on 
Butler's conjecture. See \cite{g-t1} and the introduction of \cite{os20} 
for a more detailed survey.

In the present paper, we consider $\fg^k_{2,d}$s consisting
of pairs $(\sE,V)$, where $\sE$ is a vector bundle of rank $2$ and degree
$d$, and $V$ is a $k$-dimensional space of global sections, and we extend
the range where stable $\fg^k_{2,d}$ in a range of degrees and genera. 
Our main tool for doing this is a new smoothing theorem in the case of
special determinant, and as in the classical rank-$1$ case treated by
Eisenbud and Harris, our smoothing theorem is inseparable from a careful
study of dimension. Now, the classical expected dimension $\rho$ 
generalizes to higher rank, and in some cases gives the correct dimension
for the moduli space of $\fg^k_{r,d}$s,
but is not enough to understand the general picture. One case where the
picture seems clearer is that of rank-$2$ vector bundles with fixed
canonical determinant: Bertram, Feinberg \cite{b-f2} and Mukai \cite{mu2}
observed that in this case, the behavior appears to be closer to 
that of the classical rank-$1$ case, albeit with a modified expected 
dimension. While the naive expected dimension for the fixed determinant 
case is $\rho-g$, they showed that symmetries in the canonical determinant
case forced the dimension to be at least 
$\rho_{\omega}:=\rho-g+\binom{k}{2}$, and they conjectured that this is
in fact the correct dimension, and in particular that the relevant moduli
spaces should be nonempty when $\rho_{\omega} \geq 0$. The existence 
portion of their conjecture remains open, while their results on
modified expected dimension were generalized by the first author to a wider 
class of special determinants in \cite{os16}.

In \cite{te1}, the second author developed a theory of limit linear series
in higher rank, and in \cite{te5} applied it to show that, subject to 
proving a certain smoothing theorem, one could use it to prove existence 
for a large infinite family of cases of the Bertram-Feinberg-Mukai 
conjecture. The main obstruction to proving the smoothing theorem is that
the symmetries which occur in the canonical determinant case and lead to
the larger expected dimension do not persist in an obvious way in the
definition of higher-rank limit linear series. In \cite{o-t1}, we showed
that by relating higher-rank limit linear series to an alternate
limit linear series construction first introduced in \cite{os8}, it is
possible to prove the necessary smoothing theorem for curves having two
irreducible components. However, because higher-rank limit linear series
do not have a simple inductive structure, this is not enough to prove the
desired existence results. The purpose of the present paper, in 
conjunction with \cite{os20}, is to generalize the constructions of
\cite{o-t1} to work on chains of curves with any number of components,
leading to the proof of the desired smoothing theorem. In addition, we
broaden the families of examples considered in \cite{te5}
to more general special determinants, consequently proving
existence of large families of components of moduli spaces of $\fg^k_{2,d}$s 
with fixed special determinant and having the (modified) expected dimension.

We now explain our smoothing theorem for higher-rank limit linear series.
Given nonnegative integers $g,d,k$, the naive expected dimension for the 
stack of $\fg^{k}_{2,d}$s on a general curve of genus $g$ is equal to
$$\rho-1=4g-4-k(k-d+2g-2).$$
Thus, given a line bundle of degree $d$ on such a curve $X$, the naive expected
dimension for the stack of $\fg^{k}_{2,d}$s on $X$ with fixed determinant
$\sL$ is given by $\rho-g$, but according to \cite{os16} the higher expected 
dimension $\rho-g+\binom{k}{2}$ applies when $\sL$ is special.
Obviously, this distinction is only
relevant if $k \geq 2$, so we will impose this assumption. A typical 
smoothing theorem for limit linear series says that if one has a family
of limit linear series which occurs in the expected dimension on a given
reducible curve, then it smooths out to give linear series on nearby 
smooth curves. For our smoothing theorem, we restrict to the case that
the special fiber is a chain of curves, as has been the case for every
family considered in \cite{te1} and \cite{te5}, and we work with a
certain open subset of ``chain-adaptable'' (higher-rank) limit linear 
series, introduced in \cite{os20}. In order to treat (semi)stability
conditions more easily, we work with the notion of \el-(semi)stability
introduced in \cite{os22}. Compared to the usual notion of stability
on a reducible curve, this is weaker and more canonical, and leads to
equally strong conclusions in degeneration arguments. Our main smoothing
theorem is then as follows.

\begin{thm}\label{thm:main} 
Given $g,d,k$, with $k\geq 2$, suppose that there exists a projective nodal 
curve $X_0$ with dual graph $\Gamma$ a chain, a
special line bundle $\sL_0$ of degree $d$ on $X_0$, and 
$((\sE_v,V_v)_{v \in V(\Gamma)},(\vp_e)_{e \in E(\Gamma)},\psi)$ a 
chain-adaptable limit linear series of rank $2$ and fixed determinant
$\sL_0$, such that the space of such limit linear series on $X_0$ has
the expected dimension $\rho-g+\binom{k}{2}$ at the corresponding point.

Then for a general smooth curve $X$ of genus $g$ and a general special line 
bundle $\sL$ of degree $d$, the space of $\fg^{k}_{2,d}$s with fixed 
determinant $\sL$ on $X$ is nonempty, with a component of expected dimension
$\rho-g+\binom{k}{2}$.

Furthermore, if 
$((\sE_v,V_v)_{v \in V(\Gamma)},(\vp_e)_{e \in E(\Gamma)},\psi)$ is
\el-semistable (respectively, \el-stable), then we get the same statement
for semistable (respectively, stable) $\fg^{k}_{2,d}$s on $X$. 
\end{thm}

The relevent limit linear series terminology 
is reviewed in Section \ref{sec:lls} below.

Thereom \ref{thm:main} renders unconditional the existence results in 
\cite{te5} for the case of canonical determinant. Moreover, using similar
techniques, we produce more general families of limit linear series with
special determinant, and we conclude more general existence results, 
as follows:

\begin{thm}\label{thm:main-2} 
Given $g,d,k$ nonnegative, with $k\geq 2$ and $g-2 \leq d \leq 2g-2$,
suppose that 
$$4g \geq \begin{cases} 
k^2+2k(2g-2-d): & d \text{ even}\\
k^2+4+2k(2g-2-d): & d \text{ odd}.
\end{cases}$$
Then for a general smooth curve $X$ of genus $g$ and a general special line 
bundle $\sL$ of degree $d$, the space of semistable $\fg^{k}_{2,d}$s with fixed 
determinant $\sL$ on $X$ is nonempty, with a component of expected dimension
$\rho_{\sL}:=\rho-g+\binom{k}{2}$. If further
$$(g,d,k) \neq (1,0,2), (2,2,2), (3,2,2) \text{ or } (4,6,4),$$
the same is true of the space of stable $\fg^{k}_{2,d}$s.
\end{thm}

Note that we do not have to assume that $\rho_{\sL} \geq 0$, and in fact
this is not always the case; see Remark \ref{rem:neg-dim}. Also, the
case $(g,d,k)=(2,2,2)$ is an exception to the existence of stable bundles in
the Bertram-Feinberg-Mukai 
conjecture: although $\rho_{\sL}=0$ in this case, there cannot be any 
stable $\fg^{2}_{2,2}$. See the note on page 123 of \cite{te4}. In particular,
within the imposed range for $g$, our stability results are optimal.

This existence result can be seen as validation of the first main
result of \cite{os16}, insofar as it provides many examples of components
of moduli spaces of $\fg^{k}_{2,d}$s with fixed special determinant having
dimension equal to the modified expected dimension $\rho-g+\binom{k}{2}$.
However, Theorem \ref{thm:main} sets up a more general machinery, and
indeed Zhang \cite{zh2} has already been able to use it to prove existence 
results for the canonical determinant case which improve on those of 
Theorem \ref{thm:main-2}. In this context, we also mention recent work of
Lange, Newstead and Park \cite{l-n-p1} which approaches the existence
question in the canonical determinant case via fundamental class
computations. Their results are comparable to those of Zhang, except that
due to combinatorial complications, they restrict to prime genera.

Finally, by allowing determinants to vary we immediately conclude 
existence of components of larger than the expected dimension $\rho-1$ in 
a wide range of cases.

\begin{cor}\label{cor:varying-det}
Given $g,d,k$ nonnegative, with $k\geq 2$ and $g-2 \leq d \leq 2g-2$,
suppose that 
$$4g \geq \begin{cases} 
k^2+2k(2g-2-d): & d \text{ even}\\
k^2+4+2k(2g-2-d): & d \text{ odd}.
\end{cases}$$
Then for a general smooth curve $X$ of genus $g$, the space of semistable 
$\fg^{k}_{2,d}$s is nonempty, with a component of dimension
$\rho+\binom{k}{2}-(d-g+3)$. If further
$$(g,d,k) \neq (1,0,2), (2,2,2), (3,2,2) \text{ or } (4,6,4),$$
the same is true of the space of stable $\fg^{k}_{2,d}$s.

In particular, if $\binom{k}{2}>d-g+2$, the space has a component of 
dimension strictly greater than $\rho-1$.
\end{cor}

The main theoretical tool in the proof of Theorem \ref{thm:main} involves 
the definition and study of a
generalization of the linked alternating Grassmannian introduced in 
\cite{o-t1}, together with an associated nondegeneracy condition which
replaces the notion of linked symplectic forms. When combined with the
general limit linear series theory developed in \cite{os20}, we are able
to deduce our general smoothing results. Rather than considering a
strict generalization of the linked symplectic Grassmannian of \cite{o-t1},
we consider a somewhat broader notion; this simplifies definitions and
arguments, without weakening the resulting smoothing theorems. In brief,
prelinked alternating Grassmannians are cut out inside prelinked 
Grassmannians by an isotropy condition with respect to a ``linked 
alternating form.'' Prelinked Grassmannians have an open subset consisting 
of ``simple points,'' and we can analyze the behavior of prelinked 
alternating Grassmannians on this subset. We also introduce a smaller 
subset of ``internally simple points,'' and if the prelinked alternating 
form in question is ``internally symplectic,'' we get especially good 
behavior, which leads to our 
smoothing theorem. The relevant statement is as follows.

\begin{thm}\label{thm:main-lag} If
$\LAG(r,\sE_{\bullet},\left<,\right>_{\bullet})
\subseteq \LG(r,\sE_{\bullet})$ is a prelinked alternating Grassmannian,
and $z \in \LAG(r,\sE_{\bullet},\left<,\right>_{\bullet})$ is a simple
point of $\LG(r,\sE_{\bullet})$, then locally at $z$, we have that
$\LAG(r,\sE_{\bullet},\left<,\right>_{\bullet})$ is cut out by $\binom{r}{2}$
equations inside $\LG(r,\sE_{\bullet})$.

If further $\left<,\right>$ is internally symplectic, and $z$ is an
internally simple point, then $\LAG(r,\sE_{\bullet},\left<,\right>_{\bullet})$ 
is smooth at $z$ of codimension $\binom{r}{2}$ inside $\LG(r,\sE_{\bullet})$.
\end{thm}

In generalizing the ideas of \cite{o-t1} to the present situation, the 
main challenge is to determine the suitable definitions, which are not
obvious from the previously developed two-component case. Given the 
correct definitions, the arguments of \cite{o-t1} go through in a 
transparent manner to obtain Theorem \ref{thm:main-lag}. For further
motivation of our definitions, see \cite{o-t1}, especially Remarks 3.8 and 
4.6.

\section{Preliminaries}

We begin with some definitions of a combinatorial nature. Ultimately,
$\Gamma$ will be the dual graph of a reducible curve, and the associated
graph $G$ will be used to keep track of multidegrees of vector bundles
and natural maps between them.

\begin{sit}\label{sit:basic} Let $\Gamma$ be a tree.
Let $d,k$ be positive integers, and fix also integers $b$ and
$d_v$ for each $v \in V(\Gamma)$, satisfying
$$\sum_{v \in V(\Gamma)} d_v - |E(\Gamma)|2b=d.$$

We define a directed graph $G$ as follows:
let $V(G) \subseteq \ZZ^{V(\Gamma)}$ consist of vectors
$w=(i_v)_{v \in V(\Gamma)}$ satisfying:
\begin{ilist}
\itm $\sum_v i_v=d$, 
\itm $i_v \equiv d_v \pmod{2}$ for all $v \in V(\Gamma)$,
\itm for every subtree $\Gamma'$ of $\Gamma$ obtained as a connected
component of the complement of some edge of $\Gamma$, 
$$|V(\Gamma')|2b-\sum_{v \in \Gamma'} (d_v-i_v) \geq 0.$$
\end{ilist}
Then, $G$ has an edge from $w$ to $w'$ in $V(G)$ 
if there is a vertex $v \in V(\Gamma)$ of valence $\ell$
such that $w'-w$ is $-2\ell $ in index $v$, is $2$ in index $v'$ for
each $v'$ adjacent to $v$, and is $0$ elsewhere. Given an edge $\e \in E(G)$,
let $v(\e)$ be the associated vertex of $V(\Gamma)$.

We also define the hyperplane $H_d \subseteq \ZZ^{V(\Gamma)}$ to be the
set of vectors satisfying (i) above.
\end{sit}

Note that given $w,w' \in V(G)$, we have $\frac{w+w'}{2} \in H_d$.

Without further comment, we will assume that we are in the above situation,
and refer freely to the above notation. 

\begin{figure}
\centering
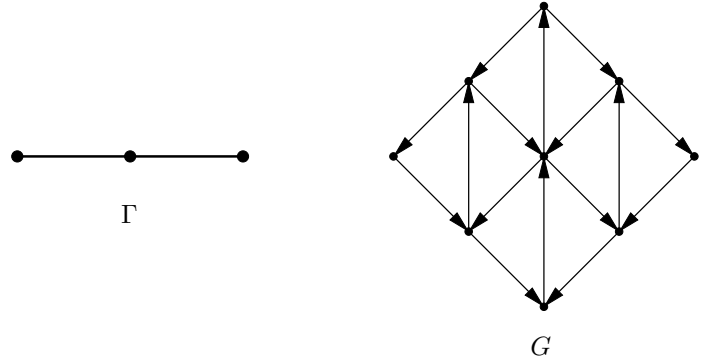
\caption{An example of Situation \ref{sit:basic}, with $b=2$.}
\end{figure}

It is convenient to introduce a distance function on $H_d$ as follows:

\begin{defn}\label{def:distance}
Given $w,w' \in H_d$, define $d(w,w')$ to be the minimal
number of operations required to get from $w$ to $w'$, where an allowable
operation is, for a choice of $v \in V(\Gamma)$ of valence $\ell$, to 
substract $\ell$ from the $v$th coordinate and add $1$ to the $v'$th 
coordinate for each vertex $v'$ adjacent to $v$.
\end{defn}

\begin{rem}
Situation \ref{sit:basic} is described more generally in Definitions
3.1.2 and 3.4.9 of \cite{os20} (where $G$ is called $\bar{G}_{\II}$),
with an additional parameter $r$ which in our case is set to
$2$. From this point of view, $d(w,w')$ is in essence the shortest path 
from $w$ to $w'$ in $G$, considered for the case $r=1$, except without the 
imposition of condition (iii) for 
the vertices of $G$. Note however that due to the directed nature of $G$,
this distance function is not symmetric.
\end{rem}

We now move on to the fundamental definitions.
We work throughout over a fixed base scheme $S$. We will use $h$ and $t$
to denote the head and tail of an edge of a graph, or a path in a graph.

\begin{defn}\label{defn:basic} Suppose we are given
data $\sE_{\bullet}$ consisting of
vector bundles $\sE_w$ of rank $2$ on $S$ for each $w \in V(G)$, and
morphisms $f_{\e}:\sE_{t(\e)} \to \sE_{h(\e)}$
for each $\e \in E(G)$. Given a directed path $P=(\e_1,\dots,\e_n)$ in $G$,
denote by $f_P$ the composition $f_{\e_n} \circ \cdots \circ f_{\e_1}$,
and by $v(P)$ the set (considered with multiplicities) 
$\{v(\e_1),\dots,v(\e_n)\}$. 

Given also $s \in \Gamma(S,\sO_S)$, 
we say that $\sE_{\bullet}$ is \textbf{$s$-prelinked} if
for $P$ and $P'$ any directed paths in $G$ with the same tail and
head and with the length of $P$ less than or equal to the length of $P'$, we
have 
$$f_{P'} = s^{c} f_{P},$$
where $c$ is the difference between the number of times $V(\Gamma)$ appears
in $v(P')$ and the number of times $V(\Gamma)$ appears in $v(P)$.
\end{defn}

This definition is slightly more restrictive than that considered in
Appendix A of \cite{os20}, reflecting that we consider here a more specific
collection of graphs $G$ onto which we may impose more structure.

Note that for two given vertices of $G$, any two minimal paths $P,P'$ have
$v(P)=v(P')$, and more generally, any path $P$ has the same 
$v(P)$ as a path obtained as a minimal path together with insertions
of copies of $V(\Gamma)$ -- see Proposition 3.1.4 of \cite{os20}.
Thus, the condition implies in
particular that if $P$ is a minimal path between two vertices, the resulting 
map $f_P$ does not depend on $P$, and more generally it
specifies precisely how $f_P$ depends on $P$. Thus, the following
notation makes sense:

\begin{notn}\label{notn:fww}
If $\sE_{\bullet}$ is $s$-prelinked, given $w,w'\in V(G)$,
we denote by $f_{w,w'}:\sE_w \to \sE_{w'}$ the map $f_P$ determined by
any minimal path $P$ from $w$ to $w'$.
\end{notn}

\begin{defn}\label{defn:simple-pt} Let $K$ be a field over $S$, and 
for $s \in K$, let $E_{\bullet}$ be $s$-prelinked on $\Spec K$. 
We say that $E_{\bullet}$ is \textbf{simple} if
there exist $w_1,\dots,w_r \in V(G)$ (not necessarily distinct) and
$v_i \in E_{w_i}$ for $i=1,\dots,r$ such that for every $w \in V(G)$,
there exist paths $P_1,\dots,P_g$ with each $P_i$ going from $w_i$ to
$w$, and such that $f_{P_1}(v_1),\dots,f_{P_r}(v_r)$ form a basis for
$E_w$.
\end{defn}

The following is a straightforward consequence of Nakayama's lemma.

\begin{prop}\label{prop:structure} Suppose that $\sE_{\bullet}$
is simple at a point $x \in S$. Then locally at $x$, for $w \in G$ there 
exist $r_w \geq 0$ and subbundles $\sW_w \subseteq \sE_w$ of rank $r_w$ 
such that:
\begin{ilist}
\itm $\sum_{w \in G} r_w = r$.
\itm The natural map
$$\bigoplus_{w' \in G} f_{w',w}(\sW_{w'}) \to \sE_w$$
is an isomorphism for each $w \in G$.
\end{ilist}

In particular, the points of $S$ at which $\sE_{\bullet}$ is simple form
an open subset.
\end{prop}

\section{Linked bilinear forms}

In this section and the next, we investigate the new definitions which
form the basis for our foundational results. In some sense, the 
definitions are the most important part, as they are calibrated so
that the arguments already used in \cite{o-t1} will go through in the
more general setting. See Remarks 3.8 and 4.6 and Examples 5.3 and 5.4
of \cite{o-t1} for discussion of the motivation behind the definitions.

As we discuss linked bilinear forms, we restrict from now on to the case
$r=2$. For arbitrary $r$, one can work instead with multilinear forms,
but nondegeneracy conditions seem much harder to understand in this 
context; compare to \cite{os19}.

\begin{defn}\label{defn:link-bilin} Given an $s$-prelinked
$\sE_{\bullet}$, and $m \in H_d$,
a \textbf{linked bilinear form} of \textbf{index} $m$ on
$\sE_{\bullet}$ is a collection of bilinear pairings for each $w,w' \in V(G)$
$$\left<,\right>_{w,w'}:\sE_w \times \sE_{w'} \to \sO_S$$
satisfying the following compatibility conditions: 
for any suitable $w,w'$, and $\e \in E(G)$ with $t(\e)=w$, we have
$$\left<,\right>_{h(\e),w'} \circ (f_{\e} \times \id) = 
s^{\delta_{w,w',\e}} \cdot \left<,\right>_{w,w'},$$
and for any $\e'$ with $t(\e')=w'$, we have
$$\left<,\right>_{w,h(\e')} \circ (\id \times f_{\e'}) = 
s^{\delta_{w,w',\e'}} \cdot \left<,\right>_{w,w'},$$
where 
$$\delta_{w,w',\e}=\begin{cases}1 : 
d(\frac{w+w'}{2},m)<d(\frac{h(\e)+w'}{2},m) \\
0: \text{otherwise,}\end{cases}$$
and 
$$\delta_{w,w',\e'}=\begin{cases}1 : 
d(\frac{w+w'}{2},m)<d(\frac{w+h(\e')}{2},m) \\
0: \text{otherwise,}\end{cases}.$$
\end{defn}

Note that if $w=w'$, so that the notation $\delta_{w,w',\e}$ is 
\textit{a priori} ambiguous, we have that both definitions coincide, so in 
fact there is no ambiguity.

The following lemma checks that our compatibility conditions are internally
consistent, and will be useful later. The notation is slightly cumbersome,
because we have to keep track of the order of applying maps on the left and
the right, but the idea is simply to see what our compatibility conditions
say when we apply sequences of maps on either side. 

\begin{lem}\label{lem:compat-consist} 
Given $w,w' \in V(G)$, and paths $P$ and $P'$ starting at $w,w'$ respectively,
with $P=(\e_1,\dots,\e_n)$ and $P'=(\e_{n+1},\dots,\e_{n+n'})$, fix 
$\sigma \in S_{n+n'}$ such that $\sigma^{-1}$ preserves the order of 
$\{1,\dots,n\}$ and $\{n+1,\dots,n+n'\}$.
For $i=1,\dots,n+n'$, set 
$$f_i=\begin{cases}
f_{\e_{\sigma(i)}} \times \id:&\sigma(i) \leq n \\
\id \times f_{\e_{\sigma(i)}}:&\sigma(i) > n,
\end{cases}$$
and let 
\begin{equation}\label{eq:genl-epsilon}
\delta_{P,P'} =\frac{1}{|V(\Gamma)|}
\left(n+n'+d\left(\frac{h(P)+h(P')}{2},m\right)
-d\left(\frac{w+w'}{2},m\right)\right).
\end{equation}
Then the compatibility conditions imposed
in Definition \ref{defn:link-bilin} give
$$\left<,\right>_{h(P),h(P')} \circ f_{n+n'} \circ \cdots \circ f_1 = 
s^{\delta_{P,P'}} \cdot \left<,\right>_{w,w'},$$
In fact, if we denote by $P_i$ for $i=1,\dots,n+n'$ the truncation of $P$
given by $(\e_j,\dots,\e_n)$ where $j$ is minimal in $\{1,\dots,n\}$ with
$\sigma^{-1}(j) \geq i$, and similarly for $P'_i$, we have
$$\sum_{i=1}^{n+n'} \delta_{t(P_i),t(P'_i),\e_{\sigma(i)}}
=\delta_{P,P'}.$$
\end{lem}

We use the convention that if $P_i$ is empty, then $t(P_i)=h(P)$, and
similarly for $P'_i$.

\begin{proof} Setting $w_i=t(P_i)$ and $w'_i=t(P'_i)$ for each $i$,
the value of each $\delta_{t(P_i),t(P'_i),e_{\sigma(i)}}$ is determined
by $d(\frac{w_{i+1}+w'_{i+1}}{2},m)-d(\frac{w_i+w'_i}{2},m)$. According
to our description of minimal paths, this difference is always either $-1$ or 
$|V(\Gamma)|-1$, and
$\sum_{i=1}^{n+n'} \delta_{w_i,w'_i,\e_{\sigma(i)}}$ can be expressed
as the number of times it is $|V(\Gamma)|-1$. Then the sum of these
differences is equal to
$\left(|V(\Gamma)|\sum_{i=1}^{n+n'} \delta_{w_i,w'_i,\e_{\sigma(i)}}\right)
-n-n'$, giving the desired formula.
\end{proof}
 
\begin{defn}\label{defn:link-alt} In the notation of Definition 
\ref{defn:link-bilin}, a linked bilinear form is a \textbf{linked alternating
form} if $\left<,\right>_{w,w}$ is an alternating form on $\sE_w$ for all 
$w$, and 
$$\left<,\right>_{w,w'}=-\left<,\right>_{w',w} \circ \sw_{w,w'}$$
for all $w \neq w'$, where $\sw_{w,w'}:\sE_w \times \sE_{w'} \to 
\sE_{w'} \times \sE_w$ is the canonical map switching factors.
\end{defn}

This definition is equivalent to requiring that the induced form on
$\bigoplus_w \sE_w$ be alternating.

Because being $s$-prelinked is preserved by base change, it makes sense to 
define moduli functors of linked bilinear forms and linked alternating 
forms. Moreover, it is clear that these functors are represented by schemes 
(compare Lemma 2.2 of \cite{os17}), and that they further have natural 
module structures (i.e., for each $S$-scheme $T$ the set of linked 
bilinear/alternating forms is a $\Gamma(T,\sO_T)$-module). Our first 
result is that for simple $s$-prelinked bundles, the moduli of linked 
bilinear forms and of linked alternating forms behave just like their 
classical counterparts, in the sense that the representing scheme is in 
fact a vector bundle of the appropriate rank.

\begin{prop}\label{prop:form-structure} Suppose
$\sE_{\bullet}=((\sE_w)_w,(f_{\e})_{\e})$ is
$s$-prelinked of rank $r$, and simple, and $m \in H_d$. Then
the moduli scheme of linked bilinear forms on $\sE_{\bullet}$ of index
$m$ is a vector bundle on $S$ of rank $r^2$, and
the moduli scheme of linked alternating forms on $\sE_{\bullet}$ of
index $m$ is a vector bundle on $S$ of rank $\binom{r}{2}$.
\end{prop}

\begin{proof} First, choose subbundles $\sW_w \subseteq \sE_w$ as provided
by Proposition \ref{prop:structure}. Clearly, a linked bilinear form on
$\sE_{\bullet}$ induces by restriction a collection of bilinear pairings
$$\left<,\right>_{w,w'}':\sW_w \times \sW_{w'} \to \sO_S,$$
or equivalently a bilinear form on $\bigoplus_w \sW_w$,
and our claim is that this restriction map induces an isomorphism of
functors from linked bilinear forms to bilinear forms on $\bigoplus_w \sW_w$.
Because $\sum_w \rk \sW_w = r$, the claim yields the first statement
of the proposition.

To prove the claim, suppose we have a collection of $\left<,\right>_{w,w'}'$
as above; we aim to construct an inverse to the restriction map. Because
$$\sE_w \cong \bigoplus_{w' \in G} f_{w',w}(\sW_{w'}),$$
in order to define $\left<v_1,v_2\right>_{w,w'}$, it is enough to do so for
$v_1$ in $f_{u,w} (\sW_{u})$ and $v_2$ in $f_{u',w'} (\sW_{u'})$ as $u,u'$ vary
in $V(G)$. Starting from 
the necessity of having $\left<,\right>_{w,w'}=\left<,\right>_{w,w'}'$ on
$\sW_w \times \sW_{w'}$, we then see that inductive application of the
compatibility conditions
of Definition \ref{defn:link-bilin} determine $\left<,\right>_{w,w'}$
uniquely: specifically, if $v \in \sW_u$ and $v' \in \sW_{u'}$, then we set
$$\left<f_{u,w}(v),f_{u',w'}(v')\right>_{w,w'}=s^{\delta_{P,P'}}
\left<v,v'\right>_{w,w'},$$
where $P$ and $P'$ are minimal paths from $u$ to $w$ and from $u'$ to $w'$
respectively, and $\delta_{P,P'}$ are as in Lemma \ref{lem:compat-consist}.
It remains to check that the result forms $\left<,\right>_{w,w'}$ satisfy
the conditions for a linked bilinear form, which is straightforward to
verify directly, again using Lemma \ref{lem:compat-consist}.
The claim then follows, as the preceding construction is visibly inverse
to the restriction map.

To obtain the second statement of the proposition, it is enough to observe
that under the isomorphism of functors constructed above, a linked
bilinear form is alternating if and only if the induced form on
$\bigoplus_w \sW_w$ is alternating. Indeed, this follows from the
symmetry of the compatibility conditions together with Lemma
\ref{lem:compat-consist}.
\end{proof}

\begin{defn}\label{defn:link-isotropy-locus} If $\sE_{\bullet}$ is
$s$-prelinked, with a linked bilinear form $\left<,\right>_{\bullet}$,
the \textbf{locus of isotropy} of $\left<,\right>_{\bullet}$ on $\sE_{\bullet}$
is the closed subscheme of $S$ representing the functor of morphisms
$T \to S$ such that $\left<,\right>_{\bullet}$ is identically zero
after restriction to $T$.
\end{defn}

The fact
that the locus of isotropy is represented by a closed subscheme is clear,
as $\sE_{\bullet}$ together with $\left<,\right>_{\bullet}$
induces a morphism from $S$ to the moduli scheme of linked bilinear forms
on $\sE_{\bullet}$, and the locus of isotropy is the preimage under
this morphism of the zero form.

Proposition \ref{prop:form-structure} thus implies:

\begin{cor}\label{cor:iso-codim} Suppose $\sE_{\bullet}$ is
$s$-prelinked and simple, and $\left<,\right>_{\bullet}$ is a linked 
bilinear (respectively, linked alternating) form on $S$. Then the locus 
of isotropy is locally cut out in $S$ by
$d^2$ (respectively, $\binom{d}{2}$) equations, and thus if $S$
is locally Noetherian, every component of this locus has codimension
at most $d^2$ (respectively, $\binom{d}{2}$) in $S$.
\end{cor}

\section{Linked symplectic forms}

In \cite{o-t1}, our definition of linked symplectic form consisted of
two nondegeneracy conditions, which we can think of as a condition on
the interior of $G$, and a condition on the boundary. With these conditions,
we were able to prove that linked symplectic Grassmannians are smooth of
the expected dimension at all simple points. While it is undoubtedly 
possible to follow the same process in our present generalized setting,
we instead simplify both the definition and arguments substantially via
the observation that for our purposes (i.e., for proving Theorem 
\ref{thm:main}), the boundary is irrelevant, and it is 
enough to restrict attention to a more restrictive class of ``internally'' 
simple points. Thus, our definitions will not precisely generalize those
of \cite{o-t1}, and accordingly we use slightly different terminology.

\begin{defn}\label{defn:link-symp} Given an $s$-prelinked
$\sE_{\bullet}$ on a scheme $S$, and $m \in H_d$,
a linked alternating form $(\left<,\right>_{w,w'})_{w,w' \in V(G)}$ of index 
$m$ on $\sE_{\bullet}$ is a \textbf{linked internally symplectic form}
if for each $w,w'\in V(G)$ with $\frac{w+w'}{2}=m$, we have 
$\left<,\right>_{w,w'}$ perfect.
\end{defn}

Here by perfect, we mean that $\left<,\right>_{w,w'}$ is nondegenerate after
restriction to any point of $S$.

\begin{defn}\label{defn:internal-simple} Let $K$ be a field over $S$, and 
for $s \in K$, let $E_{\bullet}$ be $s$-prelinked on $\Spec K$.
Given also $m \in H_d$, we say that $E_{\bullet}$ is \textbf{internally 
simple} relative to $m$ if
there exist $w_1,\dots,w_r \in V(G)$ (not necessarily distinct) and
$v_i \in F_{w_i}$ for $i=1,\dots,r$ such that:
\begin{Ilist}
\itm for every $w \in V(G)$,
there exist paths $P_1,\dots,P_g$ with each $P_i$ going from $w_i$ to
$w$, and such that $f_{P_1}(v_1),\dots,f_{P_r}(v_r)$ form a basis for
$E_w$;
\itm for each $i$, we have $2m-w_i \in V(G)$.
\end{Ilist}
\end{defn}

It follows from Nakayama's lemma that over an arbitrary base, being 
internally simple relative to $m$ is an open condition, so we may apply 
it over any base to mean simply that it holds over each point.

For our main calculation relating to linked symplectic forms, we need
to work with respect to subbundles, so we also make the following
definition.

\begin{defn}\label{defn:link-sub} Given an $s$-prelinked $\sE_{\bullet}$,
a collection $\sF_{\bullet}$ of subbundles of $\sE_{\bullet}$ of fixed
rank $r$ is a \textbf{linked subbundle} if for all $e \in G$, from $w$
to $w'$, we have $f_e (\sF_w) \subseteq \sF_{w'}$. If also we have a
linked alternating form $\left<,\right>_{\bullet}$ on $\sE_{\bullet}$, 
we say that $\sF_{\bullet}$ is \textbf{isotropic} for 
$\left<,\right>_{\bullet}$ if the restriction of $\left<,\right>_{\bullet}$
to $\sF_{\bullet}$ vanishes identically.
\end{defn}

The following construction arises in the analysis of tangent spaces of
prelinked alternating Grassmannians, and is therefore key to the proof
of Theorem \ref{thm:main-lag}.

\begin{defn}\label{defn:tangent-form} Suppose that 
$\sE_{\bullet}$ is $s$-prelinked, and 
$\left<,\right>_{\bullet}$ is a linked alternating form on $\sE_{\bullet}$
of index $m$. Let $\sF_{\bullet} \subseteq \sE_{\bullet}$ be a linked
subbundle which is internally simple relative to $m$,
and suppose that $\sF_{\bullet}$ is isotropic for 
$\left<,\right>_{\bullet}$. Finally, let $(\sW_w \subseteq \sF_w)_w$ be as 
in Lemma \ref{prop:structure}, with the nonzero $\sW_w$ supported on
$w \in V(G)$ with $2m-w \in V(G)$. Given
family of homomorphisms $(\vp_w:\sW_w \to \sE_w/\sF_w)_{w \in V(G)}$ 
define the
associated linked alternating form $\left<,\right>^{\vp_{\bullet}}_{\bullet}$
on $\sF_{\bullet}$ by applying the following formula on the $\sW_w$:
$$\left<,\right>^{\vp_{\bullet}}_{w,w'} =
\left<,\right>_{w,w'} \circ (\vp_w \times \id) + 
\left<,\right>_{w,w'} \circ (\id \times \vp_{w'}).$$ 
\end{defn}

Note that this is well-defined because $\sF_{\bullet}$ is assumed to be 
isotropic. Also, recall that by Proposition \ref{prop:form-structure},  
the pairings on the $\sW_w$ defined above uniquely determine a linked
alternating form $\left<,\right>^{\vp_{\bullet}}_{\bullet}$ on 
$\sF_{\bullet}$. The main consequence of the symplectic condition is
the following.

\begin{lem}\label{lem:symplectic-transverse} In the situation of
Definition \ref{defn:tangent-form}, suppose further that
$\left<,\right>_{\bullet}$ is an internally linked symplectic form, and 
that $S$ is a point. Then the 
map from $\bigoplus_{w \in V(G)} \Hom(\sW_w,\sE_w/\sF_w)$ to the space of
linked alternating forms on $\sF_{\bullet}$ is surjective.
\end{lem}

We will need a preliminary observation on linked bilinear forms:

\begin{prop}\label{prop:form-induced} Suppose that 
$\sE_{\bullet}$ is $s$-prelinked, and 
$\left<,\right>_{\bullet}$ is a linked bilinear form on $\sE_{\bullet}$
of index $m$. Given $w,w' \in V(G)$, suppose that $2m-w \in V(G)$. Then,
if $P$ is a minimal path from $w'$ to $2m-w$, we have
$$\left<,\right>_{w,w'}=\left<,\right>_{w,2m-w} \circ (\id \times f_P).$$
\end{prop}

\begin{proof} Observing that the length of $P$ is equal to
$d\left(\frac{w+w'}{2},m\right)$,
this is an immediate consequence of Lemma \ref{lem:compat-consist}.
\end{proof}

With this, the proof of the lemma is now straightforward.

\begin{proof}[Proof of Lemma \ref{lem:symplectic-transverse}]
If we choose bases $v^w_j$ for each $\sW_w$, it is clearly enough to prove 
that for all $w,u,p,q$ with either $w \neq u$ or $p \neq q$,
there exists a choice
of $\vp_{\bullet}$ such that the induced pairing
$\left<v^{w'}_{p'},v^{u'}_{q'}\right>^{\vp_{\bullet}}_{w',u'}$ is nonzero
precisely when $(w',u',p',q')$ is either $(w,u,p,q)$ or $(u,w,q,p)$.
Given $w,u,p,q$, we construct the desired $\vp_{\bullet}$ as follows.
First set $\vp_{w'}=0$ for all $w' \neq w$,
and set $\vp_{w}(v^{w}_{p'})=0$ for all $p' \neq p$. We then wish to show
that there exists a choice of $\vp_w(v^w_p) \in \sE_w$ such that
$\left<\vp_w(v^w_p),v^u_q\right>_{w,u} \neq 0$, but 
$\left<\vp_w(v^w_p),v^{u'}_{q'}\right>_{w,u'}=0$ for all other choices of
$u',q'$. Equivalently, if we denote by $\widehat{\sW}_{u} \subseteq \sW_u$ the 
span of the $v^u_{q'}$ for $q' \neq q$, we want 
$$\vp_w(v^w_p) 
\in (\widehat{\sW}_u)^{\perp} \cap \left(\cap_{u'\neq u} \sW_{u'}^{\perp}\right),$$
but
$$\vp_w(v^w_p) \not\in \cap_{u'} \sW_{u'}^{\perp}.$$
Here each orthogonal space should be taken with respect to the appropriate 
pairing. Now, by the internally simple hypothesis, we have that 
$2m-w \in V(G)$, 
so for each $u'$, choose a minimal path $P_{u'}$ from $u'$ to $2m-w$. By 
Proposition \ref{prop:form-induced}, the above conditions are equivalent to
having
\begin{align*} 
\vp_w(v^w_p) & \in \left(f_{P_u} \left(\widehat{\sW}_u\right)\right)^{\perp}
\cap \left(\bigcap_{u' \neq u} \left(f_{P_{u'}} 
\left(\sW_{u'}\right)\right)^{\perp}\right) \\
& =  \left(f_{P_u} \left(\widehat{\sW}_u\right) \oplus
\left(\bigoplus_{u' \neq u} f_{P_{u'}} 
\left(\sW_{u'}\right)\right)\right)^{\perp},
\end{align*}
but
\begin{align*} 
\vp_w(v^w_p) & \not\in 
\left(\bigcap_{u'} \left(f_{P_{u'}} 
\left(\sW_{u'}\right)\right)^{\perp}\right) \\
& =  \left(\bigoplus_{u'} f_{P_{u'}} \left(\sW_{u'}\right)\right)^{\perp}.
\end{align*}
The sums are direct
sums because of Proposition \ref{prop:structure}, and now all the orthogonal 
complements are relative to $\left<,\right>_{w,2m-w}$. Again by Proposition
\ref{prop:structure}, the two sums
give distinct subspaces of $\sF_{2m-w}$, so by the nondegeneracy of
$\left<,\right>_{w,2m-w}$ imposed in the definition of a linked internally
symplectic form, we conclude that a $\vp_w(v^w_p)$ satisfying the desired 
conditions exists, and the lemma follows. 
\end{proof}

\section{Applications to linked Grassmannians}

We begin with some preliminary background on (pre)linked Grassmannians.

\begin{defn}\label{def:lg} Let $\sE_{\bullet}$ be an $s$-prelinked bundle
of rank $t$. Given $r<t$, define
the \textbf{prelinked Grassmannian} $\LG(r,\sE_{\bullet})$ to be the
scheme representing the functor associating to an $S$-scheme $T$ the
set of all linked subbundles of rank $r$ of the pullback of $\sE_{\bullet}$
to $T$. 
\end{defn}

Recall that linked subbundles were defined previously in Definition
\ref{defn:link-sub}.

Thus, a point of a prelinked Grassmannian is itself an $s$-prelinked bundle,
of rank $r$. 

\begin{defn}\label{def:simple} 
We say a point of $\LG(r,\sE_{\bullet})$ is \textbf{simple} 
(respectively, \textbf{internally simple} relative to $m$) if the
induced $s$-prelinked bundle is simple (respectively, internally simple
relative to $m$).
\end{defn}

For our purposes, the main result we will need is the following:

\begin{prop}\label{prop:simple-smooth} On the locus of simple points,
$\LG(r,\sE_{\bullet})$ is smooth over $S$ of relative dimension $r(t-r)$.

When $S=\Spec K$ is a point, at a particular $K$-valued simple point 
$(F_w)_{w \in V(G)}$, we may describe the
tangent space as follows: given $W_w \subseteq F_w$ as in 
Proposition \ref{prop:structure}, the tangent space to $\LG(r,\sE_{\bullet})$
at $(F_w)_{w \in V(G)}$ is equal to
$$\bigoplus_{w \in V(G)} \Hom(W_w,\sE_w/F_w).$$
\end{prop}

See Proposition A.2.2 of \cite{os20}.

We can now proceed to define linked alternating Grassmannians, and 
we easily conclude our main result on them. 

\begin{defn}\label{def:LAG}
Given $\sE_{\bullet}$ $s$-prelinked with a linked alternating
form $\left<,\right>_{\bullet}$, the \textbf{prelinked alternating Grassmannian}
$\LAG(r,\sE_{\bullet},\left<,\right>_{\bullet})$ is the closed subscheme
of $\LG(r,\sE_{\bullet})$ parametrizing linked subbundles which are 
isotropic for $\left<,\right>_{\bullet}$.
\end{defn}

\begin{proof}[Proof of Theorem \ref{thm:main-lag}] By definition,
$\LAG(r,\sE_{\bullet},\left<,\right>_{\bullet})$ is precisely the isotropy
locus of $\left<,\right>_{\bullet}$ for the universal 
subbundle on $\LG(r,\sE_{\bullet})$. Since the statement
is local, we may restrict to the locus of simple points of 
$\LG(r,\sE_{\bullet})$. On this locus, the universal subbundle is simple,
so Corollary \ref{cor:iso-codim} implies that 
$\LAG(r,\sE_{\bullet},\left<,\right>_{\bullet})$ is
cut out locally by $\binom{r}{2}$ equations, as desired.

Next, recall that $\left<,\right>_{\bullet}$ induces a morphism
from $\LG(r,\sE_{\bullet})$ to the space of linked alternating forms on the 
universal subbundle, and
$\LAG(r,\sE_{\bullet},\left<,\right>_{\bullet})$ is precisely the 
pullback of the zero section under this morphism. By Proposition 
\ref{prop:form-structure}, this space is (on the simple locus) 
a vector bundle $\sG$ of rank $\binom{r}{2}$. Letting $s_0$ be the
zero section of $\sG$ and $s_1$ the section induced by 
$\left<,\right>_{\bullet}$, we may view 
$\LAG(r,\sE_{\bullet},\left<,\right>_{\bullet})$ as $s_0 \cap s_1$.
In order to prove the theorem, it is then
enough (see for instance Lemma 2.3 of \cite{os19}) to see that under the 
symplectic hypotheses, the fibers
of $\LAG(r,\sE_{\bullet},\left<,\right>_{\bullet})$ are smooth of the
correct dimension, or equivalently, that the tangent
spaces to $s_0$ and $s_1$ intersect transversely in the fiber over any
point of $S$. We may thus assume that $S$ is a point, and consequently the
$\sE_{\bullet}$ are simply vector spaces.

At any point of $s_0$, the tangent space of $\sG$ decomposes 
canonically as a direct
sum of the tangent space of $\LG(r,\sE_{\bullet})$ (which is described
by Proposition \ref{prop:simple-smooth}) and the tangent space to the moduli 
space of linked alternating forms on the corresponding fixed linked subspace.
Since the latter moduli space is a vector space, the tangent 
space is identified with the space itself. Then thinking of $z$ as a
point of $s_0 \cap s_1$, given a tangent vector $v$ to
$\LG(r,\sE_{\bullet})$ at a point, our tautological section $s_1$ induces a
tangent vector of $\sG$ at $z$, and hence also a tangent vector $v'$ of
the fiber of $\sG$ at $z$. We may then think of $v'$ as a linked alternating
form. One checks from the definitions that if the tangent vector $v$ is
represented by $(\vp_i:\sF_i \to \sE_i/\sW_i)_i$ for some choice of
$\sW_i$ as in Proposition \ref{prop:structure}, the resulting linked 
alternating form obtained from $v'$ is precisely the
$\left<,\right>^{\vp_{\bullet}}_{\bullet}$ of Definition
\ref{defn:tangent-form}.
Since tangent vectors to the zero section always yield the 
zero linked alternating form, transversality of the tangent spaces of
the two sections follows from the surjectivity of the map
$\vp_{\bullet} \to \left<,\right>^{\vp_{\bullet}}_{\bullet}$, given to
us by Lemma \ref{lem:symplectic-transverse}. We thus conclude the theorem.
\end{proof}

\section{Special line bundles on reducible curves}

In order to talk about higher-rank limit linear series with special
determinant, we need to develop the appropriate definition of special
line bundle on a curve of compact type. We make the following definition:

\begin{defn}\label{def:special} Let $X$ be a curve of compact type over a 
field, and $\sL$ a line bundle on $X$. We say $\sL$ is
\textbf{special} if there exists a morphism $\sL \to \omega_X$ which
is not everywhere zero on any component of $X$.
\end{defn}

The definition is justified by the following proposition:

\begin{prop}\label{prop:special}
Suppose that $\pi:X \to B$ is a smoothing family, with
$B$ the spectrum of a DVR, and with smooth generic fiber $X_{\eta}$.
Let $X_0$ be the special fiber. Then:

\begin{ilist}
\itm If $\sL_{\eta}$ is a special line bundle on $X_{\eta}$, then
it extends to a line bundle $\sL$ on $X$ such that $\sL|_{X_0}$ is
also special.
\itm If $\sL_0$ is a special line bundle on $X_0$, then after possible 
faithfully flat base change $B' \to B$ with $B'$ still the
spectrum of a DVR, there is a line bundle $\sL$ on $X'=X \times_B B'$
and a morphism $\sL \to \omega_{X'/B'}$ which does not vanish identically
on any component of the special fiber $X'_0$, and such that $\sL|_{X'_0}$
is isomorphic to the pullback of $\sL_0$.
\end{ilist}
\end{prop}

Recall the definition of smoothing family:

\begin{defn}\label{def:smoothing-fam} A morphism of schemes 
$\pi: X \rightarrow B$ constitutes a \textbf{smoothing family} if:

\begin{Ilist}
\itm $B$ is regular and connected;
\itm $\pi$ is flat and proper;
\itm The fibers of $\pi$ are genus-$g$ curves of compact type;
\itm Each connected component $\Delta'$ of the singular locus of $\pi$ maps
isomorphically onto its scheme-theoretic image $\Delta$ in $B$, and
furthermore $\pi ^{-1} \Delta$ breaks into two (not necessarily
irreducible) components intersecting along $\Delta'$;
\itm Any point in the singular locus of $\pi$ which is smoothed
in the generic fiber is regular in the total space of $X$.
\end{Ilist}
\end{defn}

We first prove the following lemma.

\begin{lem}\label{lem:g0ds} 
Let $X_0$ be a curve of compact type, and $w \in H_d$
a multidegree of total degree $d$, nonnegative on each component of $X_0$.
Then the space of 
$\fg^0_{w}$s on $X_0$ has the expected dimension $d$.
\end{lem}

Note that even though $X_0$ is reducible, we are considering linear series
in the classical sense, and not limit linear series.

\begin{proof} Let $G^0_{w}(X_0)$ be the space of 
$\fg^0_{w}$s on $X_0$, and denote by $G^{0,\nd}_{w}(X_0)$
the open subset on which the global section is not identically zero on
any component of $X_0$. Let $\Gamma$ be the dual graph of $X_0$, so that
$w$ consists of a non-negative integer $i_v$ for each 
$v \in V(\Gamma)$, and let $Y_v$ denote the component of $X_0$ corresponding 
to a given $v$.

The argument giving the lower bound on dimension of $G^0_{w}(X_0)$
is unaffected by the reducibility of $X_0$, so it suffices to prove that
the dimension is at most $d$.

Now, we have a natural restriction map
\begin{equation}\label{eq:g0ds}
G^{0,\nd}_{w}(X_0) \to \prod_{v \in V(\Gamma)} G^0_{i_v}(Y_v),
\end{equation}
and because $i_v \geq 0$ for each $v$, the target space has dimension
$\sum_{v} i_v = d$.
Because $\Gamma$ is a tree, gluing conditions are independent, so the 
image consists of tuples $(\sL_v,V_v)_{v \in V(\Gamma)}$ such that
if $V_v$ vanishes at a node corresponding to an edge $e$, and $v'$ is
the other vertex adjacent to $e$, then $V_{v'}$ likewise vanishes at the
same node. Such vanishing also determines the fibers of \eqref{eq:g0ds}:
if sections vanish at a node, \eqref{eq:g0ds} fails
to be injective due to the ability to scale independent on either side of
the node. However, this is the only ambiguity, so we see that if we have
$(\sL_v,V_v)_{v \in V(\Gamma)}$ in the image of \eqref{eq:g0ds}, and
$S \subseteq E(\Gamma)$ is the set of nodes at which the $V_v$ vanish,
then the fiber of \eqref{eq:g0ds} over the given point has dimension
$|S|$.
On the other hand, the dimension of the locus in 
$\prod_{v \in V(\Gamma)} G^0_{i_v}(Y_v)$ with given $S$ is $d-2|S|$,
so we conclude that $G^{0,\nd}_{w}$ has dimension at most, hence equal
to, $d$.

The argument on the complement of $G^{0,\nd}_{w}$ is similar;
let $S'$ be a nonempty subset of $V(\Gamma)$, and $G^{0,S'}_{w}(X_0)$ the
locally closed subscheme of $G^0_{w}(X_0)$ on which the
sections vanish identically precisely on $Y_v$ for $v \in S'$. Then we
again have a restriction map
\begin{equation}\label{eq:g0ds-2}
G^{0,S'}_{w}(X_0) \to \prod_{v \not\in S'} G^0_{i_v}(Y_v),
\end{equation}
and again using that every $i_v$ is nonnegative, we have that the target
space has dimension 
$d':=\sum_{v \not \in S'} i_v \leq d$. 
In this case, the image consists of
tuples $(\sL_v,V_v)_{v \not \in S'}$ such that if $v$ is adjacent to a
vertex in $S'$, then $V_v$ vanishes at the relevant node, and if $v,v'$
are adjacent and neither is in $S'$, then $V_v$ vanishes at the relevant
node if and only if $V_{v'}$ does. Fix a subset $S_1 \subseteq E(\Gamma)$
of edges such that neither adjacent vertex is in $S'$, and let $S_2$ 
consist of all edges such that exactly one adjacent vertex is in $S'$.
Consider the locus of tuples $(\sL_v,V_v)_{v \not \in S'}$ which vanish 
precisely on the nodes corresponding to elements of $S:=S_1 \cup S_2$. This 
locus has dimension
$d'-2|S_1|-|S_2|$.
On the other hand, the fiber dimension of \eqref{eq:g0ds-2} is equal
to $|S_1|+|S_2|-m$, where $m$ is the number of connected components of
$\cup_{v \in S'} Y_v$,
Thus, we get that $\dim G^{0,S'}_{w}(X_0) < d' \leq d$,
as desired.
\end{proof}

\begin{proof}[Proof of Proposition \ref{prop:special}] 
For (i), let $D_{\eta}$ be the divisor on $X_{\eta}$
giving the vanishing of a non-zero morphism 
$\sL_{\eta} \to \omega_{X_{\eta}}$.
Let $D$ be the closure of $D_{\eta}$ in $X$; since $X$ is assumed regular,
this is a Cartier divisor, so we can set $\sL=\omega_{X/B}(-D)$. This
visibly has the desired properties.

For (ii), the morphism $\sL_0 \to \omega_{X_0}$ can also be considered
a global section of $\sL^{-1}_0 \otimes \omega_{X_0}$, not vanishing
on any component of $X_0$. Setting $d=2g-2-\deg \sL_0$, and $w$
to be the multidegree of $\sL^{-1} \otimes \omega_{X_0}$, by Lemma
\ref{lem:g0ds} the space of $\fg^0_{w}$s on $X_0$ has 
expected dimension $d$, so our given one is the specialization of one
from $X_{\eta}$. It follows that after base change $B' \to B$,
we have a line bundle $\sM$ with a global section not vanishing identically 
on any component of $X'_0$, and such that $\sM|_{X'_0}$ is the pullback
of $\sL_0^{-1} \otimes \omega_{X_0}$, so finally setting 
$\sL=\sM^{-1} \otimes \omega_{X'/B'}$ we obtain the desired statement.
\end{proof}

\section{Limit linear series with special determinant}\label{sec:lls}

We now use Theorem \ref{thm:main-lag} to prove Theorem \ref{thm:main}.
This is in essence a more detailed presentation of the sketch given in
\cite{o-t1}, although certain aspects are more complicated.

We assume throughout that $X_0$ is a genus-$g$ curve of compact type with
dual graph $\Gamma$ over $\Spec K$, with $K$ an algebraically closed field,
and denote by $Y_v$ the component of $X_0$ corresponding to each
$v \in V(\Gamma)$. We also fix integers $d,k,b$ and $d_{\bullet}$ as in 
Situation \ref{sit:basic}, and consequent associated graph $G$. 
We first recall the definitions involved in the
statements of Theorem \ref{thm:main}, starting with 
limit linear series as generalized in \cite{te1} from
the construction of Eisenbud and Harris.

In order to facilitate precise statements of gluing conditions, we begin
by introducing notation for twisting bundles used to go between different
multidegrees.

\begin{notn}\label{not:twisting-bundles-ii}
For every vertex $v \in V(\Gamma)$, denote by $\sO_v$ the line bundle on
$X_0$ obtained as follows: let $Z_1,\dots,Z_n$ be the closures in $X$
of the connected components of $X\smallsetminus Y_v$, and for $i=1,\dots,n$,
let $\Delta'_i$ be the node of $X_0$ obtained as $Z_i \cap Y_v$.
We then define $\sO_{v}$ to be $\sO_{Y_{v}}(-\sum_i \Delta'_i)$ on
$Y_{v}$ and to be $\sO_{Z_i}(\Delta'_i)$ on each $Z_i$.

Given $w,w' \in V(G)$, let $P$ be a minimal directed path in $G$
from $w$ to $w'$. Let $v_1,v_2,\dots,v_{\ell}$ be the sequence in $V(\Gamma)$
induced by the edges making up $P$. Set
$$\sO_{w,w'}=\bigotimes_{i=1}^{\ell} \sO_{v_i}.$$
\end{notn}

Now, given $v \in V(\Gamma)$, we have the stack
$\cG^k_{r,d_v}(Y_v)$ of $\fg^k_{r,d_v}$'s on $Y_v$.
Given $e \in E(\Gamma)$ adjacent to $v$,
we have also the stack $\cM_{r}(\Delta'_e)$ of rank-$r$ vector
bundles on $\Delta'_e$, where $\Delta'_e$ is the node of $X_0$ corresponding 
to the edge $e$. There is a natural restriction map
$\cG^k_{r,d_v}(Y_v) \to \cM_r(\Delta'_e)$, but for reasons which will
soon become apparent, we instead fix a choice of $w_0 \in V(G)$ and
consider the map induced first by
twisting the universal bundle by $\left(\sO_{(w_v,w_0)}|_{Y_v}\right)$,
and then restricting to $\Delta'_e$. 

\begin{notn}\label{notn:pkrd}
Denote by $\cP^k_{r,d_{\bullet}}(X_0)$ the product of all the stacks
$\cG^k_{r,d_v}(Y_v)$ fibered over the stacks
$\cM_r(\Delta'_e)$ via the above maps.
\end{notn}

Note that while $\cM_r(\Delta'_e)$ has a single point, that point has
automorphism group $\GL_r$, so the stack is non-trivial. The purpose
of fibering over it and twisting each vector bundle to multidegree $w_0$
is that, due to the definition of $2$-fibered products,
this process precisely introduces a choice of gluing map at each node.
Thus, another way to express $\cP^k_{r,d_{\bullet}}(X)$ is as a proper
scheme over $\cM_{r,w_0}(X_0)$ 
whose fibers parametrize tuples of spaces of global sections of the
restriction of the given bundle to each component $Y_v$, twisted by 
$\sO_{(w_0,w_v)}|_{Y_v}$ to obtain degree $d_v$.
Here $\cM_{r,w_0}(X_0)$ denotes the moduli stack of vector bundles of
rank $r$ and multidegree $w_0$ on $X_0$.

We now recall the definition of Eisenbud-Harris-Teixidor limit linear
series.  To minimize confusion, we will use
superscripts to index vector bundles on irreducible components, and
subscripts to index bundles on the entire curve. For
convenience in gluing map notation, we choose directions for all
edges of $\Gamma$.

\begin{defn}\label{def:eht-lls}
Let $((\sE^v,V^v)_{v \in V(\Gamma)},(\vp_e)_{e \in E(\Gamma)})$ be a
$K$-valued point of
$\cP^k_{r,d_{\bullet}}(X)$, where $(\sE^v,V^v)$ is the corresponding
point of $\cG^k_{r,d_v}(Y_v)$, and if $e$ is an edge from $v$ to $v'$,
$$\vp_e:
\left(\sE^v \otimes \left(\sO_{(w_v,w_0)}|_{Y_v}\right)\right)|_{\Delta'_e}
\risom
\left(\sE^{v'} \otimes \left(\sO_{(w_{v'},w_0)}|_{Y_{v'}}\right)\right)
|_{\Delta'_e}$$
is the corresponding gluing isomorphism. Then $((\sE_v,V_v)_v,(\vp_e)_e)$
is an \textbf{Eisenbud-Harris-Teixidor limit linear series}
if for each $e \in E(\Gamma)$ adjacent to $v$ and $v'$, we have:
\begin{Ilist}
\itm
$H^0(Y_v,\sE_v(-(b+1)\Delta'_e))=0$ and similarly for $v'$;
\itm
if $a^{e,v}_1,\dots,a^{e,v}_k$ and $a^{e,v'}_1,\dots,a^{e,v'}_k$
are the vanishing sequences at $\Delta'_e$ for $(\sE_v,V_v)$ and
$(\sE_{v'},V_{v'})$ respectively, then for every $i$ we have
\begin{equation}\label{eq:eht-compat} a^{e,v}_i+a^{e,v'}_{k+1-i} \geq b;
\end{equation}
\itm
there exist bases $s^{e,v}_1,\dots,s^{e,v}_k$ and
$s^{e,v'}_1,\dots,s^{e,v'}_k$ of $V_v$ and $V_{v'}$ respectively such
that $s^{e,v}_i$ has vanishing order $a^{e,v}_i$ at $\Delta'_e$ for each $i$,
and similarly for $s^{e,v'}_i$, and if we have $a^{e,v}_i+a^{e,v'}_{k+1-i}=b$
for some $i$, then
\begin{equation}\label{eq:eht-gluing}\vp_e(s^{e,v}_i)=s^{e,v'}_{k+1-i};
\end{equation}
\end{Ilist}

We say that $((\sE_v,V_v)_v,(\vp_e)_e)$ is \textbf{refined} if equality
always holds in \eqref{eq:eht-compat}.
\end{defn}

In imposing \eqref{eq:eht-gluing}, we view each section in the appropriate
twist determined by its order of vanishing; see Remark 4.1.4 of \cite{os20}
for details. \textit{A priori}, the Eisenbud-Harris-Teixidor limit linear
series only form a set, but in Definition 4.2.1 of \cite{os20}
we endow it with a natural stack structure.

\begin{notn}\label{not:gkrd-eht}
We denote by $\cG^{k,\EHT}_{r,d_{\bullet}}(X_0)$ the stack of
Eisenbud-Harris-Teixidor limit linear series.
\end{notn}

A limit linear series of fixed determinant is essentially as one would
expect, with the only subtlety being that twisting is allowed after
taking the determinant. This eliminates parity restrictions on degrees, 
which turn out to be unnecessary.

\begin{defn}\label{def:eht-det}
Given a line bundle $\sL$ on $X_0$ of multidegree $w \in H_d$, 
the stack $\cG^{k,\EHT}_{r,\sL,d_{\bullet}}(X_0)$ of 
\textbf{Eisenbud-Harris-Teixidor limit linear series of
determinant $\sL$} is the (2-)fibered product of 
$\cG^{k,\EHT}_{r,d_{\bullet}}(X_0)$ with $\Spec K$ over the Picard stack
$\cPic^{w}(X_0)$, with the map 
$\cG^{k,\EHT}_{r,d_{\bullet}}(X_0) \to \cPic^{w}(X_0)$
given by taking the determinant of the underlying bundle
and twisting appropriately, and the map 
$\Spec K \to \cPic^{w}(X_0)$ obtained by considering $\sL$ as a point on
the Picard stack.

Thus, if we write $\cM_{r,\sL}(X_0)$ for the moduli stack of vector bundles
of rank $r$ and determinant $\sL$ on $X_0$, there is a forgetful morphism
$\cG^{k,\EHT}_{r,\sL,d_{\bullet}}(X_0) \to \cM_{r,\sL}(X_0)$.
\end{defn}

To make the twist in Definition \ref{def:eht-det} more precise, we 
introduce the following notation, which is essentially the rank-$1$ 
analogue of Notation \ref{not:twisting-bundles-ii}:

\begin{notn}\label{not:twist-rk-1} Given $w,w' \in H_d$, let 
$$\sO'_{w,w'}=\bigotimes_{i} \sO_{v_i},$$
where $v_1,\dots,v_n \in V(\Gamma)$ is a minimal sequence such that
for each $v \in V(\Gamma)$, the $v$th entry of $w'-w$ is given by
$$\#\{i:v_i \text{ is adjacent to } v\} - \ell \cdot \#\{i:v_i=v\},$$
where $\ell$ is the valence of $v$.
\end{notn}

Then the twist used in Definition \ref{def:eht-det} after taking the 
determinant consists of tensoring with $\sO'_{w_0,w}$.

Thus, the determinant $\sL$ condition for a given limit linear series 
requires that we fix an isomorphism between $\sL$ and the appropriate twist 
of the determinant of an underlying vector bundle. In particular, this 
includes extra data which has the effect of rigidifying the stack slightly,
and increasing the dimension by $1$. We denote the isomorphism
$\det \sE \otimes \sO'_{w_0,w} \risom \sL$ by $\psi$. Although the 
particular choice of $\psi$ will never be relevant, having made a choice
will affect dimension counts in the stack context.

We also recall the definition of chain-adaptability. 

\begin{defn}\label{def:adapted}
Let $X$ be a smooth projective curve over $\Spec F$, and $(\sE,V)$ a pair
with $\sE$ a vector bundle of rank $r$ on $X$, and $V$ a $k$-dimensional
space of global sections. Given points $P,Q \in X(F)$, let 
$a_1(P),\dots,a_k(P)$ and $a_1(Q),\dots,a_k(Q)$ be the vanishing sequences
at $P$ and $Q$ respectively. Then we say that a basis
$s_1,\dots,s_k \in V$
is \textbf{$(P,Q)$-adapted} if $\ord_P s_i = a_i(P)$ and
$\ord_Q s_i = a_{k+1-i}(Q)$ for $i=1,\dots,k$. We say that $(\sE,V)$
is \textbf{$(P,Q)$-adaptable} if there exists a $(P,Q)$-adapted basis of $V$.
\end{defn}

\begin{defn}\label{def:chain-adapted}
Let $X_0$ be a curve consisting of a chain of smooth projective curves
$X_1,\dots,X_n$ over $\Spec F$, with $P_i,Q_i \in X_i(F)$ for each $i$,
and the point $Q_i$ on $X_i$ glued to $P_{i+1}$ on $X_{i+1}$. Then a
refined Eisenbud-Harris-Teixidor limit linear series on $X_0$ is \textbf{chain
adaptable} if the pair induced by restriction to each $X_i$ ($i=2,\dots,n-1$)
is $(P_{i-1},P_i)$-adaptable.
\end{defn}

Now all of the concepts involved in Theorem \ref{thm:main} have been
defined. The reason for stating the theorem in terms of the above-defined 
limit linear series is that they are comparatively tractable in practice. 
However, they do not lend themselves easily to
theoretic arguments, such as the desired smoothing theorem. In contrast, we 
now recall the notion of linked linear series, which are not effective
computational tools, but are more amenable to theoretical constructions.
In the language of \cite{os20}, we will be using ``type II linked linear
series,'' the type I variant not being relevant for our purposes. In
fact, we will be using a slight variant which is shown to be equivalent
in Proposition 3.4.12 of \cite{os20}.

Because we wish to prove smoothing theorems, we have to work in families.
The only change to our running hypotheses is that in place of $X_0$,
we consider a smoothing family $\pi:X \to B$, and in addition to the case 
that the base $B$ is a point as considered above,
we also allow $B$ to be the spectrum of a $DVR$, with $\pi:X \to B$ a
smoothing family having smooth generic fiber. In this case, $\Gamma$ will
denote the dual graph of the special fiber, and $Y_v$ is the appropriate
component of the special fiber. Accordingly, we have to
generalize our definition of twisting bundles.

\begin{notn}\label{not:twisting-bundles-ii-2}
In the case that $B$ is the spectrum of a DVR, we have
$Y_{v}$ a Cartier divisor in $X$, and we set $\sO_{v}=\sO_X(Y_v)$.

We then define $\sO_{w,w'}$ as before.

Finally, for any $B$, given $w_0 \in V(G)$, and a $T$-valued point $\sE$ of
$\cM_{r,w_0}(X/B)$, for $w \in V(G)$ write
$$\sE_w:=\sE \otimes \sO_{w_0,w}|_{X \times_B T}.$$
\end{notn}

We will need to choose some additional data to define linked linear series,
which we use to construct maps as described below.

\begin{notn}\label{not:mors}
For each vertex $v \in V(\Gamma)$, suppose we fix a morphism
$$\iota_{v}:\sO_X \to \sO_{v}$$
vanishing precisely on $Y_v$.

Next, observe that $\bigotimes_{v \in V(\Gamma)} \sO_{v} \cong \sO_X$. Fixing
such an isomorphism (unique up to an element of $\sO_B^*$),
we obtain a induced morphism
$$\iota'_{v}:\bigotimes_{v' \neq v} \sO_{v'} \overset{\iota_v}{\to}
\bigotimes_{v' \in V(\Gamma)} \sO_{v'} \to \sO_X.$$

Finally, fix $w_0 \in V(G)$. For any edge $\e$ in $G$, let
$w$ be the tail and $w'$ the head, and let $v$ be the associated edge of
$\Gamma$. If $\sE$ is a $T$-valued point of $\cM_{r,w_0}(X/B)$
then we either have $\sE_{w'} = \sE_{w} \otimes \sO_{v}$
or $\sE_{w}=\sE_{w'} \otimes \bigotimes_{v' \neq v} \sO_{v'}$.
Thus, using $\iota_{v}$ or $\iota'_{v}$, and pushing forward
under $\pi$, we obtain a morphism
$$f_{\e}: \pi_* \sE_w \to \pi_* \sE_{w'}.$$
\end{notn}

\begin{defn}\label{def:grd-space}
Choose a vertex
$w_0 \in V(G)$. In addition, for each $v \in V(\Gamma)$, choose
a morphism $\iota_v: \sO_X \to \sO_v$ which vanishes precisely on $Y_v$.
The moduli stack $\cG^{k}_{r,w_0,d_{\bullet}}(X/B)$
of \textbf{linked linear series} is the category fibered
in groupoids over $B-\Sch$ whose
objects consist of tuples
$(S,\sE,(\sV_w)_{w \in V(G)})$, where $S$ is
a $B$-scheme, $\sE$ is a vector bundle of rank $r$ and multidegree
$w_0$ on $X \times_B S$, and $\sV_w$ is a rank-$k$ subbundle (in the
sense of Definition B.2.1 of \cite{os20})
of $\pi_* (\sE_w)$, satisfying the following conditions:
\begin{Ilist}
\itm for every $v \in V(\Gamma)$, every $e \in E(\Gamma)$ adjacent to $v$,
and every $z \in S$ mapping to the closed point of $B$,
we have
$$H^0(Y,\sE_{w_v}|_{Y}
(-(b+1)(\Delta')))=0,$$
where $Y$ is the component of the fiber $X_z$ corresponding to $v$, 
and $\Delta'$ denotes the node corresponding to $e$.
\itm For every edge $\e$ in $G$, let
$w$ be the tail and $w'$ the head. Then we require that
$$f_{\e}(\sV_w) \subseteq \sV_{w'}.$$
\end{Ilist}

Given $B' \to B$ and a line bundle $\sL$ on $X_{B'}$ of multidegree 
$w \in H_d$, the moduli stack 
$\cG^{k}_{r,\sL,d_{\bullet}}(X/B)$
of \textbf{linked linear series of determinant $\sL$} is the (2-)fibered
product of $\cG^k_{r,w_0,d_{\bullet}}(X/B)$ with $B'$ over
$\cPic^{w}(X/B)$, as in the case of limit linear series.
\end{defn}

In \cite{os20}, linked linear series are shown
to be parametrized by algebraic stacks.
It is always possible to increase $b$ and the $d_v$, in which case one
gets an open immersion into a larger moduli stack; see Proposition 6.2.1
of \cite{os20}. Although in general the stacks of
limit linear series and linked linear series
are rather different, the importance of chain-adaptability from our point
of view is that it describes an open substack on which the two 
constructions are isomorphic.

\begin{thm}\label{thm:grd-compare} Suppose $B$ is a point, and write
$X_0=X$. Then there is a forgetful morphism
\begin{equation}\label{eq:forget-comp} \cG^{k}_{r,w_0,d_{\bullet}}(X) \to
\cG^{k,\EHT}_{r,d_{\bullet}}(X)
\end{equation}
which is an isomorphism over the open substack of chain-adaptable limit
linear series.
\end{thm}

See Theorem 4.3.4 and Corollary 5.2.7 of \cite{os20}.

Thus, in order to prove Theorem \ref{thm:main}, we may restrict our
attention to linked linear series. A slight specialization of Theorem
6.1.4 of \cite{os20} says the following.

\begin{thm}\label{thm:foundation} Let $\pi:X\to B$ be a smoothing family with
$B$ being the spectrum of a DVR, and having smooth generic fiber $X_{\eta}$, 
and special fiber $X_0$. Given also
$k,r,d_{\bullet}$ as in Situation \ref{sit:basic}, let
$\sL$ be a line bundle of multidegree $w \in H_d$ on $X \times_B B'$ for
some $B$-scheme $B'$.
Then $\cG^{k}_{r,\sL,d_{\bullet}}(X/B)$ is an Artin stack over $B'$,
and the natural map
$$\cG^{k}_{r,\sL,d_{\bullet}}(X/B) \to 
\cM_{r,\sL}(X_{B'}/B')$$
is relatively representable by schemes of finite type.
Moreover, formation of the stack
$\cG^{k}_{r,\sL,d_{\bullet}}(X/B)$
is compatible with base change $B'' \to B$, and
in particular, if $\eta'$ is a point of $B'$ over the generic point of $B$,
the base change to $\eta'$ parametrizes
triples $(\sE,\sL,V)$ of
a vector bundle $\sE$ of rank $r$ and degree $d$ on $X_{\eta'}$ together with
an isomorphism $\det \sE \risom \sL|_{\eta'}$ and a $k$-dimensional vector
space $V \subseteq H^0(X_{\eta'},\sE)$.

Moreover, the simple locus of $\cG^{k}_{r,\sL,d_{\bullet}}(X/B)$ has
universal relative dimension at least
$k(d-k-r(g-1))$ over $\cM_{r,\sL}(X'/B')$, and therefore
universal relative dimension at least $\rho-g$ over $B'$.
In particular, if 
$\cG^{k}_{r,\sL,d_{\bullet}}(X_0)$ has dimension exactly $\rho-g$ at a
simple point $z$, then $\cG^{k}_{r,\sL,d_{\bullet}}(X/B)$ is universally
open over $B$ at $z$, and has fibers of pure dimension $\rho-g$ in an open 
neighborhood of $z$.
\end{thm}

The notion of universal relative dimension is introduced in \cite{os21}.

The above theorem is valid for any determinant $\sL$, but the last part
is not useful when $\sL$ is special, since the fiber will always have
dimension at least $\rho-g+\binom{k}{2}$. Thus, our 
main task is to prove the following variant which
addresses the case of special determinants.

\begin{thm}\label{thm:foundation-special} In the situation of
Theorem \ref{thm:foundation}, suppose further that $\sL$ is special.

Then the simple locus of $\cG^{k}_{r,\sL,d_{\bullet}}(X/B)$ has universal 
relative dimension at least $\rho-g+\binom{k}{2}$ over $B$.
In particular, if
$\cG^{k}_{r,\sL,d_{\bullet}}(X_0)$ has dimension exactly 
$\rho-g+\binom{k}{2}$ at a simple point $z$, then the structure morphism
is universally open in a neighborhood of $z$, with pure fiber dimension 
$\rho-g+\binom{k}{2}$.
\end{thm}

In families, we say $\sL$ is special if there is a morphism 
$\sL \to \omega_{X/B}$ which does not vanish uniformly on any component
of any fiber.

Before describing the limit case, we briefly recall a derivation for
smooth curves (see \S 2 of \cite{b-f2} or Theorem 4.2 of \cite{mu2})
of the modified expected dimension
in the special case of rank $2$ and canonical determinant, via the
following alternate construction of the moduli space in question.
Let $\cM_{2,\omega}(X)$ be the moduli stack of vector bundles
of rank $2$ and fixed canonical determinant on a smooth curve $X$ of genus
$g$; this is smooth of
dimension $3g-3$. Let $\widetilde{\sE}$ be the universal bundle on
$\cM_{2,\omega}(X) \times X$, and let $D$ be a sufficiently ample effective
divisor on $X$ (technically, we must cover $\cM_{2,\omega}(X)$ by a nested
increasing sequence of quasicompact open substacks, and carry out this
construction on each, letting $D$ grow). Let $D'$ be the pullback of $D$ to
$\cM_{2,\omega}(X) \times X$. Then $p_{1*} \widetilde{\sE}(D')$ is a vector 
bundle of rank 
$$\deg \widetilde{\sE}+\rk \widetilde{\sE} \deg D+ \rk \widetilde{\sE}(1-g)=2g-2+2 \deg D
+2-2g= 2 \deg D.$$
Let $G:=G(k,p_{1*} \widetilde{\sE}(D'))$ be the
relative Grassmannian on $\cM_{2,\omega}(X)$; our moduli space 
$\cG^k_{2,\omega}(X)$ is cut out
by the closed condition of subspaces lying in $p_{1*}\widetilde{\sE}$. We 
express this
condition in terms of the bundle $p_{1*}(\widetilde{\sE}(D')/\widetilde{\sE}(-D'))$,
which has rank $4 \deg D$. We see that because $D$ was chosen to be large,
$p_{1*} \widetilde{\sE}(D')$ is naturally a subbundle, as is 
$p_{1*}(\widetilde{\sE}/\widetilde{\sE}(-D'))$, which also has rank $2 \deg D$. 
Then the inclusion of the universal subbundle on $G$, together with the 
pullback from $\cM_{2,\omega}(X)$ of 
$p_{1*}(\widetilde{\sE}/\widetilde{\sE}(-D'))$, induces a morphism
\begin{equation}\label{eq:incid-corr}
G \to G(k,p_{1*}(\widetilde{\sE}(D')/\widetilde{\sE}(-D')) 
\times_{\cM_{2,\omega}(X)}
G(2 \deg D,p_{1*}(\widetilde{\sE}(D')/\widetilde{\sE}(-D')),\end{equation}
and our desired moduli space is precisely the preimage in $G$ of the
incidence correspondence in the product. The incidence correspondence is
smooth inside a smooth product space over $\cM_{2,\omega}(X)$, so is a
local complete intersection, and we thus obtain an upper bound on the 
codimension of $\cG^k_{2,\omega}(X)$ inside $G$.

However, this lower bound is not sharp: we next make use of the canonical 
determinant hypothesis to observe that
by choosing local representatives, using the isomorphism
${\bigwedge}^2 \widetilde{\sE} \cong p_{2}^* \omega_X$, and summing residues
over points of $D$, we obtain a symplectic form on 
$p_{1*}(\widetilde{\sE}(D')/\widetilde{\sE}(-D'))$, Moreover, both 
$p_{1*} \widetilde{\sE}(D')$ and $p_{1*}(\widetilde{\sE}/\widetilde{\sE}(-D'))$ are
isotropic for this form, with the former following from the residue theorem,
and the latter from the lack of poles. Thus, \eqref{eq:incid-corr} in fact
has its image in a product of symplectic Grassmannians, and the incidence
correspondence has smaller codimension, and is still a local complete
intersection, so we obtain the modified codimension
bound for $\cG^k_{2,\omega}(X)$ cut out in $G$. In fact, in the language of
\cite{os21}, we obtain that the moduli space has universal relative 
dimension at least $\rho-g+\binom{k}{2}$ from Corollary 7.7 of \cite{os21}.

We now prove the corresponding dimension bounds for limit linear series, by
following the above argument and applying Theorem \ref{thm:main-lag}.

\begin{proof}[Proof of Theorem \ref{thm:foundation-special}]
First, Proposition 6.2.1 of \cite{os20}
shows that suitable choices for increasing $d_{\bullet}$ and 
$b$ will always induce an open immersion of the corresponding spaces
$\cG^k_{r,\sL,d_{\bullet}}(X/B)$, so to prove a local dimension bound
we may increase $d_{\bullet}$ and $b$ as much as we wish.
In particular, if $z$ is a simple point of 
$\cG^k_{r,\sL,d_{\bullet}}(X/B)$, then we can increase $d_{\bullet}$ and
$b$ so that $z$ becomes internally simple. Now, let $\widetilde{\sE}_0$ 
be the universal bundle on $\cM_{2,\sL}(X'/B')$, and for each $w \in V(G)$, 
let $\widetilde{\sE}_w$ be the corresponding twist of 
$\widetilde{\sE}_0$. Then for every edge $\e$ of $G$, we also have 
an induced morphism 
$\tilde{f}_{\e}:\widetilde{\sE}_{t(\e)} \to \widetilde{\sE}_{h(\e)}$. 
Let $D$ and $D'$ be as in the smooth case described above, with the added 
condition that $D$ be supported on the smooth locus of $\pi$. Then we
can consider for each $w\in V(G)$ the vector bundle
$p_{1*}(\widetilde{\sE}_w(D')/\widetilde{\sE}_w(-D'))$, with morphisms 
$f_{\e}$ induced by $\tilde{f}_{\e}$. This creates an $s$-prelinked structure
$\widetilde{\sG}_{\bullet}$,
and the relative moduli space $\cG^k_{r,\sL,d_{\bullet}}(X/B)$ of linked 
linear series may be realized inside the linked Grassmannian 
$\LG(k,\widetilde{\sG}_{\bullet})$. 
As in the smooth case, $\cG^k_{r,\sL,d_{\bullet}}(X/B)$ is described as 
the locus on which the universal subbundle is contained inside both
$p_{1*}(\widetilde{\sE}_w(D'))$ and
$p_{1*}(\widetilde{\sE}_w/\widetilde{\sE}_w(-D'))$ for each $w$, so we are in 
exactly the situation of \eqref{eq:incid-corr}, with linked Grassmannians
in place of Grassmannians.

We thus need to see that the special determinant
hypothesis gives us (at least locally on $\cM_{2,\sL}(X'/B')$) a linked 
alternating form on $\widetilde{\sG}_{\bullet}$, which is internally symplectic
with respect to the multidegree $m$ of $\sL$, considered in $H_d$, and
such that the subbundles in question are isotropic for this form.
Now, by definition, we have maps
$$\sE_{w_0} \otimes \sE_{w_0} \otimes \sO'_{w_0,m} 
\twoheadrightarrow \left(\bigwedge^2 \sE_{w_0}\right) \otimes \sO'_{w_0,m}
\risom \sL \to \omega_{X'/B'}.$$
Given $w,w' \in V(G)$, the line bundle 
$\sO_{w,w_0} \otimes \sO_{w',w_0} \otimes \sO'_{w_0,m}$ 
is of the form $\bigotimes_i \sO_{v_i}$ for a sequence of $v_i \in V(\Gamma)$,
and if we let $v'_j$ be the sequence obtained from the $v_i$ by removing
all copies of $V(\Gamma)$, and set $\sO_{w,w',m}= \bigotimes_j \sO_{v'_j}$,
then we have 
$$\sO_{w,w',m} \cong \sO_{w,w_0} \otimes \sO_{w',w_0} \otimes \sO'_{w_0,m},$$ 
so we get induced morphisms
\begin{multline*} \sE_{w} \otimes \sE_{w'} 
= \sE_{w_0} \otimes \sO_{w_0,w} \otimes \sE_{w_0} \otimes \sO_{w_0,w'} 
\to \sE_{w_0} \otimes \sO_{w_0,w} \otimes \sE_{w_0} \otimes \sO_{w_0,w'} 
\otimes \sO'_{\frac{w+w'}{2},m}  \\
\risom \sE_{w_0} \otimes \sE_{w_0} \otimes \sO'_{w_0,m} 
\twoheadrightarrow \sL \to \omega_{X'/B'},\end{multline*}
which if $w=w'$ will factor through $\bigwedge^2 \sE_w$.

Applying the same method of choosing representatives and summing over
residues as described in the smooth case above, we obtain a linked 
alternating form on $\widetilde{\sG}_{\bullet}$.
Now, if $w+w'=2m$, we have that $\sO'_{\frac{w+w'}{2},m}$ is trivial,
so one checks by a local calculation that our linked alternating form
is in fact internally symplectic.
As in the smooth case, it is clear that pairing sections of
$p_{1*}(\widetilde{\sE}_w/\widetilde{\sE}_w(-D'))$ with
$p_{1*}(\widetilde{\sE}_{w'}/\widetilde{\sE}_{w'}(-D'))$ does not gives any
poles for the associated differentials, so all residues vanish, while
pairing $p_{1*}(\widetilde{\sE}_{w}(D'))$ with 
$p_{1*}(\widetilde{\sE}_{w'}(D'))$
gives a global section of $\omega_{X'/B'}(2D)$, and again the residue 
theorem implies that summing over residues gives zero. Thus, both of 
these linked subbundles are isotropic.

We can thus put together the canonical determinant construction for the
smooth case with the linked linear series construction simply by replacing 
the symplectic Grassmannians arising in the smooth case with prelinked 
alternating Grassmannians. Using that our form is internally symplectic,
Theorem \ref{thm:main-lag} tells us that these spaces are smooth 
of dimension equal to the usual symplectic Grassmannian (at least,
on the open locus of internally simple points), so the dimension count in the 
linked linear
series case goes through exactly as in the case of smooth curves, and
we obtain the desired lower bound on universal relative dimension. 
The universal openness assertion then follows from Corollary 7.4 of
\cite{os21}.
\end{proof}

We conclude by describing how our main theorem follows from Theorem
\ref{thm:foundation-special}. Recall that rather than using standard 
stability for reducible curves, we will find it more convenient to use the 
notion of \el-stability introduced in \cite{os22}. For the sake of 
propriety, we give the relevant definition prior to the proof, although it 
will not be used.

\begin{defn}\label{def:l-stable}
Let $\sE$ be a vector bundle of rank $r$ on a nodal $X$. We say that 
$\sE$ is $\boldsymbol{\ell}$-\textbf{semistable} (respectively,
$\boldsymbol{\ell}$-\textbf{stable}) if
for all proper subsheaves $\sF \subseteq \sE$ having constant rank
$r'$ on every component of $X$, we have
\begin{equation}\label{eq:defn}
\frac{\chi(\sF)}{r'} \leq \frac{\chi(\sE)}{r} \left(\text{respectively, }
\frac{\chi(\sF)}{r'} < \frac{\chi(\sE)}{r}\right).
\end{equation}
\end{defn}

Note that $\ell$-(semi)stability is trivially equivalent
to usual (semi)stability in the case that $X$ is smooth.

\begin{defn}\label{def:l-stable-lls} A limit linear series
$((\sE^v,V^v)_v,(\vp_e)_e)$ 
of degree $d$ on $X$ is $\boldsymbol{\ell}$-\textbf{semistable} 
(respectively, $\boldsymbol{\ell}$-\textbf{stable}) if the underlying
vector bundle $\sE$ of multidegree $w_0$ on $X$ is \el-semistable
(respectively, \el-stable).
\end{defn}

\begin{proof}[Proof of Theorem \ref{thm:main}] With $X_0$ as in the
theorem statement, let $\pi:X \to B$ be a smoothing family with $B$
the spectrum of a DVR, with special fiber $X_0$, and with the generic fiber 
$X_{\eta}$ being smooth. According to Theorem \ref{thm:grd-compare}, the
hypotheses of the theorem imply that the simple locus of
$\cG^{k}_{2,d_{\bullet},\sL_0}(X_0)$ has at least one component of dimension 
equal to $\rho-g+\binom{k}{2}$. By Proposition \ref{prop:special}, we have
that after possible faithfully flat base change $B' \to B$ with $B'$ still
the spectrum of a DVR, $\sL_0$ is the restriction of some $\sL$ on 
$X':=X \times _B B'$, with $\sL_{\eta}:=\sL_{X'_{\eta}}$ also special.
Then by Theorem \ref{thm:foundation-special} it
follows that $\cG^k_{2,\sL_{\eta}}(X'_{\eta})$ has at least one component of
dimension equal to $\rho-g+\binom{k}{2}$. Since $k \geq 2$, by Theorem 1.3 
of \cite{os16}, we conclude that $h^1(X'_{\eta},\sL_{\eta})$ must be equal 
to $1$. In particular, we must have $d \geq g-2$, so that on a general curve
a general special line bundle of degree $d$ has $h^1$ equal to $1$.

We have thus produced a single smooth curve $X_{\eta}$ with
line bundle $\sL_{\eta}$ such that $h^1(X_{\eta},\sL_{\eta})=1$ and
$\cG^k_{2,\sL_{\eta}}(X'_{\eta})$ has at least one component of
dimension equal to $\rho-g+\binom{k}{2}$. To conclude the statement of the
theorem, it is enough to apply the dimensional lower bounds of Corollary
1.2 of \cite{os16} to a versal family in which $X_{\eta}$ varies over
an open subset of curves of genus $g$, and $\sL_{\eta}$ varies over an
open subset of special line bundles. Although \textit{loc.\ cit.}\ is
not stated for families, the argument goes through unchanged in this context;
see also Theorem 3.4 of \cite{os19}, where essentially the same argument
is developed in families. Note that here we also use that the locus of
special line bundles is irreducible on a general curve: for $d=2g-2$, this
is because the canonical bundle is the unique choice, while for $d<2g-2$,
this is a consequence of classical Brill-Noether theory (see 
\cite{f-l1} and \cite{gi1}).

Finally, the assertions on (semi)stability follow from openness in
families, which is Proposition 6.4.4 of \cite{os20}.
\end{proof}

\section{Specific existence results}\label{sec:families}

We now explicitly produce families of limit linear series 
having special determinant and (the appropriately modified) expected 
dimension, allowing us to apply Theorem \ref{thm:main} to conclude
Theorem \ref{thm:main-2}.

We will work with the following situation throughout this section.

\begin{sit}
Let $X_1,\dots,X_g$ be genus-$1$ curves, with marked points
$P_i, Q_i$ such that $P_i-Q_i$ is not $m$-torsion for any $m$ with
$1 \leq m \leq 2d$.
Let $X$ be obtained by gluing $Q_i$ to $P_{i+1}$ for 
$i=1,\dots,g-1$. Fix $d$ with $g-2 \leq d \leq 2g-2$. For $i=1,\dots,g$ set 
$$D_i=
\begin{cases} Q_1: & i=1 \\ Q_i+P_i: & 1<i<g \\ P_g: & i=g. \end{cases}.$$
\end{sit}

In this situation, we write 
$((\sE^i,V^i)_{i=1,\dots,g},(\vp_j)_{j=1,\dots,g-1})$ to denote
a limit linear series on $X$, where 
$(\sE^i,V^i)$ is on the component $X_i$, and $\vp_j$ is a gluing map from
$X_j$ to $X_{j+1}$. 

Recall that the canonical line bundle $\sL$ on $X$ is determined by
$\sL|_{X_i} = \sO_{X_i}(D_i)$ for each $i$. 

\begin{prop}\label{prop:special-lb}
Let $\sL$ be a line bundle on $C$ of degree $d$
such that if $i \leq d+2-g$ we have
$\sL|_{X_i} \cong \cO_{X_i}(D_i)$,
and for $i>d+2-g$ we have $\deg \sL|_{X_i} = \deg D_i-1$, and 
$\sL|_{X_i} \not \cong \sO_{X_i}(D)$ for $0 \leq D < D_i$.
Then $\sL$ is special. 

Moreover, the degree-$d$ aspect $\sL_i$
of $\sL$ on $X_i$ 
for $1 \leq i \leq d+2-g$ is equal to $\sO_{X_i}(2(i-1)P_i+(d-2i+2)Q_i)$.
If we fix $d_{\bullet}=(d-1,\dots,d-1,d,d-1,\dots,d-1)$, with $d_j=d$,
then the degree-$d_i$ aspect of $\sL_i$ of $\sL$ on $X_i$ for 
$1 \leq i \leq d+2-g$ is given by
$$\sL_i \cong \begin{cases} \sO_{X_i}(2(i-1)P_i+(d-2i+1)Q_i): & i<j\\
\sO_{X_i}(2(i-1)P_i+(d-2i+2)Q_i): & i=j\\
\sO_{X_i}((2i-3)P_i+(d-2i+2)Q_i): & i>j.\end{cases}$$
\end{prop}

In the above, the ``degree-$d_i$ aspect'' is the line bundle
on $X_i$ obtained from $\sL$ as follows: if $w$ is the multidegree of $\sL$,
$b$ is determined by $\sum_i d_i-(g-1)b=d$,
and $w_i$ is given by $d_i$ in index $i$ and $d_j-b$ for all $j \neq i$,
then the degree-$d_i$ aspect of $\sL$ is the restriction to $X_i$ of
$\sL \otimes \sO'_{w,w_i}$.

\begin{proof} Given the description of the canonical line bundle on $X$,
the speciality is immediate from the definitions. Next, let $c_i$ be the 
degree of $\sL|_{X_i}$ for each $i$; then we have 
$c_i = 2-\delta_i-\epsilon_i$, where 
$\delta_i=\begin{cases} 1 : & i=1,g \\ 0: & \text{otherwise,}\end{cases}$
and 
$\epsilon_i
=\begin{cases} 0 : & i \leq d+2-g \\ 1 : & \text{otherwise.}\end{cases}$
Then for the degree-$d$ aspects, we have
$$\sL_i=\sL|_{X_i}((c_1+\dots+c_{i-1})P_i+(c_{i+1}+\dots+c_g)Q_i),$$
and it is enough to verify that $c_1+\dots+c_{i-1}+1=2(i-1)$ for 
$1<i \leq d+2-g$, which is clear.
The degree-$d_i$ aspects for the variant multidegree are computed similarly.
\end{proof}

In this situation, we can give very explicit criteria for 
$\ell$-(semi)stability as follows.

\begin{prop}\label{prop:chain-stable} Given a limit linear series 
$((\sE^i,V^i)_{i},(\vp_j)_{j})$ 
of degree $d$ on $X$, suppose that each $\sE^i$ is semistable. Then
the limit linear series is \el-semistable. 
If further there do not exist subbundles $\sF^i \subseteq \sE^i$ of equal
rank $r'<r$ such that $\chi(\sF^i)/r'=\chi(\sE^i)/r$ for all $i$, and
each $\sF^i$ glues to $\sF^{i+1}$ under $\vp_i$, then the limit linear
series is \el-stable.
\end{prop}

\begin{proof} Corollary 1.7 of \cite{os22} gives the analogous statements
for a single vector bundle on $X$ and its restrictions to the $X_i$.
Because semistability on components and the gluing maps between components
are both preserved under twists, the desired statement follows.
\end{proof}

The following notation will be convenient.

\begin{notn}\label{not:twist-down}
If $V \subseteq \Gamma(X',\sE)$ for a vector bundle $\sE$
on a smooth projective curve $X'$, and $D$ is an effective divisor on 
$X'$, set
$$V(-D):=V \cap \Gamma(X',\sE(-D)).$$
\end{notn}

As a final preparatory step, we state the following result, which is 
important for analyzing possibilities for limit linear series with
prescribed vanishing sequences, and which is an immediate consequence
of Proposition 5.1.3 of \cite{os20}.

\begin{prop}\label{prop:unique-subspaces} Given a pair $(\sE,V)$ of a 
vector bundle on a smooth projective curve $X'$ together with a 
$k$-dimensional space $V$ of global of $\sE$, and given also distinct
point $P,Q \in X'$, suppose that the vanishing sequence of $(\sE,V)$ 
at $P$ (respectively, $Q$) is given by $a_1,\dots,a_k$ (respectively,
$b_1,\dots,b_k$). Then given nonnegative integers $a,b$, we have
$$\dim V(-aP-bQ) \geq \#\{i:a_i \geq a \text{ and } b_{k+1-i} \geq b\}.$$

In particular, if we have
$$\dim \Gamma(X',\sE(-aP-bQ))
=\#\{i:a_i \geq a \text{ and } b_{k+1-i} \geq b\},$$
then $\Gamma(X',\sE(-aP-bQ)) \subseteq V$.

Additionally, for any $i$ we must have a section in $V$ vanishing to order
at least $a_i$ at $P$ and at least $b_{k+1-i}$ at $Q$.
\end{prop}

We now proceed to the limit linear series arguments which lead to
Theorem \ref{thm:main-2}. We have four cases to consider, according to
the parity of $d$ and of $k$ (thus, only the first two cases were necessary
to consider in the case of canonical determinant). Because the families
described are somewhat complicated, we have included a small example of
each family in Appendix \ref{app:examples}.

\begin{prop}\label{prop:even-even} Assume that $d=2d'$ and $k=2k'$ are even,
and set $d_i=d$ for $i=1,\dots,g$. Then the $b$ of Situation \ref{sit:basic}
is determined to be $d'$. Fix $\sL$ as in Proposition 
\ref{prop:special-lb}, and assume further that for 
$i>d+2-g$, we have $\sL_i \not\cong \sO_{X_i}(aP_i+(d-a)Q_i)$ for any
$a$ between $0$ and $d$. Suppose further that
\begin{equation}\label{eq:even-even-ineq} g\geq (k')^2+2k'(g-1-d').
\end{equation}
Then there exists a nonempty open substack of the 
\el-semistable chain-adaptable locus of
$\cG^{k,\EHT}_{2,\sL,d_{\bullet}}(X)$ having the expected dimension 
$\rho_{\sL}:=\rho-g+\binom{k}{2}$.
If further 
$$(g,d,k) \neq (1,0,2), (2,2,2), (3,2,2) \text{ or } (4,6,4),$$
the same is true of the \el-stable locus.
\end{prop}

\begin{proof} The basic strategy is that we consider one set of conditions
on the first $(k')^2$ components of $X$, different conditions on the next 
$2k'(g-1-d')$
components, and finally we allow generic behavior on the remaining 
$g-(k')^2-2k'(g-1-d')$ components. On the first set of components, we will
impose extra vanishing on two sections per component, with indices cycling 
over increasing odd numbers of components (that is, we break $(k')^2$ into
$1+3+\dots+(2k'+1)$). On the next collection of components, we impose 
extra vanishing on only a single section per component, cycling 
one by one over the $2k'$ sections determining our vector spaces. On
the remaining components, no extra vanishing is imposed. See Example
\ref{ex:example-1}. We begin by
defining sequences $a^i$ for $i=1,\dots,g+1$ as follows:
$$a^1=0,0,1,1,\dots,k'-1,k'-1, \quad\text{and}\quad 
a^{i+1}_j=a^{i}_j+1-\epsilon^i_j,$$
with $\epsilon^i_j=0$ or $1$, and the latter case occurring precisely when
one of the following holds:
\begin{itemize} 
\item ($1 \leq i \leq (k')^2$) we have $i=m^2+2c+1$ for some $0 \leq m < k'$ 
and $0 \leq c < m$, and $j=2c+1$ or $2m+1$;
\item ($1 \leq i \leq (k')^2$) we have $i=m^2+2c+2$ for some $0 \leq m < k'$ 
and $0 \leq c < m$, and $j=2c+2$ or $2m+2$;
\item ($1 \leq i \leq (k')^2$) we have $i=m^2+2m+1$ for $0 \leq m < k'$, and
$j=2m+1$ or $2m+2$;
\item ($(k')^2+1 \leq i \leq (k')^2+2k'(g-1-d')$) we have $i=(k')^2+2k'm+c$
for some $0 \leq m <(g-1-d')$ and $1 \leq c \leq 2k'$, and $j=c$.
\end{itemize}
Our arguments will make use of the following elementary observations:
\begin{ilist}
\item for $i=1,\dots,g+1$ and $j=1,\dots,k'-1$ we have
$$a^i_{2j-1} \leq a^i_{2j} < a^i_{2j+1} \leq a^i_{2j+2};$$
\item for $i=1,\dots,g$ and $j=1,\dots,k'$ we have
$$a^{i+1}_{2j}=a^i_{2j-1}+1$$
unless $i=(m+1)^2$ for $0 \leq m < k'$ and $j=m+1$,
in which case $a^{i+1}_{2j}=a^i_{2j-1}$;
\item for $i=1,\dots,(k')^2$, if $j_1<j_2$ are the indices with
$\epsilon^i_{j_1}=\epsilon^i_{j_2}=1$, we have
$$a^i_{j_1}+a^i_{j_2}=2(i-1);$$
\item for $i=(k')^2+1,\dots,(k')^2+2k'(g-1-d')$, if $j$ is the index with
$\epsilon^i_j=1$, then 
$$2(i-1)>a^i_j+a^i_{2k'}+1;$$
\item we have
$$a^{g+1}=d',d',d'-1,d'-1,\dots,d'-(k'-1),d'-(k'-1).$$
\end{ilist}

We now describe the desired open subset 
$\cU \subseteq \cG^{k,\EHT}_{2,\sL,d_{\bullet}}(X)$ in terms of three 
conditions, as the locus of tuples 
$((\sE^i,V^i)_{i},(\vp_j)_{j},\psi)$ 
satisfying the following:
\begin{Ilist} 
\itm each $\sE^i$ is semistable on $X_i$;
\itm for $i > (k')^2+2k'(g-1-d')$, there does not exist any
line subbundle $\sL'$ of $\sE^i$ with either 
$\sL' \cong \sO_{X_i}(aP_i+(d'-a)Q_i)$ for
$a=0,\dots,d'$, or with $\sL'^2 \cong \sL_i$;
\itm for $i=1,\dots,g-1$ the vanishing sequence of $V_{i+1}$ at $P_{i+1}$ is 
equal to $a^{i+1}$, and the vanishing sequence of $V_{i}$ at
$Q_{i}$ is equal to $d'-a^{i+1}_k,\dots,d'-a^{i+1}_1$. 
\end{Ilist}

Notice that because the imposed vanishing at each $Q_i$ and $P_{i+1}$ has
complementary vanishing orders adding to precisely the required $d'$, to 
impose the desired vanishing is equivalent to imposing that we have 
\textit{at most} the desired vanishing, which is an open condition.
Thus, our description makes it clear that $\cU$ is indeed an open 
substack of $\cG^{k,\EHT}_{2,\sL,d_{\bullet}}(X)$, so
in order to prove the proposition, we need to show that $\cU$
is nonempty of the expected dimension, with the necessary (semi)stability
and chain-adaptability properties. The main point will be to give 
an equivalent, more explicit description of $\cU$. Specifically, we claim 
that a point
$((\sE^i,V^i)_{i},(\vp_j)_{j},\psi)$
of $\cG^{k,\EHT}_{2,\sL,d_{\bullet}}(X)$ is in $\cU$
if and only if it satisfies the following conditions:
\begin{enumerate}
\item[(I$'$)] the $\sE^i$ are described as follows:
\begin{itemize}
\item for each $i=1,\dots,(k')^2$, writing $i=m^2+c'$ for 
$1 \leq c' \leq 2m+1$, we have
$$\sE^i \cong \sO_{X_i}(a^i_{c'} P_i + (d'-a^i_{c'}) Q_i) \oplus 
\sL_i(-a^i_{c'} P_i -(d'-a^i_{c'} )Q_i);$$
\item for each $i=(k')^2+1,\dots,(k')^2+2k'(g-1-d')$, writing
$i=(k')^2+2k'm+c$ for $1 \leq c \leq 2k'$, we have
$$\sE^i \cong \sO_{X_i}(a^i_c P_i + (d'-a^i_c) Q_i) \oplus 
\sL_i(-a^i_c P_i -(d'-a^i_c )Q_i);$$
\item for each $i>(k')^2+2k'(g-1-d')$, we have 
$$\sE^i \cong \sL'_i \oplus \left(\sL_i \otimes (\sL'_i)^{-1}\right),$$
where $\sL'_i$ has degree $d'$, and is chosen so that $(\sL'_i)^{2} \not \cong
\sL_i$, and neither $\sL'_i$ nor $\sL_i \otimes (\sL'_i)^{-1}$
is of the form 
$\sO_{X_i}(aP_i+(d'-a)Q_i)$, for $a=0,\dots,d'$;
\end{itemize}
\item[(II$'$)] for each $i$, we have
$$V^i = \bigoplus_{j=1}^{k'} 
\Gamma(X_i, \sE^i(-a^i_{2j-1} P_i-(d'-a^{i+1}_{2j})Q_i)).$$
\end{enumerate}
More precisely, we will show in particular that if $\sE^i$ is as in (I$'$), 
then the space $V^i$ described in (II$'$) has the correct dimension $k$, and
the vanishing sequences asserted in (III). In
fact, in the course of proving the claim, we will further prove that
$V^i$ has a $(P_i,Q_i)$-adapted basis for all $i$, so that 
$\cU$ consists entirely of chain-adaptable limit linear series. Moreover,
this basis may be chosen to come from the maximal summands of $\sE^i$,
making gluings easy to analyze.

The first observation is that since $k' \geq 1$, we obtain
$d+2-g \geq (k')^2$ from \eqref{eq:even-even-ineq}. Thus, 
for $1 \leq i \leq (k')^2$, we have 
$\det \sE^i \cong \sO_{X_i}(2(i-1)P_i+(d-2i+2)Q_i)$.
On the other hand, for $i>(k')^2$, if $i \leq d+2-g$ we see that
$$\sL_i(-a^i_c P_i -(d'-a^i_c )Q_i)
=\sO_{X_i}((2i-2-a^i_c) P_i+(d'-2i+2+a^i_c)Q_i).$$ 

Now, it is clear that (I$'$) implies conditions (I) and (II), and
conversely, for $i \leq (k')^2$, condition (I$'$) follows
from (I) and (III) by Proposition \ref{prop:unique-subspaces}.
For $(k')^2+1 \leq i \leq (k')^2+2k'(g-1-d')$, we similarly 
obtain from (I) and (III) that $\sO_{X_i}(a^i_c P_i + (d'-a^i_c) Q_i)$ 
must be a line subbundle of $\sE^i$, so it suffices to check that
$\sE^i$ cannot be indecomposable. But this is immediate from the
observation that
$\sL_i \not \cong \sO_{X_i}(2a^i_c P_i + 2(d'-a^i_c)Q_i)$, which is
a consequence of (iv) above for 
$i \leq d+2-g$ and is a hypothesis for $i>d+2-g$.
Finally, for $i>(k')^2+2k'(g-1-d')$,
condition (I$'$) is immediate from conditions (I) and (II).
It thus 
suffices to show that assuming the description in condition (I$'$), we
have that (III) is equivalent to (II$'$).

For each 
$i=1,\dots,g$ and $j=1,\dots,k'$, we first observe that it follows from
(i) above that
$$\#\{\ell:a^i_{\ell} \geq a^i_{2j-1} \text{ and }
 d'-a^{i+1}_{\ell} \geq d'-a^{i+1}_{2j}\}=\#\{2j-1,2j\}=2.$$
At the same time, if $\sE^i$ are as in (I$'$), we see that 
$$\dim \Gamma(X_i,\sE^i(-a^i_{2j-1}P_i-(d'-a^{i+1}_{2j})Q_i))=2,$$
and this space is generated by a pair of a sections $s_1$ and $s_2$, one 
from each summand of $\sE^i$, with $s_1$
vanishing to order $a^i_{2j-1}$ at $P_i$ and $d'-a^{i+1}_{2j-1}$ at $Q_i$,
and $s_2$ vanishing to order $a^i_{2j}$ at $P_i$ and $d'-a^{i+1}_{2j}$
at $Q_i$. Indeed, we see immediately from (ii) and
the fact that $\sE^i$ is a direct sum of two degree-$d'$ line bundles
that 
$$\dim \Gamma(X_i,\sE^i(-a^i_{2j-1}P_i-(d'-a^{i+1}_{2j})Q_i))\leq 2,$$
with equality except possibly when $i=(m+1)^2$ for $0 \leq m < k'$.
But equality in this last case follows from (I$'$) and (iii),
which together with the determinant description
imply that
$\sE^i(-a^i_{2j-1}P_i-(d'-a^{i+1}_{2j})Q_i)\cong \sO_{X_i}^{\oplus 2}$.
It thus remains to verify that we have sections $s_1$ and $s_2$ as 
described. One checks
that for our choice of the $a^i$, this is equivalent to the condition
that $\sE^i$ contain as a summand the bundle
$\bigoplus_{\ell:\epsilon^i_{\ell}=1} 
\sO_{X_i} (a^i_{\ell} P_i +(d'-a^i_{\ell})Q_i)$, with additional
summands of the form $\sO_{X_i}(aP_i+(d'-a)Q_i)$ only if $a$ is
not equal to any $a^i_{2j-1}$ or $a^i_{2j-1}+1$.
This in turn follows for $i \leq (k')^2$ from (iii), and
is trivially satisfied for $i>(k')^2+2k'(g-1-d')$. In the remaining
case, our genericity hypotheses for $\sL$ ensures the condition is
satisfied if also $i>d+2-g$, and if $i \leq d+2-g$, we conclude the
desired statement from (iv).

We thus conclude from Proposition \ref{prop:unique-subspaces}
that if (III) is satisfied, then $V^i$ contains 
$\Gamma(X_i,\sE^i(-a^i_{2j-1}P_i-(d'-a^{i+1}_{2j})Q_i))$.
On the other hand, as $j$ varies, (i) implies that
these spaces are disjoint,
so we conclude that (III) and (I$'$) imply (II$'$). Note that we have also
shown that $V^i$ has a $(P_i,Q_i)$-adapted basis, giving the asserted chain
adaptability.
Conversely, we have seen that the $V^i$ described in (II$'$) has vanishing
sequences as specified by (III),
completing the proof that $\cU$ is described equivalently by (I$'$) and
(II$'$).

We now investigate gluings. For this, it is useful to observe that the
choice of determinant isomorphism $\psi$ is equivalent to a tuple of
determinant isomorphisms $\psi_i$ for the restriction to each $X_i$,
commuting with (the determinants of) the gluings maps $\vp_j$.
Now, we claim that given bundles 
$\sE^i$ and spaces $V^i$ as specified in (I$'$) and (II$'$), together with 
choices of determinant isomorphisms $\psi_i$, there always exist gluings
$\vp_j$, which are unrestricted except for the compatibility with the
$\psi_i$ (which can always be achieved by scaling) and the conditions 
imposed by non-repeated vanishing orders. Indeed, we saw in describing the 
$(P_i,Q_i)$-adapted bases of $V^i$ that any nonrepeated vanishing orders 
are realized by sections in one of the two summands of $\sE^i$, so 
it is enough to check
that when there are nonrepeated vanishing orders they come from summands 
which can be matched consistently with the summands in the next component. 
This is easily verified,
and we obtain the following explicit conditions imposed on gluings:
\begin{itemize}
\item for $i=1,\dots,(k')^2$,
writing $i=m^2+c'$ for $1 \leq c' \leq 2m+1$, if $c'$ is odd and less than
$2m+1$ we must have
$$\vp_{i}\left(\left.
\sO_{X_i}(a^i_{c'} P_i + (d'-a^i_{c'}) Q_i)\right|_{Q_i}\right) 
= \left.
\sL_{i+1}(-a^{i+1}_{c'+1} P_{i+1} -(d'-a^{i+1}_{c'+1})Q_{i+1})
\right|_{P_{i+1}},$$
$$\text{and}$$
$$\vp_i\left(\left.
\sL_i(-a^i_{c'} P_i -(d'-a^i_{c'} )Q_i)\right|_{Q_i}\right) 
= \left.
\sO_{X_{i+1}}(a^{i+1}_{c'+1} P_{i+1} + (d'-a^{i+1}_{c'+1}) Q_{i+1})
\right|_{P_{i+1}};$$
\item for $i=(k')^2+1,\dots,(k')^2+2k'(g-1-d')$, writing
$i=(k')^2+2k'm+c$ for $1 \leq c \leq 2k'$, if $c$ is odd we must have
$$\vp_i\left(\left.
\sO_{X_i}(a^i_c P_i + (d'-a^i_c) Q_i)\right|_{Q_i}\right) 
= \left.
\sL_{i+1}(-a^{i+1}_{c+1} P_{i+1} -(d'-a^{i+1}_{c+1})Q_{i+1})
\right|_{P_{i+1}},$$
$$\text{and}$$
$$\vp_i\left(\left.
\sL_i(-a^i_c P_i -(d'-a^i_c )Q_i)\right|_{Q_i}\right) 
= \left.
\sO_{X_{i+1}}(a^{i+1}_{c+1} P_{i+1} + (d'-a^{i+1}_{c+1}) Q_{i+1})
\right|_{P_{i+1}}.$$
\end{itemize}
It then follows in particular that $\cU$ is nonempty.

It follows immediately from the definition of $\cU$ and 
Proposition \ref{prop:chain-stable} that the limit linear series in 
$\cU$ are all \el-semistable.
Next, observe that due to our genericity hypothesis on the $P_i$ and $Q_i$,
we have that $\sE^i \cong \sL'^{\oplus 2}$ for some $\sL'$ on $X_i$ if
and only if $i=m^2$ for some $m \leq k'$.
Now, if $k'\geq 2$, there are only two weakly destabilizing line subbundles
$\sM_3 \subseteq \sE^3$, which would have to restrict to one of two lines 
in the fiber at $Q_3$, and hence if we fix a gluing $\vp_3$, there are only 
two weakly destabilizing line subbundles $\sM_4 \subseteq \sE^4$ which could 
glue to a choice of $\sM_3$.
Thus if $g>4$, the unrestricted gluing between $\sE^4$ and $\sE^5$
gives us by Proposition \ref{prop:chain-stable} 
that an open dense subset of $\cU$ is \el-stable. 
On the other hand, if $k'=1$ and $g>3$, we find that we have an unrestricted 
gluing between $\sE^3$ and $\sE^4$, which likewise yields \el-stability.
Similarly, if $k'=1$, $g-1-d'=0$, and $g=3$, we obtain \el-stability from 
the unrestricted gluing between $\sE^2$ and $\sE^3$. The only remaining 
cases for which $g \geq (k')^2+2k'(g-1-d')$ are
$$(g,d,k) = (1,0,2), (2,2,2), (3,2,2) \text{ or } (4,6,4),$$
as desired.

It remains to compute the dimension of $\cU$. All the $\sE^i$ are uniquely 
determined except when $i>(k')^2+2k'(g-1-d')$. For the latter, the coarse 
space for the choices of $\sE^i$ is smooth of dimension $1$ for each $i$.
Thus, if we consider the forgetful map to the moduli stack of tuples
$(\sE^i,\psi_i)_{i=1,\dots,g}$ (i.e., the map which forgets gluings and 
subspaces), we see that
the image $\cU'$ of $\cU$ is a gerbe over a space of dimension 
$g-((k')^2+2k'(g-1-d'))$, 
with stabilizer group equal to 
the group of tuples of automorphisms of the $\sE^i$ which preserve the 
$\psi_i$. The condition of preserving the $\psi_i$ reduces the dimension
of each automorphism group by $1$, so taking this into account, each 
$\sE^i$ has a $1$-dimensional group of automorphisms, except for $\sE^{a^2}$
with $a \leq k'$, which each have $3$-dimensional automorphism groups.
We conclude that the stabilizer groups of $\cU'$ have dimension
$g+2k'$, and hence $\cU'$ is smooth of dimension
$$g-((k')^2+2k'(g-1-d'))-(g+2k')=-k'(k'+2g-2d').$$ 
The fibers of the forgetful map are described by the choices of 
gluings and subspaces, but all the $V^i$ are uniquely determined.
The gluings also have to commute with the $\psi_i$, so the
dimensions are also each reduced by $1$. We thus see from the above
description of restrictions on gluing that the dimension of each fiber is 
$$3(g-1)-2\cdot \frac{(k')^2-k'}{2}-2\cdot k'(g-1-d')=3g-3-k'(k'+2g-3-2d'),$$
so we conclude that the dimension of $\cU$ is
\begin{equation*}\begin{split}-k'(k'+2g-2d') +3g-3&-k'(k'+2g-3-2d') \\
& = 3g-3-2k'(2k'+2g-2-2d')+2(k')^2-k' \\
& =\rho_{\sL},\end{split}\end{equation*}
as desired.
\end{proof}

For the next case, it will be convenient to introduce some additional
notation and terminology.

\begin{notn}\label{notn:fiber-lines} Let $Y$ be a smooth, proper curve,
and $(\sE,V)$ a vector bundle of rank $2$ on $Y$ together with a vector 
space of global sections. Given $P \in Y$, let $a_1,\dots,a_k$ be the
vanishing sequence of $(\sE,V)$ at $P$, and suppose that $a_i$ is not
repeated in the sequence. Then $\ell(V,a_i,P)$ denotes the line in
$\sE|_P$ obtained from 
$V(-aP) \subseteq \Gamma(Y,\sE(-aP)) \to \Gamma(Y,\sE(-aP)|_P)$, 
using the canonical indentification $\PP(\sE|_P)=\PP(\sE(-aP)|_P)$.
\end{notn}

\begin{defn}\label{def:special-lines} Let $Y$ be a smooth, proper curve,
$P \in Y$ a point, $\sE$ of rank $2$ on $Y$, and $\ell \in \sE|_P$ a line.
Suppose we have an exact sequence of bundles
$$0 \to \sL \to \sE \to \sL' \to 0$$
with $\deg \sL = \deg \sL'$.
If the sequence is split and $\sL \not \cong \sL'$,
we say $\ell$ is \textbf{nondegenerate} if $\ell \neq \sL|_P,\sL'_P$.
If the sequence is nonsplit, we say $\ell$ is \textbf{nondegenerate} if
$\ell \neq \sL|_P$.
\end{defn}

Note that in either case, the notion of nondegeneracy is intrinsic to $\sE$.

\begin{prop}\label{prop:even-odd} Assume that $d=2d'$ is even and $k=2k'+1$
is odd,
and set $d_i=d$ for $i=1,\dots,g$. Then the $b$ of Situation \ref{sit:basic}
is determined to be $d'$. Fix $\sL$ as in Proposition 
\ref{prop:special-lb}, and assume further that for 
$i>d+2-g$, we have $\sL_i \not\cong \sO_{X_i}(aP_i+(d-a)Q_i)$ for any
$a$ between $0$ and $d$. Suppose further that
\begin{equation}\label{eq:even-odd-ineq} g\geq (k')^2+k'+1+(2k'+1)(g-1-d').
\end{equation}
Then there exists a nonempty open substack of the 
\el-semistable chain-adaptable locus of
$\cG^{k,\EHT}_{2,\sL,d_{\bullet}}(X)$ having the expected dimension 
$\rho_{\sL}:=\rho-g+\binom{k}{2}$.
If further $(g,d,k) \neq (3,4,3)$, the same is true of the \el-stable locus.
\end{prop}

\begin{proof} The basic strategy is the same as in Proposition 
\ref{prop:even-even}, except that after the first $(k')^2$ components,
we impose a slightly different pattern on the next $k'+1$ components
before continuing as in the previous case. This, combined with the 
additional section to consider, makes the family of limit linear series
more complicated to describe. See Example \ref{ex:example-2}.
We define sequences $a^i$ for $i=1,\dots,g+1$ as follows:
$$a^1=0,0,1,1,\dots,k'-1,k'-1,k' \quad\text{and}\quad 
a^{i+1}_j=a^{i}_j+1-\epsilon^i_j,$$
with $\epsilon^i_j=0$ or $1$, and the latter case occurring precisely when
one of the following holds:
\begin{itemize} 
\item ($1 \leq i \leq (k')^2$) we have $i=m^2+2c+1$ for some $0 \leq m < k'$ 
and $0 \leq c < m$, and $j=2c+1$ or $2m+1$;
\item ($1 \leq i \leq (k')^2$) we have $i=m^2+2c+2$ for some $0 \leq m < k'$ 
and $0 \leq c < m$, and $j=2c+2$ or $2m+2$;
\item ($1 \leq i \leq (k')^2$) we have $i=m^2+2m+1$ for $0 \leq m < k'$, and
$j=2m+1$ or $2m+2$;
\item ($(k')^2+1 \leq i \leq (k')^2+k'$) we have $i=(k')^2+c$
for some $1 \leq c \leq k'$, and $j=2c-1$ or $2k'+1$;
\item ($i = (k')^2+k'+1$) we have $j=2k'+1$;
\item ($(k')^2+k'+2 \leq i \leq (k')^2+k'+1+(2k'+1)(g-1-d')$) we have 
$i=(k')^2+k'+1+(2k'+1)m+c$
for some $0 \leq m <(g-1-d')$ and $1 \leq c \leq 2k'+1$, and $j=c$.
\end{itemize}

Let 
$\cU \subseteq \cG^{k,\EHT}_{2,\sL,d_{\bullet}}(X)$ denote the set of
$((\sE^i,V^i)_{i},(\vp_j)_{j},\psi)$ 
satisfying the following:
\begin{Ilist} 
\itm each $\sE^i$ is semistable on $X_i$;
\itm for $i > (k')^2+k'+1+(2k'+1)(g-1-d')$, there does not exist any
line subbundle $\sL'$ of $\sE^i$ with either 
$\sL' \cong \sO_{X_i}(aP_i+(d'-a)Q_i)$ for
$a=0,\dots,d'$, or with $\sL'^2 \cong \sL_i$;
\itm for $i=1,\dots,g-1$ the vanishing sequence of $V_{i+1}$ at $P_{i+1}$ is 
equal to $a^{i+1}$, and the vanishing sequence of $V_{i}$ at
$Q_{i}$ is equal to $d'-a^{i+1}_k,\dots,d'-a^{i+1}_1$; 
\itm for $i=(k')^2+k'+1$, we have
$h^0(X_i,\sE^i(-(a^i_{2k'+1})P_i-(d'-a^i_{2k'+1})Q_i)) \leq 1$.
\itm we have $\ell(V_i,a,P_i)$ and $\ell(V_{i-1},d'-a,Q_{i-1})$ 
nondegenerate in the following cases:
$$i=m^2+2c+2 \text{ for } 0 \leq m < k',
0 \leq c < m, \text{ and } a = a^i_{2k'+1},$$
$$i=(k')^2+c \text{ for } 3 \leq c \leq k'+1,
\text{ and } a = a^i_{1},$$
$$i=(k')^2+k'+1+(2k'+1)m+c \text{ for } 0 \leq m < g-1-d', 
3 \leq c \leq 2k'+1\text{ odd}, \text{ and } a=a^i_1.$$
\end{Ilist}

As in the proof of Proposition \ref{prop:even-even}, we will verify 
that 
$((\sE^i,V^i)_{i},(\vp_j)_{j},\psi)$
in $\cG^{k,\EHT}_{2,\sL,d_{\bullet}}(X)$ is in $\cU$
if and only if it satisfies the following conditions:
\begin{enumerate}
\item[(I$'$)] the $\sE^i$ are described as follows:
\begin{itemize}
\item for each $i=1,\dots,(k')^2$, writing $i=m^2+c'$ for 
$1 \leq c' \leq 2m+1$, we have
$$\sE^i \cong \sO_{X_i}(a^i_{c'} P_i + (d'-a^i_{c'}) Q_i) \oplus 
\sL_i(-a^i_{c'} P_i -(d'-a^i_{c'} )Q_i);$$
\item for each $i=(k')^2+1,\dots,(k')^2+k'$, we have
$$\sE^i \cong \sO_{X_i}(a^i_{2k'+1} P_i + (d'-a^i_{2k'+1}) Q_i) \oplus 
\sL_i(-a^i_{2k'+1} P_i -(d'-a^i_{2k'+1})Q_i);$$
\item for $i=(k')^2+k'+1$, we have $\sE^i$ the unique indecomposable vector
bundle of degree $d$ with $\sO_{X_i}(a^i_{2k'+1}P_i+(d'-a^i_{2k'+1})Q_i)$ 
occurring
as a line subbundle;
\item for each $i=(k')^2+k'+2,\dots,(k')^2+k'+1+(2k'+1)(g-1-d')$, writing
$i=(k')^2+k'+1+(2k'+1)m+c$ for $1 \leq c \leq 2k'+1$, we have
$$\sE^i \cong \sO_{X_i}(a^i_c P_i + (d'-a^i_c) Q_i) \oplus 
\sL_i(-a^i_c P_i -(d'-a^i_c )Q_i);$$
\item for each $i>(k')^2+k'+1+(2k'+1)(g-1-d')$, we have 
$$\sE^i \cong \sL'_i \oplus \left(\sL_i \otimes (\sL'_i)^{-1}\right),$$
where $\sL'_i$ has degree $d'$, and is chosen so that $(\sL'_i)^{2} \not \cong
\sL_i$, and neither $\sL'_i$ nor $\sL_i \otimes (\sL'_i)^{-1}$
is of the form 
$\sO_{X_i}(mP_i+(d'-m)Q_i)$, for $m=0,\dots,d'$;
\end{itemize}
\item[(II$'$)] The $V^i$ are described as follows:
\begin{itemize}
\item for $i \leq (k')^2$, we have
$$V^i = W^i \oplus \bigoplus_{j=1}^{k'} 
\Gamma(X_i, \sE^i(-a^i_{2j-1} P_i-(d'-a^{i+1}_{2j})Q_i)),$$
where $W^i$ is a $1$-dimensional subspace of the $2$-dimensional space
$\Gamma(X_i, \sE^i(-a^i_{2k'+1} P_i-(d'-a^{i+1}_{2k'+1})Q_i))$;
\item for $i=(k')^2+1$, we have
$$V^i = \Gamma(X_i, \sE^i(-a^i_{2k'+1} P_i-(d'-a^{i+1}_{2k'+1})Q_i)) \oplus 
\bigoplus_{j=1}^{k'} 
\Gamma(X_i, \sE^i(-a^i_{2j-1} P_i-(d'-a^{i+1}_{2j})Q_i));$$
\item for $i=(k')^2+c$, with $2 \leq c \leq k'$, we have
\begin{equation*}\begin{split}
V^i = W^i & \oplus \left(\bigoplus_{j=1}^{c-2} 
\Gamma(X_i, \sE^i(-a^i_{2j} P_i-(d'-a^{i+1}_{2j+1})Q_i))\right)\oplus
U^i \\ & \oplus \left(\bigoplus_{j=c+1}^{k'}
\Gamma(X_i, \sE^i(-a^i_{2j-1} P_i-(d'-a^{i+1}_{2j})Q_i))\right) \\
& \oplus \Gamma(X_i, \sE^i(-a^i_{2k'+1} P_i-(d'-a^{i+1}_{2k'+1})Q_i)),
\end{split}\end{equation*}
where $W^i$ is a $1$-dimensional subspace of the $2$-dimensional space
$\Gamma(X_i, \sE^i(-a^i_{1} P_i-(d'-a^{i+1}_{1})Q_i))$,
and $U^i$ is the ($3$-dimensional) sum
$$\Gamma(X_i, \sE^i(-a^i_{2c-2} P_i-(d'-a^{i+1}_{2c-1})Q_i)) +
\Gamma(X_i, \sE^i(-a^i_{2c-1} P_i-(d'-a^{i+1}_{2c})Q_i));$$
\item for $i=(k')^2+k'+1$ or $i=(k')^2+k'+1+(2k'+1)m+c$ with 
$2 \leq c \leq 2k'+1$ and $0 \leq m < g-1-d'$, or 
$i>(k')^2+k'+1+(2k'+1)(g-1-d')$ we have
$$V^i = W^i \oplus \bigoplus_{j=1}^{k'} 
\Gamma(X_i, \sE^i(-a^i_{2j} P_i-(d'-a^{i+1}_{2j+1})Q_i)),$$
where $W^i$ is a $1$-dimensional subspace of
$\Gamma(X_i, \sE^i(-a^i_{1} P_i-(d'-a^{i+1}_{1})Q_i))$;
\item for $i=(k')^2+k'+1+(2k'+1)m+1$ with 
$0 \leq m < g-1-d'$, we have
$$V^i = \Gamma(X_i, \sE^i(-a^i_1 P_i-(d'-a^{i+1}_{1})Q_i)) \oplus 
\bigoplus_{j=1}^{k'} 
\Gamma(X_i, \sE^i(-a^i_{2j} P_i-(d'-a^{i+1}_{2j+1})Q_i));$$
\end{itemize}
\item[(III$'$)] the spaces $W^i$ and $W^{i-1}$ are chosen so that 
$\ell(V^i,a,P_i)$ 
and $\ell(V^{i-1},d'-a,Q_{i-1})$ satisfy the nondegeneracy conditions of
(V).
\end{enumerate}

We again observe that since $2k'+1 \geq 2$, we obtain
$d+2-g \geq (k')^2+k'+1$ from \eqref{eq:even-odd-ineq}. Thus, 
for $1 \leq i \leq (k')^2+k'+1$, we have 
$\det \sE^i \cong \sO_{X_i}(2(i-1)P_i+(d-2i+2)Q_i)$.
The proof of the equivalence of 
(I)-(V) with (I$'$)-(III$'$) proceeds as in Proposition 
\ref{prop:even-even}. The only significant difference is that we use
condition (IV) to impose indecomposability on $\sE^i$ when $i=(k')^2+k'+1$.
In this case, all global section spaces considered in our description of 
$V^i$ involve twisting down by degree $d'-1$, so are $2$-dimensional.
Note also that for 
for $(k')+2 \leq i \leq (k')^2+k'$,
the intersection of 
$\Gamma(X_i, \sE^i(-a^i_{2c-2} P_i-(d'-a^{i+1}_{2c-1})Q_i))$ with 
$\Gamma(X_i, \sE^i(-a^i_{2c-1} P_i-(d'-a^{i+1}_{2c})Q_i))$ is precisely
the $1$-dimensional space 
$\Gamma(X_i, \sE^i(-a^i_{2c-1} P_i-(d'-a^{i+1}_{2c-1})Q_i))$.
We again obtain from the argument that in fact we have chain-adaptability,
and -- with the exceptions of the sections involved in the lines $W^i$ or
the indecomposable case $V^{(k')^2+k'+1}$ -- that bases may be chosen to lie 
in summand line bundles.

We next claim that, as in Proposition \ref{prop:even-even}, given bundles 
$\sE^i$ and spaces $V^i$ as specified in (I$'$)-(III$'$), together with 
choices of determinant isomorphisms $\psi_i$ on each $X_i$, there always 
exist gluings
$\vp_j$, which are unrestricted except for the conditions imposed by
the $\psi_i$ and by
non-repeated vanishing orders. This is largely straightforward and similar
to the proof of Proposition \ref{prop:even-even}, except that
one has to check (as can be accomplished by direct calculation)
that for the indecomposable vector bundle $\sE^i$ with $i=(k')^2+k'+1$, we
necessarily have $\ell(V^i,a^{i}_{2k'},P_i)$ equal to the 
line at $P_i$ obtained by restricting the line subbundle 
$\sO_{X_i}(a^i_{2k'+1}P_i+(d'-a^i_{2k'+1})Q_i)$.
We obtain the following explicit conditions imposed on gluings:
\begin{itemize}
\item for $i=1,\dots,(k')^2$,
writing $i=m^2+c'$ for $1 \leq c' \leq 2m+1$, we must have
$$\vp_i\left(\ell(V^i,d'-a^{i+1}_{2k'+1},Q_i)\right) =
\ell(V^{i+1},a^{i+1}_{2k'+1},P_{i+1})$$
and if $c'$ is odd and less than $2m+1$ we must also have
$$\vp_{i}\left(\left.
\sO_{X_i}(a^i_{c'} P_i + (d'-a^i_{c'}) Q_i)\right|_{Q_i}\right) 
= \left.
\sL_{i+1}(-a^{i+1}_{c'+1} P_{i+1} -(d'-a^{i+1}_{c'+1})Q_{i+1})
\right|_{P_{i+1}},$$
$$\text{and}$$
$$\vp_i\left(\left.
\sL_i(-a^i_{c'} P_i -(d'-a^i_{c'} )Q_i)\right|_{Q_i}\right) 
= \left.
\sO_{X_{i+1}}(a^{i+1}_{c'+1} P_{i+1} + (d'-a^{i+1}_{c'+1}) Q_{i+1})
\right|_{P_{i+1}};$$
\item for $i=(k')^2+1,\dots,(k')^2+k'$, we must have 
$$\vp_i\left(\ell(V^i,d'-a^{i+1}_{1},Q_i)\right) =
\ell(V^{i+1},a^{i+1}_{1},P_{i+1})$$
$$\text{and}$$
$$\vp_{i}\left(\left.
\sO_{X_i}(a^i_{2k'+1} P_i + (d'-a^i_{2k'+1}) Q_i)\right|_{Q_i}\right) 
= \left.
\sO_{X_{i+1}}(a^{i+1}_{2k'+1} P_{i+1} + (d'-a^{i+1}_{2k'+1}) Q_{i+1})
\right|_{P_{i+1}};$$
\item for $i=(k')^2+k'+1$, we have
$$\vp_i\left(\ell(V^i,d'-a^{i+1}_{1},Q_i)\right) =
\left.\sO_{X_{i+1}}(a^{i+1}_{1} P_{i+1} + (d'-a^{i+1}_{1}) Q_{i+1})
\right|_{P_{i+1}};$$
\item for $i=(k')^2+k'+2,\dots,(k')^2+k'+1+(2k'+1)(g-1-d')$,
writing $i=(k')^2+k'+1+(2k'+1)m+c$
for some $0 \leq m <g-1-d'$ and $1 \leq c \leq 2k'+1$, we require
$$\vp_i\left(\ell(V^i,d'-a^{i+1}_{1},Q_i)\right) =
\ell(V^{i+1},a^{i+1}_{1},P_{i+1}),$$
and if $c$ is even we must also have
$$\vp_i\left(\left.
\sO_{X_i}(a^i_c P_i + (d'-a^i_c) Q_i)\right|_{Q_i}\right) 
= \left.
\sL_{i+1}(-a^{i+1}_{c+1} P_{i+1} -(d'-a^{i+1}_{c+1})Q_{i+1})
\right|_{P_{i+1}},$$
$$\text{and}$$
$$\vp_i\left(\left.
\sL_i(-a^i_c P_i -(d'-a^i_c )Q_i)\right|_{Q_i}\right) 
= \left.
\sO_{X_{i+1}}(a^{i+1}_{c+1} P_{i+1} + (d'-a^{i+1}_{c+1}) Q_{i+1})
\right|_{P_{i+1}};$$
\item for $i=(k')^2+k'+1+(2k'+1)(g-1-d')+1,\dots,g-1$, we must have
$$\vp_i\left(\ell(V^i,d'-a^{i+1}_{1},Q_i)\right) =
\ell(V^{i+1},a^{i+1}_{1},P_{i+1}).$$
\end{itemize}
Note that the gluing of $\sE^i$ to $\sE^{i+1}$ in the cases 
$i=(k')^2+1$ and $i=(k')^2+k'+1+(2k'+1)m+1$ is somewhat special in that
$\ell(V^i,d'-a^{i+1}_{1},Q_i)$ coincides with the line
$\left.\sL_i(-a^i_{2k'+1} P_i -(d'-a^i_{2k'+1})Q_i)\right|_{Q_i}$ for
$i=(k')^2+1$ and with the line
$\left.\sO_{X_i}(a^i_c P_i + (d'-a^i_c) Q_i)\right|_{Q_i}$ for
$i=(k')^2+k'+1+(2k'+1)m+1$. In the other cases in which it occurs,
$\ell(V^i,d'-a^{i+1}_{1},Q_i)$ is an independent direction.
It again follows in particular that $\cU$ is nonempty.

As in Proposition \ref{prop:even-even}, we find that the limit linear series 
in $\cU$ are all \el-semistable, and that
$\sE^i \cong \sL'^{\oplus 2}$ for some $\sL'$ on $X_i$ if
and only if $i=a^2$ for some $a \leq k'$. Also as before, we see that
if $k'\geq 2$ (in which case $g>4$) we obtain \el-stability on a dense open
subset of $\cU$ by considering the gluings between $\sE^3$, $\sE^4$
and $\sE^5$.
On the other hand, if $k'=1$ and $g>3$, we find that we have a 
gluing between $\sE^3$ and $\sE^4$ which does not require gluing together
weakly destabilizing line subbundles, likewise yielding \el-stability.
The only remaining case is $(g,d,k)=(3,4,3)$, as desired.

Finally, we compute the dimension of $\cU$ as in the proof of Proposition 
\ref{prop:even-even}. Analogously to the previous case, the $\sE^i$ are 
uniquely determined except when $i>(k')^2+k'+1+(2k'+1)(g-1-d')$, and there 
is a $1$-dimensional choice for each in the latter case. Because an 
indecomposable even-degree vector bundle has automorphism group of dimension 
$2$,
we see that the sum of the dimensions of the automorphisms of $\sE^i$ 
preserving the $\psi_i$ is still $g+2k'$.
Thus, if we consider the forgetful map to the moduli stack of tuples
$(\sE^i,\psi_i)_{i=1,\dots,g}$, we see that
the image $\cU'$ of $\cU$ is a gerbe over a space of dimension
$g-((k')^2+k'+1+(2k'+1)(g-1-d'))$,
with stabilizer group of dimension $g+2k'$, and hence $\cU'$ is 
smooth of dimension
$$g-((k')^2+k'+1+(2k'+1)(g-1-d'))-(g+2k')=-(k')^2-3k'-1-(2k'+1)(g-1-d').$$ 
The fibers of the forgetful map are described by the choices of 
gluings and subspaces. The $V^i$ are always either uniquely determined,
or have a smooth one-parameter family of choices, and the latter occurs
\begin{equation*}\begin{split}
(k')^2+(k'-1)+1+& (2k')(g-1-d')\\
&+(g-((k')^2+k'+1)+(2k'+1)(g-1-d')) \\
& =g-1-(g-1-d')
\end{split}\end{equation*}
times. We see from the above description of allowable gluings that, given
a fixed choice of the $V_i$, the space of gluings is smooth of dimension
\begin{equation*}\begin{split}
3(g-1)-& \left(2\cdot \frac{(k')^2-k'}{2} +(k')^2\right) -2 k' -1 
 - ((2k'+1)+2 \cdot k') (g-1-d') \\
& \quad \quad -(g-1-((k')^2+k'+1+(2k'+1)(g-1-d'))) \\
& \quad \quad =2g-(k')^2-2-2k'(g-1-d'),
\end{split}
\end{equation*}
so we conclude that the forgetful map is smooth of relative dimension
$$3g-(k')^2-3-(2k'+1)(g-1-d'),$$
and thus that $\cU$ is smooth of dimension
\begin{equation*}\begin{split}-(k')^2-3k'-1 &-(2k'+1)(g-1-d') +3g-(k')^2-3 
-(2k'+1)(g-1-d') \\
& = 3g-3-(2k'+1)(2g-2-2d')-2(k')^2-3k'-1 \\
& = 3g-3-(2k'+1)((2k'+1)+2g-2-2d')+2(k')^2+k' \\
& =\rho_{\sL},\end{split}\end{equation*}
as desired.
\end{proof}

\begin{prop}\label{prop:odd-even} Assume that $d=2d'+1$ is odd and $k=2k'$
is even, and set $d_{(k')^2+2k'(g-2-d')+1}=d$, and $d_i=d-1$ for 
$i \neq (k')^2+2k'(g-2-d')+1$.
Then the $b$ of Situation \ref{sit:basic}
is determined to be $d'$. Fix $\sL$ as in Proposition 
\ref{prop:special-lb}, and assume further that for 
$i>d+2-g$, we have $\sL_i \not\cong \sO_{X_i}(aP_i+(d-1-a)Q_i)$ for any
$a$ between $0$ and $d-1$.
Suppose further that
\begin{equation}\label{eq:odd-even-ineq} g\geq (k')^2+1+k'(2g-2-(2d'+1)).
\end{equation}
Then there exists a nonempty open substack of the 
\el-semistable (equivalently, \el-stable) chain-adaptable locus of
$\cG^{k,\EHT}_{2,\sL,d_{\bullet}}(X)$ having the expected dimension 
$\rho_{\sL}:=\rho-g+\binom{k}{2}$.
\end{prop}

\begin{proof} The approach is similar to Proposition \ref{prop:even-even},
except that we consider an odd-degree indecomposable bundle followed by
$k'$ additional specified vector bundles immediately before passing to
generic behavior. See Example \ref{ex:example-3}. We 
define sequences $a^i$ for $i=1,\dots,g+1$ as follows:
$$a^1=0,0,1,1,\dots,k'-1,k'-1, \quad\text{and}\quad 
a^{i+1}_j=a^{i}_j+1-\epsilon^i_j,$$
with $\epsilon^i_j=0$ or $1$, and the latter case occurring precisely when
one of the following holds:
\begin{itemize} 
\item ($1 \leq i \leq (k')^2$) we have $i=m^2+2c+1$ for some $0 \leq m < k'$ 
and $0 \leq c < m$, and $j=2c+1$ or $2m+1$;
\item ($1 \leq i \leq (k')^2$) we have $i=m^2+2c+2$ for some $0 \leq m < k'$ 
and $0 \leq c < m$, and $j=2c+2$ or $2m+2$;
\item ($1 \leq i \leq (k')^2$) we have $i=m^2+2m+1$ for $0 \leq m < k'$, and
$j=2m+1$ or $2m+2$;
\item ($(k')^2+1 \leq i \leq (k')^2+2k'(g-2-d')$) we have $i=(k')^2+2k'm+c$
for some $0 \leq m <(g-2-d')$ and $1 \leq c \leq 2k'$, and $j=c$;
\item ($i=(k')^2+2k'(g-2-d')+1$) we have $j$ odd;
\item ($(k')^2+2k'(g-2-d')+2 \leq i \leq (k')^2+2k'(g-2-d')+1+k'$) we have 
$i=(k')^2+2k'(g-2-d')+1+c$ with $1 \leq c \leq k'$, and $j=2c$.
\end{itemize}

Now, consider the open subset 
$\cU \subseteq \cG^{k,\EHT}_{2,\sL,d_{\bullet}}(X)$ consisting of tuples
$((\sE^i,V^i)_{i},(\vp_j)_{j},\psi)$ 
satisfying the following:
\begin{Ilist} 
\itm each $\sE^i$ is semistable on $X_i$;
\itm for $i > (k')^2+2k'(g-2-d')+k'+1$, there does not exist any
line subbundle $\sL'$ of $\sE^i$ with either 
$\sL' \cong \sO_{X_i}(aP_i+(d'-a)Q_i)$ for
$a=0,\dots,d'$, or with $\sL'^2 \cong \sL_i$;
\itm for $i=1,\dots,g-1$ the vanishing sequence of $V_{i+1}$ at $P_{i+1}$ is 
equal to $a^{i+1}$, and the vanishing sequence of $V_{i}$ at
$Q_{i}$ is equal to $d'-a^{i+1}_k,\dots,d'-a^{i+1}_1$; 
\itm we have
$\ell(V^i,a^i_{2k'},P_i)$ and $\ell(V^i,d'-1-a^i_{2k'},Q_i)$ 
nondegenerate when $i=(k')^2+2k'(g-2-d')+1+c$ for $1 \leq c \leq k'-1$,
and further $\ell(V^i,d'-a^i_1,Q_i)\neq \ell(V^i,d'-1-a^i_{2k'},Q_i)$ 
for $i=(k')^2+2k'(g-2-d')+1$.
\end{Ilist}

We next claim that a point
$((\sE^i,V^i)_{i},(\vp_j)_{j},\psi)$
of $\cG^{k,\EHT}_{2,\sL,d_{\bullet}}(X)$ is in $\cU$
if and only if it satisfies the following conditions:
\begin{enumerate}
\item[(I$'$)] the $\sE^i$ are described as follows:
\begin{itemize}
\item for each $i=1,\dots,(k')^2$, writing $i=m^2+c'$ for 
$1 \leq c' \leq 2m+1$, we have
$$\sE^i \cong \sO_{X_i}(a^i_{c'} P_i + (d'-a^i_{c'}) Q_i) \oplus 
\sL_i(-a^i_{c'} P_i -(d'-a^i_{c'} )Q_i);$$
\item for each $i=(k')^2+1,\dots,(k')^2+2k'(g-2-d')$, writing
$i=(k')^2+2k'm+c$ for $1 \leq c \leq 2k'$, we have
$$\sE^i \cong \sO_{X_i}(a^i_c P_i + (d'-a^i_c) Q_i) \oplus 
\sL_i(-a^i_c P_i -(d'-a^i_c )Q_i);$$
\item for $i=(k')^2+2k'(g-2-d')+1$, we have that
$\sE^i$ is an indecomposable bundle of degree $d$ and determinant
$\sL_i$ on $X_i$;
\item
for each $i=(k')^2+2k'(g-2-d')+2,\dots,(k')^2+2k'(g-2-d')+k'+1$, writing
$i=(k')^2+2k'(g-2-d')+1+c$, we have
$$\sE^i \cong \sO_{X_i}(a^i_{2c} P_i + (d'-a^i_{2c}) Q_i) \oplus 
\sL_i(-a^i_{2c} P_i -(d'-a^i_{2c} )Q_i);$$
\item for each $i>(k')^2+2k'(g-2-d')+k'+1$, we have 
$$\sE^i \cong \sL'_i \oplus \left(\sL_i \otimes (\sL'_i)^{-1}\right),$$
where $\sL'_i$ has degree $d'$, and is chosen so that $(\sL'_i)^{2} \not \cong
\sL_i$, and neither $\sL'_i$ nor $\sL_i \otimes (\sL'_i)^{-1}$
is of the form 
$\sO_{X_i}(aP_i+(d'-a)Q_i)$, for $a=0,\dots,d'$;
\end{itemize}
\item[(II$'$)] the spaces $V^i$ are described as follows:
\begin{itemize}
\item for $i \leq (k')^2+2k'(g-2-d')$ and 
$i > (k')^2+2k'(g-2-d')+k'+1$, we have
$$V^i = \bigoplus_{j=1}^{k'} 
\Gamma(X_i, \sE^i(-a^i_{2j-1} P_i-(d'-a^{i+1}_{2j})Q_i));$$
\item for $i=(k')^2+k'(g-2-d')+1$ we have
$$V^i = \Gamma(X_i, \sE^i(-a^i_{1} P_i-(d'-a^{i+1}_{2k'-1})Q_i))+ W^i,$$
where $W^i$ is a $2$-dimensional subspace of the $3$-dimensional space
$\Gamma(X_i, \sE^i(-a^i_{2k'-1} P_i-(d'-a^{i+1}_{2k'})Q_i))$
which contains
$\Gamma(X_i, \sE^i(-a^i_{2k'-1} P_i-(d'-a^{i+1}_{2k'-1})Q_i))$;
\item for $i=(k')^2+2k'(g-2-d')+1+c$ with $1 \leq c \leq k'$, we have
\begin{equation*}\begin{split}
V^i = & \left(\bigoplus_{j=1}^{c-1} 
\Gamma(X_i, \sE^i(-a^i_{2j-1} P_i-(d'-a^{i+1}_{2j})Q_i))\right)\oplus
U^i 
\\ & \oplus \left(\bigoplus_{j=c+1}^{k'-1}
\Gamma(X_i, \sE^i(-a^i_{2j} P_i-(d'-a^{i+1}_{2j+1})Q_i))\right)
\oplus W^i,
\end{split}\end{equation*}
where $W^i$ is a $1$-dimensional subspace of 
$\Gamma(X_i, \sE^i(-a^i_{2k'} P_i - (d'-a^{i+1}_{2k'})Q_i))$
and $U^i$ is the ($3$-dimensional) sum
$$\Gamma(X_i, \sE^i(-a^i_{2c-1} P_i-(d'-a^{i+1}_{2c})Q_i)) +
\Gamma(X_i, \sE^i(-a^i_{2c} P_i-(d'-a^{i+1}_{2c+1})Q_i));$$
\end{itemize}
\item[(III$'$)] the spaces $W^i$ are chosen to achieve the nondegeneracy
condition of (IV).
\end{enumerate}
This equivalent description of $\cU$ is checked as in the previous 
propositions. In this case, having $k' \geq 1$ implies that
$d+2-g \geq (k')^2+1$, so we have the necessary speciality on the 
first $(k')^2$ components.

We next consider gluings. As in Proposition \ref{prop:even-even}, given 
bundles $\sE^i$ and spaces $V^i$ as specified in (I$'$)-(III$'$), together 
with choices of determinant isomorphisms $\psi_i$, there always exist gluings
$\vp_j$, which are unrestricted except for the conditions imposed by
the $\psi_i$ and by 
non-repeated vanishing orders. This imposes the following conditions:
\begin{itemize}
\item for $i=1,\dots,(k')^2$,
writing $i=m^2+c'$ for $1 \leq c' \leq 2m+1$, if $c'$ is odd and less than
$2m+1$ we must have
$$\vp_{i}\left(\left.
\sO_{X_i}(a^i_{c'} P_i + (d'-a^i_{c'}) Q_i)\right|_{Q_i}\right) 
= \left.
\sL_{i+1}(-a^{i+1}_{c'+1} P_{i+1} -(d'-a^{i+1}_{c'+1})Q_{i+1})
\right|_{P_{i+1}},$$
$$\text{and}$$
$$\vp_i\left(\left.
\sL_i(-a^i_{c'} P_i -(d'-a^i_{c'} )Q_i)\right|_{Q_i}\right) 
= \left.
\sO_{X_{i+1}}(a^{i+1}_{c'+1} P_{i+1} + (d'-a^{i+1}_{c'+1}) Q_{i+1})
\right|_{P_{i+1}};$$
\item for $i=(k')^2+1,\dots,(k')^2+2k'(g-2-d')$, writing
$i=(k')^2+2k'm+c$ for $1 \leq c \leq 2k'$, if $c$ is odd we must have
$$\vp_i\left(\left.
\sO_{X_i}(a^i_c P_i + (d'-a^i_c) Q_i)\right|_{Q_i}\right) 
= \left.
\sL_{i+1}(-a^{i+1}_{c+1} P_{i+1} -(d'-a^{i+1}_{c+1})Q_{i+1})
\right|_{P_{i+1}},$$
$$\text{and}$$
$$\vp_i\left(\left.
\sL_i(-a^i_c P_i -(d'-a^i_c )Q_i)\right|_{Q_i}\right) 
= \left.
\sO_{X_{i+1}}(a^{i+1}_{c+1} P_{i+1} + (d'-a^{i+1}_{c+1}) Q_{i+1})
\right|_{P_{i+1}};$$
\item for $i=(k')^2+2k'(g-2-d')+1$, we must have
$$\vp_i\left(\ell(V^i,d'-a^i_1,Q_i)\right)
= \left.
\sL_{i+1}(-a^{i+1}_{2} P_{i+1} -(d'-a^{i+1}_{2})Q_{i+1})
\right|_{P_{i+1}}$$
$$\text{and}$$
$$\vp_i\left(\ell(V^i,d'-a^{i+1}_{2k'},Q_i)\right)=
\left(\ell(V^{i+1},a^{i+1}_{2k'},P_{i+1})\right);$$
\item for $i=(k')^2+2k'(g-2-d')+1+c$ with $1 \leq c \leq k'-1$, we must have
$$\vp_i\left(\left.
\sL_{i}(-a^{i}_{2c} P_{i} -(d'-a^{i}_{2c})Q_{i})
\right|_{Q_i}\right)
= \left.
\sL_{i+1}(-a^{i+1}_{2(c+1)} P_{i+1} -(d'-a^{i+1}_{2(c+1)})Q_{i+1})
\right|_{P_{i+1}}$$
$$\text{and}$$
$$\vp_i\left(\ell(V^i,d'-a^{i+1}_{2k'},Q_i)\right)=
\left(\ell(V^{i+1},a^{i+1}_{2k'},P_{i+1})\right);$$
\end{itemize}
It then follows in particular that $\cU$ is nonempty.

It follows immediately from the definition of $\cU$ and 
Proposition \ref{prop:chain-stable} that the limit linear series in 
$\cU$ are all \el-semistable, and since $d$ is odd, this is equivalent
to \el-stability.
It thus remains to compute the dimension of $\cU$. All the $\sE^i$ are uniquely 
determined except when $i=(k')^2+2k'(g-2-d')+1$ or 
$i>(k')^2+2k'(g-2-d')+k'+1$. For the former, there are only finitely many
choices of an indecomposable bundle of specified determinant, while for the
latter, the coarse 
space for the choices of $\sE^i$ is smooth of dimension $1$ for each $i$.
Thus, if we consider the forgetful map to the moduli stack of tuples
$(\sE^i,\psi_i)_{i=1,\dots,g}$ (i.e., the map which forgets gluings and 
subspaces), we see that
the image $\cU'$ of $\cU$ is a gerbe over a space of dimension 
$g-((k')^2+2k'(g-2-d')+k'+1)$.
The stabilizer group is just as in Proposition \ref{prop:even-even},
except that an odd-degree indecomposable vector bundle, being stable,
has automorphism group of dimension $1$,
and hence determinant-fixing automorphism group of dimension $0$,
so the stabilizer groups of $\cU'$ have dimension
$g+2k'-1$, and hence $\cU'$ is smooth of dimension
$$g-((k')^2+2k'(g-2-d')+k'+1)-(g+2k'-1)=-k'(k'+2g-1-2d').$$ 
The $V^i$ are each either uniquely determined, or have a smooth one-parameter
space of choices, with the latter occuring $k'$ times.
The choices of gluings have dimension
$$3(g-1)-2\cdot \frac{(k')^2-k'}{2}-2\cdot k'(g-2-d')-2-2(k'-1)
=3g-3-k'(k'+2g-3-2d'),$$
so the dimension of each fiber of the forgetful map is 
$$3g-3-k'(k'+2g-3-2d')+k'=3g-3-k'(k'+2g-4-2d'),$$
and we conclude that the dimension of $\cU$ is
\begin{equation*}\begin{split}-k'(k'+2g-1-2d') +3g-3&-k'(k'+2g-4-2d') \\
& = 3g-3-2k'(2k'+2g-2-2d'-1)+2(k')^2-k' \\
& =\rho_{\sL},\end{split}\end{equation*}
as desired.
\end{proof}

\begin{prop}\label{prop:odd-odd} Assume that $d=2d'+1$ and $k=2k'+1$ are 
odd, and set $d_{(k')^2+2k'+2}=d$, and $d_i=d-1$ for $i \neq (k')^2+2k'+2$.
Then the $b$ of Situation \ref{sit:basic}
is determined to be $d'$. Fix $\sL$ as in Proposition 
\ref{prop:special-lb}, and assume further that for 
$i>d+2-g$, we have $\sL_i \not\cong \sO_{X_i}(aP_i+(d-1-a)Q_i)$ for any
$a$ between $0$ and $d-1$. Suppose further that
\begin{equation}\label{eq:odd-odd-ineq} g\geq (k')^2+1+(2k'+1)(g-1-d').
\end{equation}
Then there exists a open substack of the 
\el-semistable (equivalently, \el-stable) chain-adaptable locus of
$\cG^{k,\EHT}_{2,\sL,d_{\bullet}}(X)$ having the expected dimension 
$\rho_{\sL}:=\rho-g+\binom{k}{2}$.
\end{prop}

\begin{proof} The construction in this case is a combination of those
of Propositions \ref{prop:even-odd} and \ref{prop:odd-even}. The
first $(k')^2+k'$ components are just as in Proposition \ref{prop:even-odd},
but then the next $k'$ components have specified decomposable bundles
with extra vanishing imposed on a single section per component. On the next
two components, we take an even-degree indecomposable bundle as in
Proposition \ref{prop:even-odd}, and an odd-degree indecomposable bundle 
as in Proposition
\ref{prop:odd-even}, and then the next $(2k'+1)(g-d-2)$ components are as
in the case of Proposition \ref{prop:even-odd}, before reverting to generic
behavior for any remaining components. See Example \ref{ex:example-4}. We 
define sequences $a^i$ for $i=1,\dots,g+1$ as follows:
$$a^1=0,0,1,1,\dots,k'-1,k'-1,k' \quad\text{and}\quad 
a^{i+1}_j=a^{i}_j+1-\epsilon^i_j,$$
with $\epsilon^i_j=0$ or $1$, and the latter case occurring precisely when
one of the following holds:
\begin{itemize} 
\item ($1 \leq i \leq (k')^2$) we have $i=m^2+2c+1$ for some $0 \leq m < k'$ 
and $0 \leq c < m$, and $j=2c+1$ or $2m+1$;
\item ($1 \leq i \leq (k')^2$) we have $i=m^2+2c+2$ for some $0 \leq m < k'$ 
and $0 \leq c < m$, and $j=2c+2$ or $2m+2$;
\item ($1 \leq i \leq (k')^2$) we have $i=m^2+2m+1$ for $0 \leq m < k'$, and
$j=2m+1$ or $2m+2$;
\item ($(k')^2+1 \leq i \leq (k')^2+k'$) we have $i=(k')^2+c$
for some $1 \leq c \leq k'$, and $j=2c-1$ or $2k'+1$;
\item ($(k')^2+k'+1 \leq i \leq (k')^2+2k'$) we have $i=(k')^2+k'+c$
for some $1 \leq c \leq k'$, and $j=2c$;
\item ($i = (k')^2+2k'+1$) we have $j=2k'+1$;
\item ($i = (k')^2+2k'+2$) we have $j$ odd;
\item ($(k')^2+2k'+3 \leq i \leq (k')^2+2k'+2+(2k'+1)(g-2-d')$) we have 
$i=(k')^2+2k'+2+(2k'+1)m+c$
for some $0 \leq m <(g-2-d')$ and $1 \leq c \leq (2k'+1)$, and $j=c$.
\end{itemize}

Now, consider the open subset 
$\cU \subseteq \cG^{k,\EHT}_{2,\sL,d_{\bullet}}(X)$ consisting of tuples
$((\sE^i,V^i)_{i},(\vp_j)_{j},\psi)$ 
satisfying the following:
\begin{Ilist} 
\itm each $\sE^i$ is semistable on $X_i$;
\itm for $i > (k')^2+2k'+2+(2k'+1)(g-2-d')$, there does not exist any
line subbundle $\sL'$ of $\sE^i$ with either 
$\sL' \cong \sO_{X_i}(aP_i+(d'-a)Q_i)$ for
$a=0,\dots,d'$, or with $\sL'^2 \cong \sL_i$;
\itm for $i=1,\dots,g-1$ the vanishing sequence of $V_{i+1}$ at $P_{i+1}$ is 
equal to $a^{i+1}$, and the vanishing sequence of $V_{i}$ at
$Q_{i}$ is equal to $d'-a^{i+1}_k,\dots,d'-a^{i+1}_1$; 
\itm for $i=(k')^2+2k'+1$, we have
$h^0(X_i,\sE^i(-(a^i_{2k'+1})P_i-(d'-a^i_{2k'+1})Q_i)) \leq 1$;
\itm the spaces $V^i$ and $V^{i+1}$ are chosen so that 
$\ell(V^i,d'-a,Q_i)$ 
and $\ell(V^{i+1},a,Q_{i+1})$ are nondegenerate in the following cases:
$$i=m^2+2c+1 \text{ for } 0 \leq m < k',
0 \leq c < m, \text{ and } a = a^{i+1}_{2k'+1},$$
$$i=(k')^2+c \text{ for } 2 \leq c \leq k',
\text{ and } a = a^{i+1}_{1},$$
$$i=(k')^2+k'+c \text{ for } 1 \leq c \leq k'-2,
\text{ and } a = a^{i+1}_{2k'+1},$$
$$i=(k')^2+2k'+2+(2k'+1)m+c \text{ for } 0 \leq m <(g-2-d'), 
2 \leq c \leq 2k'\text{ even}, \text{ and } a=a^{i+1}_1.$$
Also, the spaces $V^i$ are chosen so that 
$\ell(V^i,d'-a,Q_i)$ is nondegenerate for:
$$i=(k')^2+k' \text{ and } a=a^{i+1}_1,$$
$$i=(k')^2+2k'-1 \text{ and } a=a^{i+1}_{2k'+1},$$
and so that $\ell(V^i,a,Q_{i})$ is nondegenerate for:
$$i=(k')^2+k'+1 \text{ and } a = a^i_{2k'+1}.$$
\end{Ilist}

We next claim that a point
$((\sE^i,V^i)_{i},(\vp_j)_{j},\psi)$
of $\cG^{k,\EHT}_{2,\sL,d_{\bullet}}(X)$ is in $\cU$
if and only if it satisfies the following conditions:
\begin{enumerate}
\item[(I$'$)] the $\sE^i$ are described as follows:
\begin{itemize}
\item for each $i=1,\dots,(k')^2$, writing $i=m^2+c'$ for 
$1 \leq c' \leq 2m+1$, we have
$$\sE^i \cong \sO_{X_i}(a^i_{c'} P_i + (d'-a^i_{c'}) Q_i) \oplus 
\sL_i(-a^i_{c'} P_i -(d'-a^i_{c'} )Q_i);$$
\item for each $i=(k')^2+1,\dots,(k')^2+k'$, we have
$$\sE^i \cong \sO_{X_i}(a^i_{2k'+1} P_i + (d'-a^i_{2k'+1}) Q_i) \oplus 
\sL_i(-a^i_{2k'+1} P_i -(d'-a^i_{2k'+1})Q_i);$$
\item for each $i=(k')^2+k'+1,\dots,(k')^2+2k'$, writing
$i=(k')^2+k'+c$, we have
$$\sE^i \cong \sO_{X_i}(a^i_{2c} P_i + (d'-a^i_{2c}) Q_i) \oplus 
\sL_i(-a^i_{2c} P_i -(d'-a^i_{2c} )Q_i);$$
\item for $i=(k')^2+2k'+1$, we have that $\sE^i$ is the unique 
indecomposable bundle of degree $d-1$ containing 
$$\sO_{X_i}(a^i_{2k'+1} P_i + (d'-a^i_{2k'+1}) Q_i)$$
as a line subbundle;
\item for $i=(k')^2+2k'+2$, we have that
$\sE^i$ is an indecomposable bundle of degree $d$ and determinant
$\sL_i$ on $X_i$;
\item for each $i=(k')^2+2k'+3,\dots,(k')^2+2k'+2+(2k'+1)(g-2-d')$, writing
$i=(k')^2+2k'+2+(2k'+1)m+c$ for $1 \leq c \leq 2k'+1$, we have
$$\sE^i \cong \sO_{X_i}(a^i_c P_i + (d'-a^i_c) Q_i) \oplus 
\sL_i(-a^i_c P_i -(d'-a^i_c )Q_i);$$
\item for each $i>(k')^2+2k'+2+(2k'+1)(g-2-d')$, we have 
$$\sE^i \cong \sL'_i \oplus \left(\sL_i \otimes (\sL'_i)^{-1}\right),$$
where $\sL'_i$ has degree $d'$, and is chosen so that $(\sL'_i)^{2} \not \cong
\sL_i$, and neither $\sL'_i$ nor $\sL_i \otimes (\sL'_i)^{-1}$
is of the form 
$\sO_{X_i}(aP_i+(d'-a)Q_i)$, for $a=0,\dots,d'$;
\end{itemize}
\item[(II$'$)] The spaces $V^i$ are described as follows:
\begin{itemize}
\item for $i \leq (k')^2$ or $i=(k')^2+2k'$, we have
$$V^i = W^i \oplus \bigoplus_{j=1}^{k'} 
\Gamma(X_i, \sE^i(-a^i_{2j-1} P_i-(d'-a^{i+1}_{2j})Q_i)),$$
where $W^i$ is a $1$-dimensional subspace of the $2$-dimensional space
$\Gamma(X_i, \sE^i(-a^i_{2k'+1} P_i-(d'-a^{i+1}_{2k'+1})Q_i))$;
\item for $i=(k')^2+1$ or $i=(k')^2+2k'+1$, we have
$$V^i = \Gamma(X_i, \sE^i(-a^i_{2k'+1} P_i-(d'-a^{i+1}_{2k'+1})Q_i)) \oplus 
\bigoplus_{j=1}^{k'} 
\Gamma(X_i, \sE^i(-a^i_{2j-1} P_i-(d'-a^{i+1}_{2j})Q_i));$$
\item for $i=(k')^2+c$, with $2 \leq c \leq k'$, we have
\begin{equation*}\begin{split}
V^i = W^i & \oplus \left(\bigoplus_{j=1}^{c-2} 
\Gamma(X_i, \sE^i(-a^i_{2j} P_i-(d'-a^{i+1}_{2j+1})Q_i))\right)\oplus
U^i 
\\ & \oplus \left(\bigoplus_{j=c+1}^{k'}
\Gamma(X_i, \sE^i(-a^i_{2j-1} P_i-(d'-a^{i+1}_{2j})Q_i))\right) \\
& \oplus \Gamma(X_i, \sE^i(-a^i_{2k'+1} P_i-(d'-a^{i+1}_{2k'+1})Q_i)),
\end{split}\end{equation*}
where $W^i$ is a $1$-dimensional subspace of the $2$-dimensional space
$\Gamma(X_i, \sE^i(-a^i_{1} P_i-(d'-a^{i+1}_{1})Q_i))$,
and $U^i$ is the ($3$-dimensional) sum
$$\Gamma(X_i, \sE^i(-a^i_{2c-2} P_i-(d'-a^{i+1}_{2c-1})Q_i)) +
\Gamma(X_i, \sE^i(-a^i_{2c-1} P_i-(d'-a^{i+1}_{2c})Q_i));$$
\item for $i=(k')^2+k'+c$, with $1 \leq c \leq k'-1$, we have
\begin{equation*}\begin{split}
V^i & = \left(\bigoplus_{j=1}^{c-1} 
\Gamma(X_i, \sE^i(-a^i_{2j-1} P_i-(d'-a^{i+1}_{2j})Q_i))\right)\oplus
U^i 
\\ & \oplus \left(\bigoplus_{j=c+2}^{k'}
\Gamma(X_i, \sE^i(-a^i_{2j-2} P_i-(d'-a^{i+1}_{2j-1})Q_i))\right)
\oplus W^i_1 \oplus W^i_2,
\end{split}\end{equation*}
where $U^i$ is the ($3$-dimensional) sum
$$\Gamma(X_i, \sE^i(-a^i_{2c-1} P_i-(d'-a^{i+1}_{2c})Q_i)) +
\Gamma(X_i, \sE^i(-a^i_{2c} P_i-(d'-a^{i+1}_{2c+1})Q_i)),$$
and $W^i_1$ and $W^i_2$ are $1$-dimensional subspaces of the 
$2$-dimensional spaces
$\Gamma(X_i, \sE^i(-a^i_{2k'} P_i-(d'-a^{i+1}_{2k'})Q_i))$ and
$\Gamma(X_i, \sE^i(-a^i_{2k'+1} P_i-(d'-a^{i+1}_{2k'+1})Q_i))$
respectively, chosen to satisfy the condition that 
$\ell(V^i,a^i_{2k'},P_i)=\ell(V^i,a^i_{2k'+1},P_i)$;
\item for $i=(k')^2+2k'+2$ we have
$$V^i = \Gamma(X_i, \sE^i(-a^i_{1} P_i-(d'-a^{i+1}_{2k'+1})Q_i));$$
\item for $i=(k')^2+2k'+2+(2k'+1)m+1$ with 
$0 \leq m < g-2-d'$, we have
$$V^i = \Gamma(X_i, \sE^i(-a^i_1 P_i-(d'-a^{i+1}_{1})Q_i)) \oplus 
\bigoplus_{j=1}^{k'} 
\Gamma(X_i, \sE^i(-a^i_{2j} P_i-(d'-a^{i+1}_{2j+1})Q_i));$$
\item for $i=(k')^2+2k'+2+(2k'+1)m+c$ with $2 \leq c \leq 2k'+1$
and $0 \leq m < g-2-d'$, or $i>(k')^2+1+(2k'+1)(g-1-d')$ we have
$$V^i = W^i \oplus \bigoplus_{j=1}^{k'} 
\Gamma(X_i, \sE^i(-a^i_{2j} P_i-(d'-a^{i+1}_{2j+1})Q_i)),$$
where $W^i$ is a $1$-dimensional subspace of
$\Gamma(X_i, \sE^i(-a^i_{1} P_i-(d'-a^{i+1}_{1})Q_i))$;
\end{itemize}
\item[(III$'$)] the spaces $W^i$ are chosen to achieve the same nondegeneracy
conditions as in (V).
\end{enumerate}
This equivalent description of $\cU$ is mostly checked as in the previous 
propositions. The key new points involve the fact that 
\eqref{eq:odd-odd-ineq} implies that $d-g+2 \geq (k')^2+2k'+1$,
so we have 
$\sL_i \cong \sO_{X_i}(2(i-1)P_i+2(d'+1-i)Q_i)$ for $i \leq (k')^2+2k'+1$.
This is used for $i=(k')^2+2k'+1$, where we have $a^i_{2k'+1}=i-1$, so 
(given condition (III)) we have that condition (IV) is
equivalent to indecomposability of $\sE^i$. Moreover, we verify that
for $i=(k')^2+k'$ we necessarily have 
$$\ell(V^i,d'-a^{i+1}_{2k'}, Q_i)= \ell(V^i,d'-a^{i+1}_{2k'+1}, Q_i),$$
and it also follows from the description of $\sL_i$ that for 
$i=(k')^2+2k'$ we necessarily have 
$$\ell(V^i,a^{i}_{2k'}, P_i)= \ell(V^i,a^{i}_{2k'+1}, P_i).$$
For $(k')^2+k'<i<(k')^2+2k'$, Lemma 4.2 of \cite{zh2}
asserts that because we have $\sE^i$ of the form
$$\sO_{X_i}(aP_i+(d'-a)Q_i)\oplus \sO_{X_i}((2i-2-a)P_i+(d'-2i+2+a)Q_i)$$
for $a$ satisfying 
$$a<a^i_{2k'}<a^i_{2k'+1}<(2i-2-a)-1,$$
we have $\ell(V^i,a^{i}_{2k'}, P_i)= \ell(V^i,a^{i}_{2k'+1}, P_i)$
if and only if 
$\ell(V^i,d'-a^{i+1}_{2k'}, Q_i)= \ell(V^i,d'-a^{i+1}_{2k'+1}, Q_i)$.
The desired description of $\cU$ then follows.

We next consider gluings. As in Proposition \ref{prop:even-even}, given 
bundles $\sE^i$ and spaces $V^i$ as specified in (I$'$)-(III$'$), together 
with choices of determinant isomorphisms $\psi_i$, there always exist gluings
$\vp_j$, which are unrestricted except for the conditions imposed by
the $\psi_i$ and by
non-repeated vanishing orders. This imposes the following conditions:
\begin{itemize}
\item for $i=1,\dots,(k')^2$,
writing $i=m^2+c'$ for $1 \leq c' \leq 2m+1$, we must have
$$\vp_i\left(\ell(V^i,d'-a^{i+1}_{2k'+1},Q_i)\right) =
\ell(V^{i+1},a^{i+1}_{2k'+1},P_{i+1})$$
and if $c'$ is odd and less than $2m+1$ we must also have
$$\vp_{i}\left(\left.
\sO_{X_i}(a^i_{c'} P_i + (d'-a^i_{c'}) Q_i)\right|_{Q_i}\right) 
= \left.
\sL_{i+1}(-a^{i+1}_{c'+1} P_{i+1} -(d'-a^{i+1}_{c'+1})Q_{i+1})
\right|_{P_{i+1}},$$
$$\text{and}$$
$$\vp_i\left(\left.
\sL_i(-a^i_{c'} P_i -(d'-a^i_{c'} )Q_i)\right|_{Q_i}\right) 
= \left.
\sO_{X_{i+1}}(a^{i+1}_{c'+1} P_{i+1} + (d'-a^{i+1}_{c'+1}) Q_{i+1})
\right|_{P_{i+1}};$$
\item for $i=(k')^2+1,\dots,(k')^2+k'-1$, we must have 
$$\vp_i\left(\ell(V^i,d'-a^{i+1}_{1},Q_i)\right) =
\ell(V^{i+1},a^{i+1}_{1},P_{i+1})$$
$$\text{and}$$
$$\vp_{i}\left(\left.
\sO_{X_i}(a^i_{2k'+1} P_i + (d'-a^i_{2k'+1}) Q_i)\right|_{Q_i}\right) 
= \left.
\sO_{X_{i+1}}(a^{i+1}_{2k'+1} P_{i+1} + (d'-a^{i+1}_{2k'+1}) Q_{i+1})
\right|_{P_{i+1}};$$
\item for $i=(k')^2+k'$, we must have 
$$\vp_i\left(\ell(V^i,d'-a^{i+1}_{1},Q_i)\right) =
\left(\left.\sL_{i+1}
(-a^{i+1}_{2} P_{i+1} -(d'-a^{i+1}_{2} )Q_{i+1})\right|_{P_{i+1}}\right)$$
$$\text{and}$$
$$\vp_{i}\left(\left.
\sO_{X_i}(a^i_{2k'+1} P_i + (d'-a^i_{2k'+1}) Q_i)\right|_{Q_i}\right) 
= \ell(V^{i+1},a^{i+1}_{2k'+1},P_{i+1});$$
\item for $i=(k')^2+k'+c$, with $1 \leq c \leq k'-1$ we must have 
$$\vp_{i}\left(\left.
\sL_i(-a^i_{2c} P_i - (d'-a^i_{2c}) Q_i)\right|_{Q_i}\right) 
= \left(\left.\sL_{i+1}
(-a^{i+1}_{2c+2} P_{i+1} -(d'-a^{i+1}_{2c+2} )Q_{i+1})\right|_{P_{i+1}}\right)$$
$$\text{and}$$
$$\vp_i\left(\ell(V^i,d'-a^{i+1}_{2k'+1},Q_i)\right) =
\ell(V^{i+1},a^{i+1}_{2k'+1},P_{i+1});$$
\item for $i=(k')^2+2k'$, we must have 
$$\vp_i\left(\ell(V^i,d'-a^{i+1}_{2k'+1},Q_i)\right) =
\left(\left.\sO_{X_{i+1}}(a^{i+1}_{2k'+1} P_{i+1} 
+ (d'-a^{i+1}_{2k'+1} )Q_{i+1})\right|_{P_{i+1}}\right);$$
\item for $i=(k')^2+2k'+1$, we must have 
$$\vp_{i}\left(\left.
\sO_{X_i}(a^i_{2k'+1} P_i + (d'-a^i_{2k'+1}) Q_i)\right|_{Q_i}\right) 
= \ell(V^{i+1},a^{i+1}_{2k'+1},P_{i+1});$$
\item for $i=(k')^2+2k'+2$ or $i>(k')^2+1+(2k'+1)(g-1-d')$, we must have 
$$\vp_{i}\left(\ell(V^{i},d'-a^{i+1}_{1},Q_{i})\right)
= \ell(V^{i+1},a^{i+1}_{1},P_{i+1});$$
\item for $i=(k')^2+2k'+3,\dots,(k')^2+2k'+2+(2k'+1)(g-2-d')$,
writing $i=(k')^2+2k'+2+(2k'+1)m+c$
for some $0 \leq m <g-2-d'$ and $1 \leq c \leq 2k'+1$, we require
$$\vp_i\left(\ell(V^i,d'-a^{i+1}_{1},Q_i)\right) =
\ell(V^{i+1},a^{i+1}_{1},P_{i+1}),$$
and if $c$ is even we must also have
$$\vp_i\left(\left.
\sO_{X_i}(a^i_c P_i + (d'-a^i_c) Q_i)\right|_{Q_i}\right) 
= \left.
\sL_{i+1}(-a^{i+1}_{c+1} P_{i+1} -(d'-a^{i+1}_{c+1})Q_{i+1})
\right|_{P_{i+1}},$$
$$\text{and}$$
$$\vp_i\left(\left.
\sL_i(-a^i_c P_i -(d'-a^i_c )Q_i)\right|_{Q_i}\right) 
= \left.
\sO_{X_{i+1}}(a^{i+1}_{c+1} P_{i+1} + (d'-a^{i+1}_{c+1}) Q_{i+1})
\right|_{P_{i+1}};$$
\end{itemize}
It then follows in particular that $\cU$ is nonempty.

It follows immediately from the definition of $\cU$ and 
Proposition \ref{prop:chain-stable} that the limit linear series in 
$\cU$ are all \el-semistable, and since $d$ is odd, this is equivalent
to \el-stability.
It thus remains to compute the dimension of $\cU$. All the $\sE^i$ are uniquely 
determined except when $i=(k')^2+2k'+2$ or $i>(k')^2+(2k'+1)(g-1-d')+1$. For 
the former, there are only finitely many choices for an indecomposable bundle
with specified determinant, while for the latter, the coarse 
space for the choices of $\sE^i$ is smooth of dimension $1$ for each $i$.
Thus, if we consider the forgetful map to the moduli stack of tuples
$(\sE^i,\psi_i)_{i=1,\dots,g}$ (i.e., the map which forgets gluings and 
subspaces), we see that
the image $\cU'$ of $\cU$ is a gerbe over a space of dimension 
$g-((k')^2+(2k'+1)(g-1-d')+1)$.
The stabilizer group is just as in Proposition \ref{prop:odd-even},
so the stabilizer groups of $\cU'$ have dimension
$g+2k'-1$, and hence $\cU'$ is smooth of dimension
$$g-((k')^2+(2k'+1)(g-1-d')+1)-(g+2k'-1)=-(k')^2-2k'-(2k'+1)(g-1-d').$$ 
The $V^i$ are each either uniquely determined, or have smooth one-parameter
spaces of choices. The latter occurs
$$(k')^2+(k'-1)+(k'-1)+1+(2k')(g-2-d')+(g-((k')^2+1+(2k'+1)(g-1-d')))
=d'-1$$
times.
The choices of gluings have dimension
\begin{multline*} 3(g-1)-
\left(2\cdot \frac{(k')^2-k'}{2}+(k')^2\right)
-2(k'-1) -2-2(k'-1) \\
-1-1-1 -(2k'+1)(g-2-d')-2\cdot k'(g-2-d') \\
-\left(g-1-((k')^2+1+(2k'+1)(g-1-d'))\right) \\
=3g-3-(k')^2-(2k'+1)(g-2-d')-d'-k',
\end{multline*}
so the dimension of each fiber of the forgetful map is 
\begin{multline*} 3g-3-(k')^2-(2k'+1)(g-2-d')-d'-k'+d'-1 \\
=3g-3-(k')^2-(2k'+1)(g-2-d')-k'-1,\end{multline*}
and we conclude that the dimension of $\cU$ is
\begin{equation*}\begin{split}
-(k')^2-2k'& -(2k'+1)(g-1-d')
+ 3g-3-(k')^2-(2k'+1)(g-2-d')-k'-1\\
& = 3g-3-(2k'+1)(2g-2-2d'-1)-2(k')^2-3k'-1 \\
& =\rho_{\sL},\end{split}\end{equation*}
as desired.
\end{proof}

We now easily complete the proof of our second main theorem.

\begin{proof}[Proof of Theorem \ref{thm:main-2}] In light of 
Proposition \ref{prop:special-lb}, this follows almost immediately from
Theorem \ref{thm:main} and Propositions \ref{prop:even-even},
\ref{prop:even-odd}, \ref{prop:odd-even} and \ref{prop:odd-odd}.
Indeed, the inequality of the theorem statement in the even 
(respectively, odd) degree case is the same as that in
Proposition \ref{prop:even-even} (respectively, Proposition 
\ref{prop:odd-even}), and the equivalence of the inequality of Proposition
\ref{prop:even-odd} (respectively, Proposition \ref{prop:odd-odd}) follows
from the integrality of the quantities involved.
For the stability assertion, the only remaining case is $(g,d,k)=(3,4,3)$,
which can be proved by direct analysis on smooth curves -- see \S 3 of 
\cite{b-f2} or page 105 of \cite{te5}.
\end{proof}

\section{Further discussion}

We begin by discussing the ranges of $g,k,d$ for which Theorem
\ref{thm:main-2} and Corollary \ref{cor:varying-det} apply.

\begin{rem}\label{rem:ranges} 
In some sense, the range $g-2 \leq d \leq 2g-2$ is misleading, in
that in fact for any given $g$, the degree cannot be too close to $g-2$
without violating the main inequalities of Theorem \ref{thm:main-2}.
For instance, if $d=g-2$, these inequalities are never satisfied with
$k \geq 2$, and if $d=g-1$, the only case which occurs is $k=2$ and
$g$ odd. Despite these limitations, if we set $m=-\chi(\sL)=2g-2-d$,
then for any fixed value of $m \geq 0$, for sufficiently large $g$ we obtain 
(increasingly large) ranges of $k$ for which both Theorem \ref{thm:main-2} 
and Corollary \ref{cor:varying-det} apply, and similarly, for sufficiently 
large $k$ we obtain (increasingly large) ranges of $g$ for which both 
results apply.
In particular, we have produced a large infinite family of examples
of components of the stable locus of $G^k_{2,d}(X)$ having strictly larger
than the expected dimension $\rho-1$.

We also mention that, while our existence results are certainly not optimal, 
and are extended by Zhang in \cite{zh2} using our smoothing theorem,
unlike the canonical determinant case one should not expect the existence 
results of Theorem \ref{thm:main-2} to extend to all cases for which 
$\rho_{\sL} \geq 0$. Indeed, when $d=g-2$ and $k=2$, we have 
$\rho_{\sL} = 
g-6 \geq 0$ 
for $g \geq 6$, but in this case, the locus of stable $\fg^2_{2,d}$s is
supported over determinants $\sL$ with $h^1(\sL)>1$;
see Example 6.1 of \cite{os16}. Thus, if $h^1(\sL)=1$
and $d=g-2$, there are no stable $\fg^k_{2,d}$s with determinant $\sL$ for
any $k \geq 2$. Accordingly, the limitations of our results in the case that 
$d$ is close to $g-2$ are due in part to actual failure of existence in this
range. There is no evidence against the possibility that, by studying more 
complicated families of limit linear series, the machinery of Theorem 
\ref{thm:main} could eventually produce sharp existence statements 
for the case that $h^1(\sL)=1$.
\end{rem}

We next discuss cases of Theorem \ref{thm:main-2} in which $\rho_{\sL}$ is
negative.

\begin{rem}\label{rem:neg-dim}
In Theorem \ref{thm:main-2}, the cases 
$(g,d,k)=(1,0,2), (3,2,2) \text{ or } (4,6,4)$
have $\rho_{\sL}<0$; in fact, it is elementary that these are precisely the 
cases for which the hypothesized inequality is satisfied, but $\rho_{\sL}<0$.
This does not lead to a contradiction, because in these cases we are only 
claiming to produce semistable limit linear series, which by virtue of larger 
automorphism groups may belong to a negative-dimensional family.
\end{rem}

Finally, for the sake of broader context we briefly discuss an observation 
of Gregorczyk and Newstead.

\begin{rem}\label{rem:nonexistence} Even restricting to the case
$h^1(\sL)=1$, the picture afforded by our modified expected dimension
$\rho_{\sL}$ is incomplete, in the sense that we can have 
the stable locus of $G^k_{2,\sL}(X)$ nonempty for $X$ general even when 
$\rho_{\sL}<0$. Indeed, Farkas and Ortega studied the case
$g=2a+1$, $d=2a+4$, and $k=4$ in \cite{f-o1}, and found that 
the locus of stable $\fg^k_{2,d}$s is supported over determinants $\sL$
with $h^1(\sL)=1$. However, Grzegorczyk and Newstead observed in 
Example 4.2 of \cite{g-n1} that in this case $\rho_{\sL}=7-g$, which is 
negative for $g>7$. Farkas and Ortega actually show that the 
$\fg^k_{2,d}$s are supported over a smaller locus of determinants satisfying 
an additional Koszul condition,
demonstrating that in order to understand the behavior of the fixed 
determinant case, it is necessary to consider invariants subtler than $h^1$.
Interestingly, in the case studied by Farkas and Ortega, the total
dimension (for varying determinant) still works out as predicted by our
$\rho_{\sL}$; see Example 5.2 of \cite{os19}.
\end{rem}

\appendix

\section{Examples}\label{app:examples}

As a companion to the rather complicated descriptions of the families of
limit linear series described in \S \ref{sec:families}, we give small 
examples for each of the four cases considered. On each component, we
list the vanishing sequences at $P_i$ and $Q_i$, as well as the underlying
vector bundles. We use $\sO(a,b)$ as shorthand for $\sO_{X_i}(aP_i+bQ_i)$,
and if the determinant of $\sE^i$ is fixed to be $\sL_i$, we similarly
write $\sL_i(a,b)$ for $\sL_i(-aP_i-bQ_i)$. We also use $\eind(a,b)$ to 
denote the unique (even-degree) indecomposable vector bundle of rank $2$ 
and degree $2a+2b$ containing $\sO(a,b)$ as a line subbundle, and 
$\oind(\sL_i)$ to denote an (odd-degree) indecomposable vector bundle of 
rank $2$ and determinant $\sL_i$.

All cases follow the same pattern on the first $(k')^2$ components,
and on the final components when $g$ is strictly greater than required.
In addition, when $d<2g-3$, each case involves an interval of cycling 
through all $k$ sections in succession. In order to keep examples at
a manageable size, we keep $g$ minimal in each case, and we only have
$d<2g-3$ in the first example.

\begin{subex}\label{ex:example-1}
Consider the case $g=8$, $k=4$, and $d=12$, treated by
Proposition \ref{prop:even-even}. In this case, the determinant is
uniquely determined on the first six components, and generic (but fixed)
on the remaining ones. However, the special form of the determinant will
not be relevant on the fifth and sixth components. The constructed family 
of limit linear series has the following form.

\begin{center}
\begin{tabular}{lr|lr|lr|lr|lr|lr|lr|lr}
0 & 6 & 0 & 6 & 0 & 5 & 1 & 4 & 2 & 4 & 2 & 3 & 3 & 2 & 4 & 1 \\
0 & 6 & 0 & 5 & 1 & 5 & 1 & 4 & 2 & 3 & 3 & 3 & 3 & 2 & 4 & 1 \\
1 & 4 & 2 & 4 & 2 & 3 & 3 & 3 & 3 & 2 & 4 & 1 & 5 & 1 & 5 & 0 \\
1 & 4 & 2 & 3 & 3 & 3 & 3 & 3 & 3 & 2 & 4 & 1 & 5 & 0 & 6 & 0 \\
\hline\\[-.1in]
\multicolumn{2}{c}{$\begin{matrix}\sO(0,6)\\\oplus\\\sO(0,6)\end{matrix}$} &
\multicolumn{2}{c}{$\begin{matrix}\sO(0,6)\\\oplus\\\sO(2,4)\end{matrix}$} &
\multicolumn{2}{c}{$\begin{matrix}\sO(1,5)\\\oplus\\\sO(3,3)\end{matrix}$} &
\multicolumn{2}{c}{$\begin{matrix}\sO(3,3)\\\oplus\\\sO(3,3)\end{matrix}$} &
\multicolumn{2}{c}{$\begin{matrix}\sO(2,4)\\\oplus\\\sO(6,0)\end{matrix}$}
&
\multicolumn{2}{c}{$\begin{matrix}\sO(3,3)\\\oplus\\\sO(7,-1)\end{matrix}$}
&
\multicolumn{2}{c}{$\begin{matrix}\sO(5,1)\\\oplus\\\sL_7(5,1)\end{matrix}$}
&
\multicolumn{2}{c}{$\begin{matrix}\sO(6,0)\\\oplus\\\sL_8(6,0)\end{matrix}$}
\end{tabular}
\end{center}
\end{subex}

\begin{subex}\label{ex:example-2}
Consider the case $g=7$, $k=5$, and $d=12$, treated by
Proposition \ref{prop:even-odd}. Now, the determinant is uniquely
determined on all seven components, and the constructed family of limit
linear series has the following form.

\begin{center}
\begin{tabular}{lr|lr|lr|lr|lr|lr|lr}
0 & 6 & 0 & 6 & 0 & 5 & 1 & 4 & 2 & 4 & 2 & 3 & 3 & 2 \\
0 & 6 & 0 & 5 & 1 & 5 & 1 & 4 & 2 & 3 & 3 & 2 & 4 & 1 \\
1 & 4 & 2 & 4 & 2 & 3 & 3 & 3 & 3 & 2 & 4 & 2 & 4 & 1 \\
1 & 4 & 2 & 3 & 3 & 3 & 3 & 3 & 3 & 2 & 4 & 1 & 5 & 0 \\
2 & 3 & 3 & 2 & 4 & 1 & 5 & 0 & 6 & 0 & 6 & 0 & 6 & 0 \\
\hline\\[-.1in]
\multicolumn{2}{c}{$\begin{matrix}\sO(0,6)\\\oplus\\\sO(0,6)\end{matrix}$} &
\multicolumn{2}{c}{$\begin{matrix}\sO(0,6)\\\oplus\\\sO(2,4)\end{matrix}$} &
\multicolumn{2}{c}{$\begin{matrix}\sO(1,5)\\\oplus\\\sO(3,3)\end{matrix}$} &
\multicolumn{2}{c}{$\begin{matrix}\sO(3,3)\\\oplus\\\sO(3,3)\end{matrix}$} &
\multicolumn{2}{c}{$\begin{matrix}\sO(2,4)\\\oplus\\\sO(6,0)\end{matrix}$}
&
\multicolumn{2}{c}{$\begin{matrix}\sO(4,2)\\\oplus\\\sO(6,0)\end{matrix}$}
&
\multicolumn{2}{c}{$\eind(6,0)$}
\end{tabular}
\end{center}
\end{subex}

\begin{subex}\label{ex:example-3}
Consider the case $g=7$, $k=4$, and $d=11$, treated by
Proposition \ref{prop:odd-even}. Now, the determinant is uniquely
determined on the first six components, and the constructed family of limit
linear series has the following form.

\begin{center}
\begin{tabular}{lr|lr|lr|lr|lr|lr|lr}
0 & 5 & 0 & 5 & 0 & 4 & 1 & 3 & 2 & 3 & 2 & 2 & 3 & 1 \\
0 & 5 & 0 & 4 & 1 & 4 & 1 & 3 & 2 & 2 & 3 & 2 & 3 & 1 \\
1 & 3 & 2 & 3 & 2 & 2 & 3 & 2 & 3 & 2 & 3 & 1 & 4 & 0 \\
1 & 3 & 2 & 2 & 3 & 2 & 3 & 2 & 3 & 1 & 4 & 0 & 5 & 0 \\
\hline\\[-.1in]
\multicolumn{2}{c}{$\begin{matrix}\sO(0,5)\\\oplus\\\sO(0,5)\end{matrix}$} &
\multicolumn{2}{c}{$\begin{matrix}\sO(0,5)\\\oplus\\\sO(2,3)\end{matrix}$} &
\multicolumn{2}{c}{$\begin{matrix}\sO(1,4)\\\oplus\\\sO(3,2)\end{matrix}$} &
\multicolumn{2}{c}{$\begin{matrix}\sO(3,2)\\\oplus\\\sO(3,2)\end{matrix}$} &
\multicolumn{2}{c}{$\oind(\sO(8,3))$}
&
\multicolumn{2}{c}{$\begin{matrix}\sO(3,2)\\\oplus\\\sO(6,-1)\end{matrix}$}
&
\multicolumn{2}{c}{$\begin{matrix}\sO(5,0)\\\oplus\\\sL_7(5,0)\end{matrix}$}
\end{tabular}
\end{center}

Note that in this case, the odd-degree bundle is on the fifth component.
\end{subex}

\begin{subex}\label{ex:example-4}
Consider the case $g=10$, $k=5$, and $d=17$, treated by
Proposition \ref{prop:odd-odd}. Now, the determinant is uniquely
determined on the first nine components, and the constructed family of limit
linear series has the following form.

\begin{center}
\resizebox{\linewidth}{!}{%
\begin{tabular}{lr|lr|lr|lr|lr|lr|lr|lr|lr|lr}
0 & 8 & 0 & 8 & 0 & 7 & 1 & 6 & 2 & 6 & 2 & 5 & 3 & 4 & 4 & 3 & 5 & 2 & 6 & 2 \\
0 & 8 & 0 & 7 & 1 & 7 & 1 & 6 & 2 & 5 & 3 & 4 & 4 & 4 & 4 & 3 & 5 & 2 & 6 & 1 \\
1 & 6 & 2 & 6 & 2 & 5 & 3 & 5 & 3 & 4 & 4 & 4 & 4 & 3 & 5 & 2 & 6 & 1 & 7 & 1 \\
1 & 6 & 2 & 5 & 3 & 5 & 3 & 5 & 3 & 4 & 4 & 3 & 5 & 2 & 6 & 2 & 6 & 1 & 7 & 0 \\
2 & 5 & 3 & 4 & 4 & 3 & 5 & 2 & 6 & 2 & 6 & 2 & 6 & 1 & 7 & 0 & 8 & 0 & 8 & 0 \\
\hline\\[-.1in]
\multicolumn{2}{c}{$\begin{matrix}\sO(0,8)\\\oplus\\\sO(0,8)\end{matrix}$} &
\multicolumn{2}{c}{$\begin{matrix}\sO(0,8)\\\oplus\\\sO(2,6)\end{matrix}$} &
\multicolumn{2}{c}{$\begin{matrix}\sO(1,7)\\\oplus\\\sO(3,5)\end{matrix}$} &
\multicolumn{2}{c}{$\begin{matrix}\sO(3,5)\\\oplus\\\sO(3,5)\end{matrix}$} &
\multicolumn{2}{c}{$\begin{matrix}\sO(2,6)\\\oplus\\\sO(6,2)\end{matrix}$} &
\multicolumn{2}{c}{$\begin{matrix}\sO(4,4)\\\oplus\\\sO(6,2)\end{matrix}$} &
\multicolumn{2}{c}{$\begin{matrix}\sO(4,4)\\\oplus\\\sO(8,0)\end{matrix}$} &
\multicolumn{2}{c}{$\begin{matrix}\sO(6,2)\\\oplus\\\sO(8,0)\end{matrix}$} &
\multicolumn{2}{c}{$\eind(8,0)$}
&
\multicolumn{2}{c}{$\oind(\sL_{10})$}
\end{tabular}%
}
\end{center}

Note that $\sL_{10}$ has odd degree $17$ in this case.
\end{subex}

\section*{Index of notation and terminology}

\begin{minipage}[t]{.55\textwidth}
\noindent $G$, Situation \ref{sit:basic}

\noindent $H_d$, Situation \ref{sit:basic}

\noindent $v(\e)$, Situation \ref{sit:basic}

\noindent $v(P)$, Definition \ref{defn:basic} 

\noindent $d(w,w')$, Definition \ref{def:distance}

\noindent $f_P$, Definition \ref{defn:basic} 

\noindent $f_{w,w'}$, Notation \ref{notn:fww}

\noindent $f_{\e}$, Notation \ref{not:mors}

\noindent $\left<,\right>^{\vp_{\bullet}}_{\bullet}$, Definition 
\ref{defn:tangent-form} 

\noindent $\sO_v$, Notation \ref{not:twisting-bundles-ii} and
Notation \ref{not:twisting-bundles-ii-2}

\noindent $\sO_{w,w'}$, Notation \ref{not:twisting-bundles-ii} and
Notation \ref{not:twisting-bundles-ii-2}
\end{minipage}
\begin{minipage}[t]{.45\textwidth}
\noindent $\sO'_{w,w'}$, Notation \ref{not:twist-rk-1}

\noindent $\sE_w$, Notation \ref{not:twisting-bundles-ii-2}

\noindent $\cG^{k,\EHT}_{r,d_{\bullet}}(X_0)$, Notation \ref{not:gkrd-eht}

\noindent $\cG^{k,\EHT}_{r,\sL,d_{\bullet}}(X_0)$, Definition 
\ref{def:eht-det}

\noindent $\cG^{k}_{r,w_0,d_{\bullet}}(X/B)$, Definition \ref{def:grd-space}

\noindent $\cG^{k}_{r,\sL,d_{\bullet}}(X/B)$, Definition \ref{def:grd-space}

\noindent $\cM_{r,\sL}(X_0)$, Definition \ref{def:eht-det}

\noindent $\cP^k_{r,d_{\bullet}}(X_0)$, Notation \ref{notn:pkrd}

\noindent $V(-D)$, Notation \ref{not:twist-down}

\noindent $\ell(V,a_i,P)$, Notation \ref{notn:fiber-lines}
\end{minipage}
\vspace{10pt}

\noindent $(P,Q)$-adaptable $\fg^k_{r,d}$, Definition \ref{def:adapted}

\noindent $(P,Q)$-adapted basis, Definition \ref{def:adapted}

\noindent chain adaptable (limit linear series), Definition 
\ref{def:chain-adapted}

\noindent Eisenbud-Harris-Teixidor limit linear series, Definition 
\ref{def:eht-lls}

\noindent Eisenbud-Harris-Teixidor limit linear series of determinant 
$\sL$, Definition \ref{def:eht-det}

\noindent index (of a linked bilinear form), Definition \ref{defn:link-bilin}

\noindent internally simple ($s$-prelinked bundles), Definition 
\ref{defn:internal-simple}

\noindent internally simple point, Definition \ref{def:simple}

\noindent isotropic (linked subbundle), Definition \ref{defn:link-sub}

\noindent $\ell$-(semi)stable (vector bundle), Definition \ref{def:l-stable}

\noindent $\ell$-(semi)stable (limit linear series), Definition 
\ref{def:l-stable-lls}

\noindent linked alternating form, Definition \ref{defn:link-alt}

\noindent linked bilinear form, Definition \ref{defn:link-bilin}

\noindent linked internally symplectic form, Definition \ref{defn:link-symp}

\noindent linked linear series, Definition \ref{def:grd-space}

\noindent linked linear series of determinant $\sL$, Definition 
\ref{def:grd-space}

\noindent linked subbundle, Definition \ref{defn:link-sub}

\noindent locus of isotropy, Definition \ref{defn:link-isotropy-locus}

\noindent nondegenerate (line in a fiber), Definition \ref{def:special-lines}

\noindent $s$-prelinked (vector bundles), Definition \ref{defn:basic}

\noindent prelinked alternating Grassmannian, Definition \ref{def:LAG}

\noindent prelinked Grassmannian, Definition \ref{def:lg}

\noindent refined (Eisenbud-Harris-Teixidor limit linear series), Definition 
\ref{def:eht-lls}

\noindent simple ($s$-prelinked bundles), Definition \ref{defn:simple-pt}

\noindent simple point, Definition \ref{def:simple}

\noindent smoothing family, Definition \ref{def:smoothing-fam}

\noindent special (line bundle on a reducible curve), Definition 
\ref{def:special}

\bibliographystyle{amsalpha}
\bibliography{gen}

\newcommand{\noopsort}[1]{} \newcommand{\printfirst}[2]{#1}
  \newcommand{\singleletter}[1]{#1} \newcommand{\switchargs}[2]{#2#1}
\providecommand{\bysame}{\leavevmode\hbox to3em{\hrulefill}\thinspace}
\providecommand{\MR}{\relax\ifhmode\unskip\space\fi MR }
\providecommand{\MRhref}[2]{%
  \href{http://www.ams.org/mathscinet-getitem?mr=#1}{#2}
}
\providecommand{\href}[2]{#2}
\begin{thebibliography}{{Tei}91b}

\bibitem[BBPN]{b-b-n2}
Usha Bhosle, Leticia Brambila-Paz, and Peter Newstead, \emph{On linear systems
  and a conjecture of {D}.\ {C}.\ {B}utler}, preprint.

\bibitem[BF98]{b-f2}
Aaron Bertram and Burt Feinberg, \emph{On stable rank two bundles with
  canonical determinant and many sections}, Algebraic Geometry (Catania
  1993/Barcelona 1994), Lecture Notes in Pure and Applied Mathematics, vol.
  200, Dekker, 1998, pp.~259--269.

\bibitem[FL81]{f-l1}
William Fulton and Robert Lazarsfeld, \emph{On the connectedness of degeneracy
  loci and special divisors}, Acta Mathematica \textbf{146} (1981), 271--283.

\bibitem[FO11]{f-o1}
Gavril Farkas and Angela Ortega, \emph{The maximal rank conjecture and rank two
  {B}rill-{N}oether theory}, Pure and Applied Mathematics Quarterly \textbf{7}
  (2011), no.~4, 1265--1295, Special issue in honor of Eckart Viehweg.

\bibitem[Gie82]{gi1}
D.~Gieseker, \emph{Stable curves and special divisors: {P}etri's conjecture},
  Inventiones Mathematicae \textbf{66} (1982), 251--275.

\bibitem[GN]{g-n1}
Ivona Grzegorczyk and Peter Newstead, \emph{On coherent systems with fixed
  determinant}, preprint.

\bibitem[GT09]{g-t1}
Ivona Grzegorczyk and Montserrat {Teixidor i Bigas}, \emph{Brill-{N}oether
  theory for stable vector bundles}, Moduli Spaces and Vector Bundles (Leticia
  Brambila-Paz, Steven Bradlow, Oscar Garc\'ia-Prada, and S.\ Ramanan, eds.),
  London Mathematical Society Lectures Notes Series, vol. 359, Cambridge
  University Press, 2009, pp.~29--50.

\bibitem[LNP]{l-n-p1}
Herbert Lange, Peter Newstead, and Seong~Suk Park, \emph{Non-emptiness of
  {B}rill-{N}oether loci in $m(2,k)$}, preprint.

\bibitem[Muk95]{mu2}
Shigeru Mukai, \emph{Vector bundles and {B}rill-{N}oether theory}, Current
  topics in complex algebraic geometry, MSRI Publications, vol.~28, Cambridge
  University Press, 1995, pp.~145--158.

\bibitem[Muk01]{mu6}
\bysame, \emph{Non-abelian {B}rill-{N}oether thoery and {F}ano $3$-folds},
  Sugaku Expositions \textbf{14} (2001), no.~2, 125--153.

\bibitem[Oss06]{os8}
Brian Osserman, \emph{A limit linear series moduli scheme}, Annales de
  l'Institut Fourier \textbf{56} (2006), no.~4, 1165--1205.

\bibitem[Oss11]{os17}
\bysame, \emph{Linked {H}om spaces}, Mathematical Research Letters \textbf{18}
  (2011), no.~2, 329--335.

\bibitem[Oss13a]{os16}
\bysame, \emph{Brill-{N}oether loci in rank $2$ with fixed determinant},
  International Journal of Mathematics \textbf{24} (2013), no.~1350099, 24
  pages.

\bibitem[Oss13b]{os21}
\bysame, \emph{Relative dimension of morphisms and dimension for algebraic
  stacks}, preprint, 2013.

\bibitem[Oss13c]{os19}
\bysame, \emph{Special determinants in higher-rank {B}rill-{N}oether theory},
  International Journal of Mathematics \textbf{24} (2013), no.~1350084, 20
  pages.

\bibitem[Oss14a]{os20}
\bysame, \emph{Limit linear series moduli stacks in higher rank}, preprint,
  2014.

\bibitem[Oss14b]{os22}
\bysame, \emph{Stability of vector bundles on curves and degenerations},
  preprint, 2014.

\bibitem[OT14]{o-t1}
Brian Osserman and Montserrat {Teixidor i Bigas}, \emph{Linked alternating
  forms and linked symplectic {G}rassmannians}, International Mathematics
  Research Notices \textbf{2014} (2014), no.~3, 720--744.

\bibitem[{Tei}91a]{te1}
Montserrat {Teixidor i Bigas}, \emph{Brill-{N}oether theory for stable vector
  bundles}, Duke Mathematical Journal \textbf{62} (1991), no.~2, 385--400.

\bibitem[{Tei}91b]{te4}
\bysame, \emph{Brill-{N}oether theory for vector bundles of rank 2}, Tohoku
  Mathematical Journal \textbf{43} (1991), no.~1, 123--126.

\bibitem[{Tei}04]{te5}
\bysame, \emph{Rank 2 vector bundles with canonical determinant}, Mathematische
  Nachrichten \textbf{265} (2004), 100--106.

\bibitem[Zha14]{zh2}
Naizhen Zhang, \emph{Towards the {B}ertram-{F}einberg-{M}ukai conjecture},
  preprint, 2014.

\end{thebibliography}

\end{document}